\documentclass[10pt,twoside]{amsart}
 \usepackage{amsmath,amssymb,amsthm,amscd,latexsym,graphicx}
 \usepackage{epigraph}
 \usepackage[OT2,OT1]{fontenc}

\usepackage{graphicx}
\usepackage[all]{xy}

 \font\caps=cmcsc10                    
 \font\Caps=cmcsc10 scaled \magstep1   

 \makeatletter
 \@specialpagefalse\headheight=8.5pt
 \makeatother

 \def\TSkip{\medskip}
 \newbox\TheTitle{\obeylines\gdef\GetTitle #1
 \ShortTitle  #2
 \SubTitle    #3
 \Author      #4
 \ShortAuthor #5
 \EndTitle
 {\setbox\TheTitle=\vbox{\baselineskip=20pt\let\par=\cr\obeylines%
 \halign{\centerline{\Caps##}\cr\noalign{\medskip}\cr#1\cr}}%
         \copy\TheTitle\TSkip\TSkip%
 \def\next{#2}\ifx\next\empty\gdef\STitle{#1}\else\gdef\STitle{#2}\fi%
 \def\next{#3}\ifx\next\empty%

 \else\setbox\TheTitle=\vbox{\baselineskip=20pt\let\par=\cr\obeylines%
     \halign{\centerline{\caps##} #3\cr}}\copy\TheTitle\TSkip\TSkip\fi%
 \centerline{\caps #4}\TSkip\TSkip%
 \def\next{#5}\ifx\next\empty\gdef\SAuthor{#4}\else\gdef\SAuthor{#5}\fi%
 \catcode'015=5}}\def\Title{\obeylines\GetTitle}

 \def\Abstract{\begingroup\narrower
     \parskip=\medskipamount\parindent=0pt{\caps Abstract. }}

 \long\def\MSC#1\EndMSC{\def\arg{#1}\ifx\arg\empty\relax\else
      {\par\narrower\noindent%
      2000 Mathematics Subject Classification: #1\par}\fi}

 \long\def\KEY#1\EndKEY{\def\arg{#1}\ifx\arg\empty\relax\else
         {\par\narrower\noindent Keywords and Phrases: #1\par}\fi\TSkip}

 \long\def\DATE#1\EndDATE{\def\arg{#1}\ifx\arg\empty\relax\else
         {\par\narrower\noindent \center{\textit{#1}}\par}\fi\TSkip\TSkip\TSkip}

 \font\bf= cmbx10 at 10pt

 \newcommand{\Sym}{\hbox{\ptitgot S}}

 \newcommand{\lra}{\longrightarrow}

 \newcommand{\A}{\mathbb A}

 \newcommand{\T}{\mathcal T} 
  
 \newcommand{\C}{\mathbb{C}}
 \newcommand{\F}{\mathbb{F}}
 \newcommand{\G}{\mathbb{G}}

 \newcommand{\N}{\mathbb{N}}
 \renewcommand{\P}{\mathbb P}
 \newcommand{\Q}{\mathbb{Q}}
 \newcommand{\R}{\mathbb{R}}
 \newcommand{\Z}{\mathbb{Z}}

 \newcommand{\Ind}{\mathrm{Ind}}

 \newcommand{\Cor}{\mathrm{Cor}}
 \renewcommand{\i}{\iota}

   \newcommand{\Res}{\mathrm{Res}}

 \newcommand{\GL}{\mathrm{GL}}

 \newcommand{\Coker}{\mathrm{Coker}}

  \newcommand{\Det}{\mathrm{Det}}
 
  \newcommand{\ev}{\mathrm{ev}}

  \newcommand{\End}{\mathrm{End}}
 \newcommand{\Ext}{\mathrm{Ext}}
\newcommand{\Gal}{\mathrm{Gal}}
 
 \newcommand{\Hom}{\mathrm{Hom}}
 
  \newcommand{\Ver}{\mathrm{Ver}}
   \newcommand{\Frob}{\mathrm{Frob}}

 \newcommand{\M}{\mathcal{M}}

\newcommand{\Id}{\mathrm{Id}}

 \newcommand{\KW}{\mathcal{KW}}
 
\renewcommand{\Im}{\mathrm{Im}}
 \newcommand{\tr}{\mathrm{tr}}

\newcommand{\Ker}{\mathrm{Ker}}
 
 \newcommand{\Spec}{\mathrm{Spec}}

 \newcommand{\Gr}{\mathrm{Gr}}

\renewcommand{\S}{\mathcal{S}}
\renewcommand{\M}{\mathcal{M}}
\newcommand{\W}{\mathbf{W}}
\newcommand{\E}{\mathcal{E}}
\newcommand{\frob}{\mathrm{frob}}

\newcommand{\K}{\mathcal{K}}

\renewcommand{\Sym}{\mathrm{Sym}}

 \theoremstyle{plain}
 \newtheorem{thm}{Theorem}[section]
 \newtheorem{defi}[thm]{Definition}

 \newtheorem{prop}[thm]{Proposition}

 \newtheorem{lem}[thm]{Lemma}
 \newtheorem{coro}[thm]{Corollary}
\newtheorem{conj}[thm]{Conjecture.}
 \theoremstyle{remark}
 \newtheorem{rem}[thm]{Remark}
 
  \newtheorem{pb}[thm]{Problem}
 
 \newtheorem{ex}[thm]{Example}
 
  \newtheorem{exo}[thm]{Exercise}
 \newenvironment{dem}{{\bf Proof.}}{\hfill$\square$}

\date{March 2010}

\begin{document}


 \Title
Lifting Theorems and Smooth Profinite Groups 
 \ShortTitle
Smooth
\SubTitle
 \Author
 Charles De Clercq and Mathieu Florence\footnote{Florence was partially supported by the French National Agency 
(Project GeoLie ANR-15-CE40-0012).}
 \ShortAuthor
 \EndTitle

\DATE
October 2017
\EndDATE

\address{ Charles De Clercq, Equipe Topologie Alg\'ebrique, Laboratoire Analyse, G\'eom\'etrie et Applications, Universit\'e Paris 13, 93430 Villetaneuse.}
\address{Mathieu Florence, Equipe de Topologie et G\'eom\'etrie Alg\'ebriques, Institut de Math\'ematiques de Jussieu,  Universit\'e Pierre et Marie Curie, 4, place Jussieu, 75005 Paris. }

\epigraph{La simplicit\'e est la r\'eussite absolue. Apr\`es avoir jou\'e une grande quantit\'e de notes, toujours plus de notes, c'est la simplicit\'e qui \'emerge, comme une r\'ecompense venant couronner l'art.}{Fr\'ed\'eric Chopin.}


\Abstract
This work is motivated by the search for an "explicit" proof of the Bloch-Kato conjecture in Galois cohomology, proved by Voevodsky in \cite{Vo}. Our concern here is to lay the foundation for a theory that, we believe, will lead to such a proof- and to further applications. \\
Let $p$ be a prime number. Let $k$ be  a perfect field of characteristic $p$. Let $m$ be a positive integer. Our first goal is to provide a canonical process for "lifting" a module $M$, over the ring of Witt vectors $\W_m(k)$ (of length $m$), to a $\W_{m+1}(k)$-module, in a way that deeply respects Pontryagin duality. These are our big, medium and small Omega powers (cf. Definitions  \ref{mediumbig} and \ref{small}), each of which naturally occurs as a direct factor of the previous one. In the case where $M$ is a $k$-vector space, they come  equipped with Verschiebung and Frobenius operations (cf. section \ref{FrobVerOmega}). If moreover the field $k$ is finite, Omega powers are endowed with a striking extra operation: the Transfer, to shifted Omega powers of finite-codimensional linear subspaces (cf. section \ref{sectranssmall}). To show how this formalism fits into Galois theory, we  first offer an axiomatized approach to  Hilbert's Theorem $90$ (or  more precisely, to its consequence for cohomology with finite coefficients: Kummer theory).  In the context of profinite group cohomology, we thus define the notions  a cyclotomic $G$-module (Definition \ref{defismoothmodule}), and of a smooth profinite group  (Definition \ref{defismooth}). We bear in mind that the fundamental example is that of an absolute Galois group, together with the Tate module of roots of unity. We then define the notion of exact sequences of $G$-modules \textit{of Kummer type}- see Definition \ref{Kummertype}. \\
To finish, we give applications of this formalism. The first ones are the Stable Lifting Theorems (Theorems \ref{StableLift} and \ref{Lifting2}), enabling the lifting to higher torsion  in the cohomology of smooth profinite groups, with $p$-primary coefficients. To illustrate their meaning, we translate the second one in the concrete context of Galois cohomology, with values in two-dimensional Galois representations (Corollary \ref{concrete}). We finish by an application to $p$-adic deformations. We state and prove a (perhaps unusual)   general descent statement, for the quotient map $\Z/ p^2 \Z \lra \Z /p \Z$. It is Proposition \ref{wrinklelift}.

\newpage

\tableofcontents

\section{Introduction.}

The main goal of this paper is to define and study a very low level smoothness notion. By "low level", we mean that most of our constructions just depend  on a given prime number $p$. In a subsequent paper, we plan to strenghten  the Stable Lifting Theorems for Galois cohomology that we present here, in order to achieve an explicit proof of the Bloch-Kato conjecture. This was our starting motivation.  Let us simply say that, tough we did not achieve a full proof yet, we are confident that the present approach will eventually be successful- when the foundations are solid enough. \\

It is  folklore  that mod $p$ cohomology often encodes subtle phenomena in higher $p$-primary torsion. For instance,  section \ref{alternategamma} presents a mod $p$ Hochschild $2$-cocycle, canonically attached to a mod $p^2$ algebraic object.  The general philosophy of this paper is, whenever possible,   to lift mod $p$ algebraic data to higher $p$-primary torsion in a natural way.  This data may be a module, an algebra, a Hopf algebra,  a cohomology class, etc...  We try to do so using, as a basic algebraic tool, divided powers of  finite modules over the ring of Witt vectors $\W(k)$.  It is often enlightening to think of these modules as the coefficients of some cohomology theory. Hence, in practice, the field $k$  will often be finite- a seemingly crucial requirement to define the Transfers (see Section \ref{sectrans}).

  We now wish to share some ideas which motivated this work. In what follows,  $k$ is a perfect field of characteristic $p$, and $F$ is a field  of characteristic not $p$.  Let $F_{sep}/F$ be a separable closure of $F$. We put $G:=\Gal(F_{sep}/F )$.  We denote by $d=p^s$  a power of the prime $p$.
\subsection{First idea.}

Let us start by giving a purely Galois-theoretic statement of the Norm-Residue Isomorphism Theorem, proved by Rost, Suslin and Voevodsky. Our point here is to make clear that this Theorem presents a deep connection with the general notion of smoothness.\\
For each integer $i \geq 0,$ denote by  $\mathcal H^i$ the Galois cohomology group $H^i(F, \mu_p^{\otimes i})$. Then the  cup-product operation, being $\F_p$-multilinear,  yields a homomorphism of graded $\F_p$-algebras $$h: \bigoplus_{i=0}^\infty T^i( \mathcal H^1) \lra \bigoplus_{i=0}^\infty \mathcal H^i,  $$ where $  T^i( \mathcal H^1) = \bigotimes_{\F_p}^{i}\mathcal H^1 $ is the $i$-th fold tensor power of $\mathcal H^1 $ (note that this map actually factors through the exterior power $\Lambda^i_{\F_p}(\mathcal H^1)$, if $p$ is odd or if $p=2$ and $(-1)$ is a square in $F$). Then the Norm-Residue Isomorphism Theorem states that $h$ is surjective, and that its kernel is generated, in degree two, by  all pure tensors $$ a \otimes b \in \mathcal H^1 \otimes_{\F_p} \mathcal H^1, $$ such that $a \cup b =0 \in  \mathcal H^2$.

Indeed, pure tensors arising from Steinberg's relation are a particular instance of these, and it is not too hard to show directly that all pure quadratic tensors in $\Ker (h)$ are combinations of pure tensors arising from Steinberg's relation.\\

Now, let  $(A,M_x)$ be the local ring of  variety $X$ over $\F_p$, at a \textit{smooth} rational point $x \in X(\F_p)$.  Then the natural homomorphism of graded $\F_p$-algebras $$  \bigoplus_{i=0}^\infty T^i( M_x/M_x^2) \lra \bigoplus_{i=0}^\infty \mathcal (M_x^i/M_x^{i+1}),  $$  is surjective, with kernel generated by degree two tensors of the shape $a\otimes b -b \otimes a.$\\

It is clear that these two results present strong similarities- though the second one, of commutative nature, is much easier to prove...\\
What is more, it is known to expert that the hard part of the Bloch-Kato conjecture is to show that the natural map $$H^i(F,\mu_{p^s}^{\otimes i}) \lra H^i(F,\mu_{p}^{\otimes i}) $$ is surjective, for every $i$ and every field $F$, of characteristic not $p$. Thinking (abusively) of $$H^i(F,\mu_{p^s}^{\otimes i})$$ as points modulo $p^s$  of some  algebraic object defined over $\Z_p$,  the surjectivity in question is, again, the definition of (formal) smoothness. Note that, for $i=1$, surjectivity is given by usual Kummer theory. Our Smoothness Conjecture \ref{smoothconj} states, in particular, that surjectivity for  $i$ arbitrary should "formally" follow from the $i=1$ case.  \\

Let us now briefly explain how one can hope to apply general lifting results in Galois cohomology, to prove the Bloch-Kato conjecture. For simplicity, we concentrate on surjectivity part, in the $i=2$ case  (the Merkurjev-Suslin Theorem). \\

Let $e$ be a class in $ H^2(F,\mu_p^{\otimes 2})$. By a general fact from group cohomology, one can find a finite discrete $G$-module $V$, which is an $\F_p$-vector space, classes   $$ a \in H^1(F, V \otimes_{\F_p} \mu_p) ~~~\mathrm{and} ~~~ b \in H^1(F, V^\vee \otimes_{\F_p} \mu_p),$$  such that $$e =a \cup b,$$  where the cup-product is relative   to the canonical pairing $$(V \otimes_{\F_p} \mu_p )\times (V^\vee \otimes_{\F_p}  \mu_p) \lra \mu_{p}^{\otimes 2} .$$  (here $V^\vee=\Hom_{\F_p}(V,\F_p)$). \\This  may seem obscure at first. Assuming that $\mu_p= \F_p$ for simplicity -which is a harmless assumption- it  becomes clear if we adopt the viewpoint of (Yoneda) extensions. Indeed, a   $2$-extension in $\Ext^2_{(\F_p,G)}(\F_p, \F_p)$ ($=H^2(G,\F_p)$) can always be seen as the cup-product of two $1$-extensions.\\
If $V=\F_p^N$ was equipped with the   \textit{trivial} $G$-action, then $e$ would be a sum of $N^2$ symbols, and the job would be done. More generally, imagine we can  find a surjective morphism of $(\F_p,G)$-modules $$ f: W \lra V,$$ such that the following two conditions hold.

1) the induced map $H^1(G,W) \rightarrow H^1(G,V)$ is onto.\\
2) the $(\F_p,G)$-module $W$ is  permutation, i.e. has an $\F_p$-basis which is permuted by $G$.

Using Shapiro's Lemma, we see that the (projection) formula $$f_*(c) \cup b = c \cup (f^\vee)_*(b), $$ valid for any $c \in H^1(G,W),$ would then present $e$ as a sum of \textit{corestrictions} of symbols. An   input from Milnor $K$-theory (namely, the existence of the norm, and its compatibility with the norm-residue homomorphism), shows that a corestriction of a symbol (in Galois cohomology) is a sum of symbols, and we would be done. \\
A map $f$ with the properties above does, as one may guess, not exist in general.\\ Meanwhile, our Stable Lifting Theorem \ref{Lifting2} (for $k=\F_p$) does a very similar job. However, it applies only in higher $p$-primary torsion. Indeed, note that the module $$\bigoplus_{L \in \P(V)}  \underline\Omega^n(L)(s)$$ is not a fearsome beast at all: it is a honest module which is \textit{induced from dimension one} (in the sense of Definition \ref{indurank1}, with the $H_i$'s being the stabilizers of lines of $V$). For example, if $n=0$ and $G$ is a pro-$p$-group, it is a permutation $(\F_p,G)$-module in the sense of $2)$ above.  \\ Note that the Stable Lifting Theorem can nonetheless apply to lift our cohomology classes $a$ and $b$, but only after pushing them by a power of the Verschiebung $$ \Ver^n: V \lra \underline \Omega^n(V).$$

In Exercise \ref{Requi}, we explain why the Stable Lifting Theorems do not hold for $n=0$, using Manin's $R$-equivalence  (in the context of Galois cohomology).\\

\subsection{Second idea.} The theory of Witt vectors associates, to every perfect field $k$ of characteristic $p$, a discrete valuation ring $\W(k)$, whose basic properties we shall recall in the next section. We develop here the theory of divided powers for torsion modules over Witt vectors, whose  purpose is, somehow, to categorify  Witt's construction. To do so, we use the divided powers functors $\Gamma^{p^n}_{\W(k)},$ applied to torsion $\W(k)$-modules. We view them as \textit{representing polynomial laws} (cf. \cite{Ro}, or the nice and short paper \cite{Fe}). Note that truncated Witt vectors themselves, through a simple recursive process, can be defined just using $\Gamma^{p}_{\Z}$, see Proposition \ref{newwitt}.\\ We try to proceed as functorially as possible. We eventually offer three ways of lifting a $\W_m(k)$-module to a  $\W_{m+n}(k)$-module: the big, medium and small Omega powers (respectively, $\overline \Omega^n$, $ \Omega^n$ and $\underline \Omega^n$). The composition formula $\overline \Omega^{n+n'}= \overline \Omega^n \circ \overline \Omega^{n'}$
only holds for big Omega powers. However, we prove (cf. Propositions \ref{mediumbig} and \ref{smallmedium}) that medium (resp. small) Omega powers  canonically embed, in a duality-preserving way, in big (resp. medium) Omega powers. In the course of the proofs of these two Propositions, mysterious $p$-adic constants appear- see Definition \ref{strange1} and Lemma \ref{strange2}. We do not know  understand their meaning. Do they  have one?\\
Note that, in dimension one, small, medium and big Omega powers  all coincide, and  are indeed a "categorification" of the multiplicative Teichm\"uller section $$\tau: k^\times \lra \W(k)^\times.$$ In the recent preprints \cite{K} and \cite{K2}, related constructions are proposed, in a different language. It would be interesting to explore the connections.\\
When the field $k$ is finite, we  introduce a $\W(k)$-linear map "in the wrong direction": the Transfer, notably for small Omega powers. We believe that it presents connections to the (algebraic) Steenrod algebra, as defined in \cite{Sm}. Indeed, our formula for the hyperplane Transfer  (cf. Lemma \ref{transcod1}), is very close to the formula defining $P (\xi)(l)$, in Smith's paper.\\
In section \ref{axioH90}, we spend some time to axiomatize  Kummer theory- a consequence of Hilbert's Theorem 90 for (Galois) cohomology with values in roots of unity.  In the general context of profinite group cohomology, we define the notion of a smooth profinite group $G$, of a cyclotomic $G$-module and of a Kummer-type exact sequence. We state the Smoothness Conjecture  \ref{smoothconj}, implying the Bloch-Kato conjecture.\\
Combining our formalism with classical techniques from group cohomology (restriction, corestriction and Shapiro's Lemma), we are finally able to prove very general results for the cohomology of smooth profinite groups: the Stable Lifting Theorems (Theorem \ref{StableLift} and  Theorem \ref{Lifting2}). They are not yet sufficient to prove the Smoothness Conjecture, but we strongly believe that they will- after some improvement.\\

The formalism we develop here presents  potential for applications to  other topics. We now venture to list five of these. On this matter, we deeply welcome  comments, (constructive) criticism,  suggestions and collaborative work.   

1) The first one is $p$-adic deformation theory. The descent statement that we offer in Proposition \ref{wrinklelift} is most likely an explicit description of an abstract \textit{non linear} (degree $p$) descent statement for algebraic structures mod $p^n$. It would be interesting to identify it.\\
To illustrate what we mean by "explicit", we proceed with an analogy in the classical context. This analogy is for sure well known to many experts- though the explicit computations that make it precise are hard to find in the litterature.\\ Grothendieck's faithfully flat descent theory for Modules (say, for simplicity, in the affine case) has many concrete incarnations. For $G$-Galois algebras, it specializes to  Galois descent (this is Speiser's Lemma, cf. \cite{GS}, Lemma 2.3.8). In characteristic $p$, for purely inseparable field extensions  of height one,  it specializes to Cartier's descent  (\textit{loc. cit.}, Theorem 9.3.6). If $X=\Spec(A)$ is an affine smooth variety over any field $k$ of chacteristic $p$, then the Frobenius morphism $$\Frob:  A \otimes_{\frob} k \lra A,$$ $$x\otimes \lambda \mapsto \lambda x^p $$ is finite and flat, and   Cartier's Frobenius descent is,  again, an explicit geometric description of  faithfully flat descent.

2) The second one  concerns the question of lifting (mod $p$) Galois representations to mod $p^2$ Galois representations- and perhaps even to higher torsion. This important topic has already been  investigated  by many authors, notably in the context of local or global fields. We already have some (new) results for arbitrary fields, and plan to publish them in a dedicated paper.

3)  We believe that Omega powers could perhaps be used in modular representation theory- of finite groups, or of algebraic groups. For instance, if $V$ is a finite-dimensional $k$-vector space, the  quotient $$\Omega^n(V)/p$$ is a $k$-linear representation of the algebraic group $\GL_k(V)$ (in the sense of \cite{J}). We believe that it cannot, in general, be obtained as a subquotient of  a tensor power of $V$. In other words, divided powers over $\W(k)$ are  required in its construction- though, in the end, it is a mod $p$ object.

4) It is likely that our "gentle" machinery can help to say something about  resolution of singularities. A strong reason for this fact is the following.  The $\W(k)$-module $$\Gamma_{\W(k)}^p(\W_m(k))(\simeq \W_{m+1}(k))$$ is an elementary algebraic blowup of $\W_m(k)$, that lifts (the exponent of) its torsion by one. Forming Omega powers of $\W_m(k)$-algebras (cf. section \ref{Omegaalg}) is thus a way of performing a vast amount of these small blowups. This point of view is connected to the notion of  Rees algebra of a module, as investigated in the recent paper \cite{St}.

5) Last, but not least, let us remark that our approach here is purely local, at a given prime $p$. Once polished,  it could be fruitfully  globalized, considering all  functors $\Gamma^n_\Z$ at once...\\

The paper is organized as follows. We first recall some classical facts about profinite groups, representation theory and cohomology. We then explain, in section \ref{shapiro}, a categorical formulation of the induction process from open subgroups and of Shapiro's Lemma. Though elementary, it plays an important role in this paper, where most properties concerning a profinite group $G$ (eg. $n$-surjectivity) involve all open subgroups of $G$ at once. We then emphasize the importance of Pontryagin duality (in algebra). Though invisible, it is omnipresent in all cohomological theories: the injective Abelian group $\Q / \Z$ naturally occurs when building canonical injective resolutions of sheaves. After that, it is (unfortunately...)  often  disregarded or forgotten. In section \ref{tenseproduct},  we begin with recalling that Pontryagin duality does not commute to the tensor product- even in the  category of $(\Z/ p^n \Z)$-modules, for $n \geq 2$. Note that, if it did, the  topological issue of tensor completions would be much simpler. We make a short attempt to define the "Tense Product", a symmetric monoidal operation on the  category of $(\Z/ p^n \Z)$-modules, that commutes to Pontryagin duality.  It behaves well with the Omega power functors, that we define later on.

In section \ref{gamma}, we recall (mostly well-known) facts about divided powers. We  see them as representing homogeneous polynomial laws- as explained in \cite{Ro}. We mainly concentrate on the case of modules over Witt vectors. Along the way, we give a simple presentation of truncated Witt vectors themselves, as a quotient of a divided power module over $\Z$ (cf. Proposition \ref{newwitt}).  In section \ref{frob}, we introduce the Frobenius and Verschiebung operators, for divided powers. In section \ref{omega}, we  investigate the lack of commutation between Pontryagin duality and divided powers- reminiscent of Section \ref{tenseproduct}.   We introduce (big and medium) Omega powers,  as the "correct" quotients of divided powers, commuting to duality. We show that medium Omega powers occur as a direct factor of big Omega powers. We study their first functorial properties. In section \ref{sectrans}, we define the Transfer, a fundamental gadget to prove Lifting Theorems in cohomology, by induction on the dimension of the  coefficients.  In section \ref{small}, we introduce small Omega powers. We show that they are a direct factor of medium Omega powers. Small Omega powers enjoy rich functorial properties (notably through the Transfer) which we begin to investigate in Section \ref{sectranssmall}. See, in particular, Proposition \ref{smallomegafun}. We prove the Integral Formulas for the Frobenius and the Verschiebung. They are, perhaps, connected to motivic integration. In Section \ref{axioH90}, we present a possible  axiomatization of Kummer theory. We define the notions of cyclotomic module, of smooth profinite group and of Kummer-type extension. We bear in mind that the   fundamental example of a smooth profinite group is that of an absolute Galois group, equipped with the Tate module of roots of unity. Section \ref{Hilbert} is a short digression, to stress the importance of Hilbert's Theorem 90 in our approach- perhaps the purest of all descent statements. In section \ref{stablesec}, we prove two first applications of our formalism to Galois theory: the Stable Lifting Theorems. We present a concrete corollary of the second one (Corollary \ref{concrete}). \\
We conclude by an application of our point of view to deformations: Proposition \ref{wrinklelift}, which is a descent statement for the quotient map $\Z /p^2 \Z \lra \Z /p \Z$. \\

This paper contains numerous remarks and exercises,  the goal of which is to help the reader getting familiar with our approach- especially for those wishing to read  it "linearly". Note that, though we decided to treat the case of   an arbitrary perfect field $k$ of characteristic $p$  when possible, the case $k=\F_p$ is the essential one.

\section{Notation and basic facts.}
Throughout this paper,  $p$ is a prime number. For obvious historical reasons, we could have chosen to denote $p$ by $l$: the prime "$p$" here is the "$l$" of $l$-adic cohomology. A few months ago, we thus made an attempt to replace  $p$ by $l$ everywhere in the text. The resulting formulae were esthetically questionable (if not ugly), and we decided to go back to the previous notation... \\ 

For any integer $n$, we denote by $v_p(n)$ the $p$-adic valuation of $n$. We denote by $\S_n$ the symmetric group on $n$ letters.

If $M$ is an Abelian group and $n\geq 1$ is an integer, we denote by $M[n]$ the $n$-torsion of $M$.
Let $A$ be a ring. 
If $M$ is an $A$-module, we denote by $$M^*=\Hom(M,A)$$ the $A$-dual of $M$. We denote by $$\Sym_A(M)=\bigoplus_{i=0}^\infty \Sym^i_A(M)$$ the symmetric algebra of $M$.  We denote by $$\Lambda_A(M)=\bigoplus_{i=0}^\infty \Lambda^i_A(M)$$ the exterior algebra of $M$. We have pairings $$ \Lambda^i_A(M) \times \Lambda^{i}_A(M^*) \lra A,$$ $$(x_1 \wedge \ldots \wedge x_i, \phi_1 \wedge \ldots \wedge \phi_i) \mapsto \det (\phi_b(x_a))_{1 \leq a,b \leq i}, $$ and $$ \Lambda^i_A(M) \times \Lambda^{j}_A(M) \lra  \Lambda^{i+j}_A(M),$$ $$ (x,y) \mapsto x \wedge y.$$  These are perfect if $M$ is a finite locally free $A$-module of (constant) rank $d$, and $i+j=d$. In that case, we put $$ \Det(M):=\Lambda_A^d(M);$$ it is an invertible $A$-module.\\If the $A$-module $M$ is locally free of finite rank, we denote by $$\A_A(M):=\Spec(\Sym_A(M^*))$$ the affine space of $M$; it is an affine variety over $\Spec(A)$. On the level of the functor of points, we have $$\A_A(M)(B)=M \otimes_A B,$$ for every commutative $A$-algebra $B$.

Let $k$ be a field.
Let $V$ be finite-dimensional $k$-vector space. We denote by $\delta(V)$ the dimension of $V$. We denote by $\P_k(V)$ the projective space of $V$, consisting of lines $L \subset V$ (when needed, these shall be identified with  hyperplanes in $V^*$).  It can, of course, be viewed as a $k$-variety. However, in this work (where in most cases $k$ and $V$ will be finite), it will only be considered as a set. Note that,  if $V$ is a linear representation of a group $G$,   $\P_k(V)$ is naturally  endowed with an action of $G$.\\

\subsection{Witt vectors.}
If $k$ is a perfect field of characteristic $p >0$, we denote by $\W(k)$ the ring of Witt vectors built out from $k$. It is, up to isomorphism,  the unique  complete discrete valuation ring whose maximal ideal is generated by $p$, and with residue field $k$.  Its construction is functorial in $k$.  For any positive integer $n$, we denote by $$\W_n(k):=\W(k)/ p^n$$ the truncated Witt vectors of size $n$.\\ Note that a simple (and perhaps new) recursive formula, presenting $\W_{n+1}(k)$ as  a quotient of the $p$-th divided power of the $\Z$-module $\W_{n}(k)$, shall be given later on  (Proposition \ref{newwitt}). \\We  put $$K:=\mathrm{Frac}(\W(k)).$$ We shall often use the natural arrow $$ \W_n(k) \lra K/ \W(k),$$ $$1 \mapsto \frac 1 {p^n} $$ to identify the $p^n$-torsion in  $K/ \W(k)$ with $\W_n(k)$. However, one has to be careful in doing so- see for instance section \ref{tenseproduct}.\\ For any $\W(k)$-module $M$, we put $$ M^\vee := \Hom_{\W(k)}(M, K/\W(k));$$ it is the Pontryagin dual of $M$. Note that Pontryagin duality extends the linear duality of $k$-vector spaces to all $\W(k)$-modules. More precisely, if $M$ is seen as a $\W_n(k)$-module, one has a canonical isomorphism $$ M ^\vee \simeq \Hom_{\W_n(k)}(M, \W_n(k)) = M^*.$$ The Frobenius morphism $$ k \lra k,$$ $$ x \mapsto x^p $$ lifts to  a ring homomorphism $$\frob: \W(k) \lra \W(k).$$  For any $\W(k)$-module $M$,  and any integer $i\geq 0,$ we put $$M ^{( i)} :=M \otimes_{\W(k)} \W(k);$$ where the tensor product is taken with respect to  $\frob^i$.

\subsection{Profinite groups and cohomology.}

Let $G$ be a profinite group. By definition, a $G$-set is a set $X$, equipped with a continuous action of $G$ (i.e.  such that the stabilizer of every element of $X$ is open in $G$).\\
Let $M$ be a discrete $G$-module; that is, an Abelian group $M$, equipped with the structure of a $G$-set, for which the action of $G$ is $\Z$-linear. We then denote by $H^n(G,M)$ the cohomology groups, defined by Serre in \cite{Se}. At our disposal, we have the restriction maps $$\Res: H^n(G,M) \lra H^n(G',M),$$ for any closed subgroup $G' \subset G$, and the corestriction maps $$\Cor: H^n(G',M) \lra H^n(G,M),$$ for any open subgroup  $G' \subset G$.\\ If $G' \subset G$ is an open subgroup, of index $n$ in $G$, then $\Cor \circ \Res$ equals multiplication by $n$. 

\begin{rem}
	In the course of proving results involving a profinite group $G$, we shall often reduce to the case  where $G$ is pro-$p$-group, using the standard "restriction-corestriction" argument. More precisely, imagine that the discrete $G$-module $M$ is of $p$-primary torsion, and that we have to show that a class in $H^n(G,M)$ is zero. Then, it is enough to show that its restriction to $H^n(G_p,M)$ vanishes, where $G_p$ is a pro-$p$-Sylow of $G$.\\
	
\end{rem}

 \subsection{Categories of representations. }

Let $G$ be a profinite group. Let $k$ be a perfect field of characteristic $p$, often finite in our applications.

\begin{defi}
 A $( \W(k),G)$-module is a torsion $\W(k)$-module $M$ of finite-type, endowed with a continuous $\W(k)$-linear action of $G$ (i.e. factoring through a nontrivial open subgroup of $G$).   A $(k,G)$-module is a  $( \W(k),G)$-module which is a $k$-vector space. 
\end{defi}

\begin{rem}
Assume that $k$ is finite, and that $F_{sep}/F$ is a separable closure of a field $F$. Then a $(k,\Gal(F_{sep}/F))$-module is nothing but a  Galois representation over the field $k$.
\end{rem}

\begin{rem}
	if $G$ is a pro-$p$-group, we shall, in many places, use the following classical facts. \\
	(i) Every one-dimensional $(k,G)$-module is trivial, i.e. isomorphic to $k$, equipped with the trivial action of $G$.\\
	(ii) Let $V$ be a nonzero $(k,G)$-module. Then, it admits a one-dimensional sub-$(k,G)$-module. Equivalently, we have $V^G \neq \{0\}.$
	
\end{rem}

\begin{defi}
We denote by  $\M( \W(k),G)$ (resp. $\M(k,G)$) the category of  $( \W(k),G)$-modules (resp. of $(k,G)$-modules), with morphisms being $\W(k)$-linear maps respecting the action of $G$. These categories are  Abelian. They come equipped with a tensor product  $$\otimes=\otimes_{\W(k)}.$$  They are, moreover, equipped with a perfect duality $$ M \mapsto M^\vee = \Hom_{\W(k)}(M, K/\W(k)).$$

\end{defi}

Among $( \W(k),G)$-modules, the simplest are those who come from an action of $G$  on a finite set- the permutation modules. Let us give a precise Definition.
 
	 \begin{defi}\label{indurank1}
	 Let $M$ be a $( \W(k),G)$-module. It is said to be \textit{induced from rank one} if it is isomorphic  to a finite direct sum  $$\bigoplus_i \Ind_{H_i}^G(L_i),$$ where $H_i \subset G$ are open subgroups, and $L_i$ are $( \W(k),H_i)$-modules, which are  (free) $\W_{n_i}(k)$-module of rank one, for some positive integers $n_i$.\\
	 If moreover all $L_i$ 's are equipped with the trivial $H_i$-action, $M$ is said to be a permutation module.
	 
	 \end{defi}
	\begin{rem}
	    If $G$ is a $p$-group,  all one-dimensional $(k,G)$-modules are trivial. Hence, a $(k,G)$-module is induced from rank one if and only if  it is permutation. Through the usual 'restriction-corestriction' argument, for $G$ arbitrary, $(k,G)$-modules which are induced from rank one may often be assumed to be permutation.
	\end{rem}

	\section{On induction from subgroups, and Shapiro's Lemma.}\label{shapiro}
	
	Shapiro's Lemma is a fundamental basic tool in group cohomology. We now briefly explain how we can view it.
	
	\begin{defi}
		Let $G$ be a profinite group. Let $X$ be a finite $G$-set. A discrete $G$-module over $X$ is the data of  $$\mathcal M=(M_x, \phi_{g,x}),$$ consisting of an Abelian group $M_x$, for each $x \in X$, and of additive maps $$\phi_{g,x}: M_x \lra M_{gx},$$ for each $x \in X$ and $g \in G$,  subject to the following conditions.\\
		(i) For all $x \in X$, and all $m \in M_x$,   the map $$G \lra  \bigsqcup_{g \in G}   M_{gx},$$ $$g \mapsto \phi_{g,x}(m),$$ is continuous (=locally constant). \\
		ii) For all $x \in X,$ we have $$\phi_{e,x} = \Id.$$
		(iii) For all  $x \in X$ and $g,h \in G$, we have $$ \phi_{g,hx} \circ \phi_{h,x}=\phi_{gh,x}.$$

	\end{defi}
	
	\begin{rem}
		In the particular case of a one-element set, it is clear that a discrete $G$-module over $\{*\}$ is simply a discrete $G$-module.
	\end{rem}

	\begin{rem}
		Discrete $G$-modules over $X$ form an Abelian category in the obvious way. More precisely, a morphism  $$\mathcal M=(M_x, \phi_{g,x}) \lra \mathcal M'=(M'_x, \phi'_{g,x})$$ is the data of additive maps $$f_x: M_x \lra M'_x,$$ one for each $x \in X$, such that $$\phi'_{g,x} \circ f_x = f_{gx} \circ \phi_{g,x},$$ for all $x \in X$ and all $g \in G$.
	\end{rem}
	If $\mathcal M=(M_x, \phi_{g,x})$ is a discrete $G$-module over $X$, we can form the direct sum $$N(\mathcal M):=\bigoplus_{x \in X} M_x;$$ it is a $G$-module in an obvious way, given by applying the $\phi_{g,x}$'s.
	The association $$\mathcal M \mapsto N(\mathcal M)$$ is a functor, from the category of discrete $G$-modules over $X$ to that of discrete $G$-modules. It plays the r\^ole of a trace map,  and is a categorical formulation of the  usual induction process, from open subgroups of $G$. We now explain why.\\
	
	Assume that $$X = G/H,$$ for $H \subset G$ a nontrivial open subgroup.  Denote by $x_0 \in X$ the class of the neutral element. \\
	Then  we have a functor $$\mathcal M=(M_x, \phi_{g,x}) \lra M_{x_0},$$ from the category of discrete $G$-modules over $X$ to that of discrete $H$-modules, where  $M_{x_0}$ is considered as an  $H$-module via the maps $\phi_{h,x_0}$, for $h \in H=\mathrm{Stab(x_0)}. $ \\It is not hard to see that this functor is an equivalence of categories. The proof is left to the reader as an exercise.
	
	\begin{rem}
		What precedes is a concrete example of the following philosophical statement:  if $X=G/H$, a $G$-equivariant structure over the base $X$ is nothing but an $H$-equivariant structure.
	\end{rem}
	
Now, let $\mathcal M=(M_x, \phi_{g,x})$  be a $G$-module over $X$. Put $M:=M_{x_0}$, seen  as a discrete $H$-module. Then  $$N(\mathcal M)=\bigoplus_{x \in X} M_x$$ is canonically isomorphic to  the induced module $\mathrm{Ind}_H^G(M)$. Note that, since $H$ has finite index in $G$, this induced module can be defined either by the formula $$\mathrm{Ind}_H^G(M)=M \otimes_{\Z[H]} \Z[G],$$ or by $$\mathrm{Ind}_H^G(M) = \mathrm{Maps}_H(G,M),$$ the group of $H$-equivariant maps from $G$ to $M$ ('induction=coinduction' in this case). \\
	Now, recall Shapiro's Lemma -which is elementary but of crucial importance in this paper- asserting that the cohomology groups $H^n(G,\mathrm{Ind}_H^G(M))$ and $H^n(H,M)$ are canonically isomorphic. Putting what we just said together, we get the following statement.
	
	\begin{prop}
		Put $$X = G/H,$$ for $H \subset G$ a nontrivial open subgroup.  Denote by $x_0 \in X$ the  neutral class. \\ Let $\mathcal M=(M_x, \phi_{g,x})$ be a discrete $G$-module over $X$. Then $M_{x_0}$ is canonically a discrete $H$-module, and Shapiro's lemma yields canonical isomorphisms $$H^n(G,\bigoplus_{x \in X} M_x) \stackrel \sim \lra H^n(H,M_{x_0}),$$ for each $n \geq 0$.
	\end{prop}

	\begin{rem}
		If $X$ is an arbitrary finite $G$-set and $\mathcal M=(M_x, \phi_{g,x})$ is a discrete $G$-module over $X$, we can adapt the preceding Proposition, yielding canonical isomorphisms $$H^n(G,\bigoplus_{x \in X} M_x) \stackrel \sim \lra  \bigoplus_{i=1}^m H^n(G_i,M_{x_i}),$$ where the $x_i's$ form a system of representatives of $G$-orbits in $X$, and where $G_i$ is the stabilizer of $x_i$.
		
	\end{rem}%

	To finish this section, let us give a typical example of how this Remark will be applied.\\
	Let $k$ be a finite field of characteristic $p$. Let $V$ be a $(k,G)$-module. Put $$X:=\P(V);$$ it is obviously a finite $G$-set.\\
	There is a 'tautological' discrete $G$-module over $X$, which is $\mathcal M$, defined by $$\mathcal{M}_L:= V/L,$$ for each line $L \in X$, and where the map $$\phi_{g,L}: V/L \lra V/g(L)$$ is induced by the linear  map $v \mapsto g.v$. Shapiro's Lemma  then yields canonical isomorphisms $$H^n(G,N(\mathcal M))=H^n(G,\bigoplus_{L \in \P(V)} V/L) \stackrel \sim \lra  \bigoplus_{i=1}^m H^n(G_i, V/L_i),$$ like we just discussed in the Remark above. This fundamental fact will be crucial in the proof  the Stable Lifting Theorems.
	
\section{The Tense product.}\label{tenseproduct}

In this section, we elaborate on the lack of commutation between the (usual) tensor product, and Pontryagin duality. This phenomenon is at the heart of our approach. We thus chose to introduce some related notions (the valuation of a module, its trace and its Chern character)  that we find colorful- though they will not be used in the rest of this paper.

Let  $k$ be a perfect field of characteristic $p >0$. We denote by $\{\W_n(k)-Mod\}$ the category of (arbitrary) $\W_n(k)$-modules. In what follows, the symbol $\otimes$ means $\otimes_{\W_n(k)}$. \\  Recall that, for any  $\W_n(k)$-module $M$, we have an 'equality' $$ M^\vee= \Hom_{\W(k)}(M, K/\W(k))=  \Hom_{\W_n(k)}(M, \W_n(k)).$$ We will now see that it strongly depends on $n$.\\
Let $M$ and $N$ be two $\W_n(k)$-modules. Through the preceding identification, we have a pairing $$\Phi_{n,M,N}: (M \otimes N) \times (M^\vee \otimes N^\vee) \lra \W_n(k),$$ $$ ((m \otimes n),(\phi \otimes \psi)) \mapsto \phi(m) \psi(n).$$ It is perfect if (and only if) $M$ or $N$ is a finite and free $\W_n(k)$-module. What is more, it is somewhat badly behaved, in the sense that it strongly depends on $n$: if $m \geq n$ is another integer, viewing $M$ and $N$ as $\W_m(k)$-modules yields the formula $$\Phi_{m,M,N}= p^{m-n} \Phi_{n,M,N}, $$ as pairings with values in $K/ \W(k)$. In particular, if  $M$ and $N$ are actually $k$-vector spaces, the pairing $\Phi_{2,M,N}$ is zero!\\
Thinking further, we see that the category of  $\W_n(k)$-modules is actually equipped with (at least) \textit{two} tensor product structures: the usual one, and the one given by the formula $$M \tilde \otimes N:= (M^\vee \otimes N^\vee)^\vee .$$ If $n\geq 2$, there is no canonical isomorphism (of bifunctors) $M \tilde \otimes N \simeq M  \otimes N$. 
This fact is a simple algebraic analogue of Grothendieck's theory of tensor products of topological vector spaces. Note that the Tense Product defined shortly will be ubiquitous later, when we mod out the kernel of Pontryagin duality for divided powers (see for instance Proposition \ref{polytense}).  \\

\begin{defi}(Tense Product.)
Let  $M$ and $N$ be two  $\W_n(k)$-modules. We put $$M \overline \otimes_n N:= (M \otimes N)/\Ker(\Phi_{n,M,N}) ,$$ where $$\Ker(\Phi_{n,M,N})=\{x \in M \otimes N,~~~\Phi_{n,M,N}(x,y)=0,~~~ \forall y \in M^\vee \otimes N^\vee \}. $$ It is a  $\W_n(k)$-module, called the Tense Product of the $\W_n(k)$-modules $M$ and $N$. If the dependence in $n$ is clear from the context, we shall denote it simply by $M \overline \otimes N$.

\end{defi}

\begin{rem}\label{tenseconcrete}

Assume that $M=\W_a(k)$ and $N=\W_b(k)$, with $1\leq a,b \leq n$. We then have $$ M \otimes_{\W_n(k)} N = \W_{\min(a,b)}(k),$$ whereas $$ M \overline\otimes_n N = \W_{a+b-n}(k),$$ with the convention that $\W_i(k)=0$ for negative $i$.

\end{rem}
\begin{lem} 

The category $\{\W_n(k)- Mod\}$, equipped with the Tense Product, is a symmetric monoidal category, with coherence axioms inherited from those of the usual tensor product. The Tense Product commutes with Pontryagin duality of finite modules:  for two finite $\W_n(k)$-modules, we have a canonical isomorphism $$(M \overline \otimes N) ^\vee \simeq M^\vee \overline \otimes N^\vee.$$ 
\end{lem}

 \begin{dem}
     This is routine check. Let us perhaps explain why the Tense Product is (bi)functorial. Let $f:M\lra N$ and $f':M' \lra N'$ be morphisms between $\W_n(k)$-modules. Then, for $m \in M$, $m' \in M'$, $\phi \in N^\vee$ and $\phi' \in {N'}^\vee$, one has $$ \Phi_{n,N,N'}(f(m) \otimes f'(m') , \phi \otimes  \phi')= \phi(f(m)) \phi'(f'(m'))$$ $$= \Phi_{n,M,M'}(m \otimes m , f^\vee(\phi) \otimes  {f'}^\vee(\phi')).$$ This adjunction formula shows that $$f \otimes f': M \otimes M' \lra N \otimes N' $$ passes to the quotient by the kernel of the duality $\Phi_n$, yielding a linear map  $$f  \overline \otimes f': M \overline\otimes M' \lra N \overline\otimes N'. $$ 
 \end{dem}

 We can now state the following definition.
 
\begin{defi} 

We will denote by $\mathcal W_n$ the symmetric monoidal category of $W_n(k)$-modules, with monoidal structure given by the Tense Product.  It is equipped with  the perfect duality $M \mapsto M^\vee$, induced by Pontryagin duality.
\end{defi}

As one may guess, reducing mop $p^n$ yields a (lax monoidal) functor from $\mathcal W_{n+1}$ to $\mathcal W_n$. 

\begin{lem}
Let $M$ and $N$ be $\W_{n+1}(k)$ -modules. Then the natural map $$ (M/p^n) \otimes_{\W_{n}(k)} (N/p^n)  \lra (M \otimes_{\W_{n+1}(k)} N)/p^n  $$  induces by passing to the quotient a $\W_{n}(k)$-linear map $$ (M/p^n)  \overline \otimes_n (N/p^n)  \lra (M \overline \otimes_{n+1} N)/p^n.$$

\end{lem}
\begin{dem}
To check this, one may assume that  $M=\W_a(k)$ and $N=\W_b(k)$ are of rank one, and use Remark  \ref{tenseconcrete}.

\end{dem}
\begin{defi}
   The functor $$ \mathcal W_{n+1} \lra \mathcal W_{n} $$ $$M \lra M /p^n $$ will be denoted by $\Theta$. Thanks to the previous Lemma, it is a lax monoidal functor.
    
\end{defi}

  \begin{defi}(Tense Algebra, Tense Symmetric Powers, Tense Exterior Powers.)

Let  $M$  be a  $\W_n(k)$-module. For every nonnegative integer $i$, we put $$  M^{\overline \otimes_n^i}:= \underbrace { M \overline \otimes_n M \overline \otimes_n \ldots \overline \otimes_n M}_{i \> \mathrm{times}}.$$ We set $$\overline T_n(M):=\bigoplus_{i=0}^\infty   M^{\overline \otimes_n^i}; $$ it is naturally a $\W_n(k)$-algebra, the Tense Algebra of $M$. We define $$\overline \Sym_n(M):= \bigoplus_{i=0}^\infty \overline \Sym_n^i(M) $$ to be the largest commutative quotient of $\overline T_n(M)$.  As usual, it is obtained by modding out the ideal spanned by the elements $$x \otimes y - y \otimes x, $$ for $x, y \in M$. It is the Tense Symmetric Algebra of $M$.  Similarly, we define $$\overline \Lambda_n(M):= \bigoplus_{i=0}^\infty \overline \Lambda_n^i(M) $$ to be the   quotient of $\overline T_n(M)$ obtained by modding out the ideal spanned by the elements $$x \otimes x,$$ for $x \in M$.  It is the Tense Exterior Algebra of $M$.
\end{defi}
 
 \begin{defi}
 A tense algebra is an algebra in the category $\mathcal W_n$. In other words, it is a $\W_n(k)$-algebra, such that the multiplication map $$\mu: A \otimes_{\W_n(k)} A \lra A $$ factors through the natural quotient map $$  A \otimes_{\W_n(k)} A  \lra A  \overline \otimes_n A .$$
	 \end{defi}
	 
	 \begin{ex}
	 The tense symmetric (resp. exterior) algebra of an arbitrary $\W_n(k)$-module is  of course a tense algebra.\\
	 If $A$ is is a usual $\W_n(k)$-algebra which is flat (=free) as a $\W_n(k)$-module, it is automatically tense.
	 \end{ex}
	 
	 \begin{exo}
	 Let $A \in \mathcal W_n$ be a tense $\W_n(k)$-algebra. Show that the unit $1 \in A$ has (additive) order $p^n$. In other words, it spans a free direct summand of rank one of $A$, as a $\W_n(k)$-module.
	 \end{exo}

\begin{defi}
    
    Let $M$ be a (finite) $\W_n(k)$-module. Consider the largest  number $$i\in \{-n, -n+1, \ldots, -1,0\}$$ such that $M /p^{n+i} $ is a free $\W_{n+i}(k)$-module.\\ We put $$v_n(M)=i . $$ It is the valuation of $M$.
\end{defi}

\begin{rem}
The $\W_n(k)$-module $M$ has  valuation $0$ if and only if it is free. In case it is not, its valuation is strictly negative, and can be though of as the highest 'pole' of $M$.

\end{rem}
\begin{lem}
     Let $M$ and $N$ be (finite) $\W_n(k)$-modules. Then we have $$ v_n(M \overline \otimes_n N) = v_n(M)+ v_n(N), $$ and  $$ v_n(M \oplus N) = \min (v_n(M), v_n(N))$$ with the convention that  all integers $\leq -n$ are identified to $-n$.
\end{lem}

\begin{dem}
Straightforward.
\end{dem}

\begin{lem}
    
    Let $M$ be a (finite) $\W_n(k)$-module. Then  $v_n(M)$  is the largest  negative number  $i$, such that the composite $$ M \otimes_{\W_n(k)}( M^\vee) \stackrel {\ev}\lra \W_n(k) \lra \W_{n+i}$$ passes to the quotient by $\Ker(\Phi_{n,M,M^\vee})$.
\end{lem}

\begin{dem}

Assume first that  $M=\W_a(k)$, where $a$ is an integer, with $1 \leq a \leq n$.  Then the evaluation arrow under consideration is $$ \W_a(k) \lra \W_n(k),$$ $$ 1 \mapsto p^{n-a}.$$ On the other hand, $\Ker(\Phi_{n,M,M^\vee})$ is generated by $ p^{2a-n}$, or is everything if $a \leq \frac n 2.$ We then see that the composite under consideration factors through $\Ker(\Phi_{n,M,M^\vee})$ if, and only if, $$p^{2a-n+
n-a}=p^a=0 \in \W_{n+i}(k),$$ meaning that $i \leq a-n= v_n(M)$. The general case follows.

\end{dem}

\begin{defi}
    
    Let $M$ be a (finite)  $\W_n(k)$-module. By the preceding Lemma, the composite $$ M \otimes_{\W_n(k)} M^\vee \stackrel {\ev}\lra \W_n(k) \lra \W_{n+v_n(M)}(k) $$ induces an arrow  $$ M \overline \otimes_n M^\vee \lra \W_{n+v_n(M)}(k),$$ which we denote by $\tr_M$. It is the trace of $M$.
   
\end{defi}

We conclude this section by defining the Chern character of a finite $\W_n(k)$-module, and stating a first property. 
\begin{defi}
Let  $M$ be a finite  $\W_n(k)$-module. We put $$ Ch_n(M):=a_n + a_{n-1} X^{-1}+  a_{n-2} X^{-2}+ \ldots + a_1 X^{-n+1} \in \Z [X^{-1}].$$  

\end{defi}

\begin{lem}
Let  $M$ and $N$ be two  finite $\W_n(k)$-modules. Then  $$ v_n(M)=  v_X(Ch_n(M)).$$ Modulo $X^{-n}\Z[X^{-n}]$, we have
  $$ Ch_n(M \overline \otimes^n N)=Ch_n(M)  Ch_n(N).$$

\end{lem}

\begin{dem}

Follows from  Remark \ref{tenseconcrete}.
\end{dem}

 \begin{rem}
 We believe it could be interesting, in the future, to investigate the behaviour of the  Chern character of Omega powers of $k$-vector spaces (see section \ref{omega}). For instance, if $V$ is a finite-dimensional $k$-vector space,  what can we say about the map $$ n \in \N \mapsto Ch_{n+1}(\Omega^n(V)) ?$$ We should, of course, also extend these considerations to modules over global rings.
 \end{rem}

\section{Divided powers. }\label{gamma}

For a nice and short account on properties of divided powers,  we refer the reader to \cite{Fe}. A more comprehensive study of divided powers can be found in  \cite{Ro}, which contains all the proofs of the Propositions which we state here without proof.

In this section, $A$ is a commutative ring.


\begin{defi}
Let $M$ be an  $A$-module. We denote by $\Gamma_A(M)$ (or simply by $\Gamma(M)$ if the dependence in $A$ is clear) the graded divided power algebra of $M$,  defined as follows. It is generated by degree $i$ symbols $[x]_i$, for each $i \in \N$ and each $x\in M$, with relations:\\

\noindent i) $[x]_0=1$,\\
\noindent ii)$ [x+x']_n= \sum_0^n [x]_i [x']_{n-i}$,\\
\noindent iii)$ [\lambda x]_n=\lambda^n[x]_n$,\\
\noindent iv) $[x]_n[x]_m= {n+m \choose n}[x]_{n+m}$.\\

\noindent We define $\Gamma^n(M)$ to be the homogeneous component of degree $n$ of  $ \Gamma(M)$. We put $\Gamma^+(M):=\oplus_{n \geq 1} \Gamma^n(M)$; it is an ideal of  $\Gamma(M)$.
\end{defi}

\begin{rem}
As it is well-known, the symbol $[x]_n$ plays the r\^ole of $ \frac 1 {n!}  x^n$. More precisely, if $n!$ is invertible in $A$ (which typically happens if $A$ has prime characteristic $p$ and $n=p-1$), the  natural map $$ \Sym_A^n(M) \lra \Gamma^n_A(M),$$ $$ x_1 \otimes \ldots\otimes x_n \mapsto  [x_1]_1 \ldots [x_n]_1$$ is an isomorphism, with inverse given by     $$ \Gamma^n_A(M) \lra \Sym_A^n(M),$$ $$ [x]_n \mapsto \frac 1 {n!}  x^n .$$ At this point, the reader may wonder whether the last formula makes any sense. Why does it yield a well-defined $A$-linear map? This will become apparent in a moment, using the viewpoint of polynomial laws.\\

\end{rem}

\begin{rem}\label{multi}
Equality iv), applied several times, yields the formula $$[x]_{n_1} \ldots [x]_{n_r} = {{n_1+ \ldots + n_r } \choose {n_1, \ldots, n_r}}[x]_{n_1+\ldots+n_r} ,$$ where $$  {{n_1+ \ldots + n_r } \choose {n_1, \ldots, n_r}} = \frac {(n_1+ \ldots+ n_r)!} {n_1 ! \ldots n_r!}$$ is the usual multinomial coefficient.
\end{rem}

For each positive integer $i$, the ideal $\Gamma^+(M)$ is moreover equipped with an operator $$\gamma_i : \Gamma^+(M) \lra \Gamma^+(M),$$ $$ x \mapsto \gamma_i(x),$$ playing the role of $x \mapsto x^i / i! $, which endows $(  \Gamma(M),  \Gamma(M)^+)$ with the structure of an $A$-algebra with divided powers. Let us be more precise.

\begin{prop}\label{defigammai}
    Let $M$  be an $A$-module.
For each positive integer $i \geq 0$, the polynomial law $$M \lra  \Gamma^+(M),$$ $$x \mapsto [x]_i,$$ uniquely extends to a polynomial law  $$\gamma_i : \Gamma^+(M) \lra \Gamma^+(M),$$ which is homogeneous of degree $i$, such that the following conditions hold. 

1) The $\gamma^i$'s are functorial in $A$ and $M$.\\
2) We have $\gamma_1= \Id$.\\
3) We have $$\gamma_i(x+y)= \sum_{a+b=i} \gamma_a(x) \gamma_b(y),$$ identically.\\
4) We have $$\gamma_j \circ \gamma_i=  \frac {(ij)!} {j! (i!)^j} \gamma_{ij}.$$

Moreover, the four properties above uniquely determine the $\gamma^i$'s.
\end{prop}

\begin{prop} \label{polyfunctor}
Let $M$, $N$ be $A$-modules.  We have a canonical isomorphism $$\Gamma^n(M \oplus N) \simeq \oplus_{i=0}^n (\Gamma^i(M) \otimes_A \Gamma^{n-i}(N)).$$
\end{prop}

\begin{rem}
	The previous Proposition says that divided power functors are strictly polynomial, in the sense of  \cite{FFSS}.
\end{rem}
\begin{prop}\label{divbase}
Let $M$ be an $A$-module, and let $B$ be a commutative $A$-algebra.  We have a canonical isomorphism of graded rings $$\Gamma_A(M) \otimes B  \simeq  \Gamma_B(M \otimes_A B).$$
\end{prop}

\subsection{Polynomial laws.}

Let $A$ be a commutative ring.

\begin{defi}
	If $M$ 	is an $A$-module, we denote by $\underline M$ the functor 
	$$ R \mapsto M \otimes_A R,$$ from the category of commutative $A$-algebras to that of  sets.
	
\end{defi}

\begin{defi}
	Let $M,N$ be $A$-modules. A polynomial law from $M$ to $N$ is a morphism 
	of functors $$F: \underline M \lra \underline N.$$ We shall say that $F$ is homogeneous of degree $n \geq 0$ if, for every commutative $A$-algebra $R$ and every $t \in R$ and $ m \in M \otimes_A R,$ we have $$F(tm)=t^nF(m).$$
	
\end{defi}

\begin{rem}
	One can show that a degree $0$ (resp. degree $1$) polynomial law is obtained from a constant (resp. $A$-linear) map $M \lra N$.
	
\end{rem}

\begin{rem}
	Slightly abusing notation, we will sometimes denote a polynomial law $$F: \underline M \lra \underline N$$ simply by $$F: M \lra  N,$$  dropping the underscore. We shall do so only if there is no chance of confusing $F$ with a mere map.
	
\end{rem}
\begin{rem}
Let $M$ be an $A$-module. Let $i,j$ be positive integers.
 The divided power operation $$\gamma_i : \Gamma_A^j (M) \lra  \Gamma_A^{ij} (M)$$ is a polynomial law, which is homogeneous, of degree $i$. It will often be considered as such in the sequel.
	
\end{rem}

\begin{rem}
	If $V$ and $W$ are locally free $A$-modules of finite rank, then a  polynomial law from $V$ to $W$  is nothing but a morphism of affine $A$-schemes $$\A_A(V) \lra \A_A(W).$$

\end{rem}

The next Proposition is crucial. Its content is that the functor $\Gamma$ represents the functor of polynomial laws.

\begin{prop}
Let $M$, $N$ be $A$-modules. Then $\Hom_A(\Gamma^n(M),N)$ is canonically isomorphic to the group of polynomial laws from $M$ to $N$, which are homogeneous of degree $n$.

\end{prop}

The previous Proposition admits an obvious generalization,  as follows.

\begin{prop}\label{multihomog}
Let $M_1, M_2, \ldots, M_r$ and $N$ be $A$-modules. Let $n_1, \ldots, n_r$ be positive integers. Then $$\Hom_A(\Gamma^{n_1}(M_1) \otimes_A \Gamma^{n_2}(M_2) \otimes_A \ldots \otimes_A \Gamma^{n_r}(M_r),N)$$ is canonically isomorphic to the  ($A$-module of) polynomial laws $$\underline {M_1} \times  \underline {M_2} \times \ldots \times \underline  {M_r} \lra \underline N,$$ which are homogeneous of degree $n_i$ in $M_i$ (for $i=1 \ldots r$).

\end{prop}

For $M$ an $A$-module, the association $$M \lra (M^{ \otimes n})^{\mathcal S_n},$$ $$ x \mapsto x^{\otimes n},$$ is obviously a polynomial law, which is homogeneous of degree $n$. It thus induces an $A$-linear morphism $$ F_n(M): \Gamma_A^n(M) \lra  (M^{ \otimes n})^{\mathcal S_n}.$$

\begin{prop}\label{gammafree}
If $M$ is locally  free of finite rank, the morphism $F_n(M)$ above  is an isomorphism.

\end{prop}
	
	\begin{rem}\label{symdualgamma}
	If $M$ is locally free of finite rank, the $A$-dual of $ (M^{ \otimes n})^{\mathcal S_n}$ is nothing but the symmetric power $\Sym_A^n(M^*)$. Thus, the formation of divided powers, for finite locally  free modules, is dual to that of symmetric powers.  
	\end{rem}
	
	\subsection{Divided powers and duality.}
 
Let us now mention a nice compatibility of divided powers, with respect to duality.\\

	\begin{defi}\label{duality}
		Let $M$ be an $A$-module.   Let $n \geq 1$ be an integer. The formula $$(v,\phi)  \mapsto \phi(v)^n$$ defines a polynomial law (of $A$-modules) $$\underline M \times \underline{M^*}  \lra  A,$$ which is bihomogeneous, of bidegree $(n,n)$. By Proposition \ref{multihomog}, this law corresponds to a pairing $$ \Gamma_A^{n}(M) \times  \Gamma_A^{n}(M^*) \lra A .$$ We  will  denote this pairing by $\Delta_{n}$, or by $\Delta$, or  even by $<.,.>$, if the context is clear enough. On the level of pure symbols, we have $$<[\phi]_{n},[v]_{n}>=\phi(v)^n.$$
	\end{defi}

\begin{lem}\label{formuladual}
    Let $M$ be an $A$-module.  The following assertions are true.\\
    
    1) For $x_1 , \ldots, x_r \in M$, $\phi \in M^*$ and $i_1, \ldots i_r$  positive integers, we have $$ <[x_1]_{i_1}\ldots [x_r]_{i_r},[\phi]_{i_1+ \ldots + i_r}> = {{i_1+ \ldots + i_r} \choose {i_1,\ldots,i_r}} \phi(x_1)^{i_1}\ldots  \phi(x_r)^{i_r}$$ $$={{i_1+ \ldots + i_r} \choose {i_1,\ldots,i_r}} <[x_1]_{i_1},[\phi]_{i_1}>\ldots <[x_r]_{i_r},[\phi]_{i_r}> .$$

    2) Let $m$ and $n$ be positive integers. For $X \in \Gamma^n(M)$ and $\phi \in M^*$ the formula  $$ <\gamma_m(X), [\phi]_{mn}>= \frac{ (mn)!} {(n!)^m m!} <X, [\phi]_{n}>^m \in A $$ holds (note that the integer coefficient appearing here is the number of partitions of a set of cardinality $mn$ in $m$ subsets of cardinality $n$).

\end{lem}
	
	\begin{dem}
	We can assume that $M$ is an $A$-module of finite type.\\
	If $f: L \lra M$ is a surjective linear map between $A$-modules, then the Lemma holds for $M$ if it holds for $L$. Hence, we are reduced to the case where $M=A^n$ is free. Let $g: B \lra A $ be a surjective homomorphism of commutative rings. Let $N$ be a $B$-module, for which the Lemma holds (over the ring $B$, of course). Then the Lemma holds for $M=N \otimes_B A$ (use Proposition \ref{divbase}). Altogether, we can  assume that $M$ is a free module and that $A$ is a domain of characteristic zero. But then, divided powers of $M$ are free $A$-modules, and we can check everything after extending scalars to the fraction field $F$ of $A$. We are thus reduced to the case where $A=F$ is an algebraically closed field of characteristic zero, in which case pure symbols additively span divided power modules. We can then  identify $\Gamma_A^n(M)$ and $\Sym_A^n(M)$, through the symmetrizing operator ($[x]_n$ corresponds to $\frac {x^n} {n!}$). Then the duality $\Delta$ is given by the usual (reassuring, but awkward) formula $$ \Sym_A^n(M) \times \Sym_A^n(M^*) \lra A,$$ $$ (x_1 \otimes \ldots \otimes x_n , \phi_1 \otimes \ldots \otimes \phi_n ) \mapsto n! \sum_{\sigma \in S_n}  \phi_1(x_{\sigma(1)}) \ldots \phi_n(x_{\sigma(n)}) .$$  Checking point 1) is then straightforward. \\
	For 2), write $$X=\sum_i [x_i]_n. $$ We compute, using 1):  $$ <\gamma_m(X), [\phi]_{mn}>=  <\frac {X^m} {m!}, [\phi]_{mn}>=\frac 1 {m!} <(\sum_i [x_i]_n )^m, [\phi]_{mn}>$$ $$ = \frac 1 {m!} <\sum_{i_1, \ldots, i_m} [x_{i_1}]_n \ldots  [x_{i_m}]_n , [\phi]_{mn}> $$ $$= \frac 1 {m!} {{mn} \choose {n,n, \ldots, n}}  \sum_{i_1, \ldots, i_m} \phi(x_{i_1})^n \ldots  \phi(x_{i_m})^n$$ $$= \frac 1 {m!} {{mn} \choose {n,n, \ldots, n}}  \sum_{i_1, \ldots, i_m} <[x_{i_1}]_n, [\phi]_n>  \ldots <[x_{i_m}]_n, [\phi]_n>$$ $$= \frac 1 {m!} {{mn} \choose {n,n, \ldots, n}}  (<[x_{i_1}]_n, [\phi]_n> + \ldots+ <[x_{i_m}]_n, [\phi]_n>)^m $$  $$=  \frac{ (mn)!} {(n!)^m m!} <X, [\phi]_{n}>^m$$ (here we write sums  over all integers $i_1, \ldots, i_m$, repetitions allowed, to avoid  multinomial coefficients). The proof is over.

	\end{dem} 
	
\subsection{Divided powers versus symmetric powers.}\label{PDVSSIM}

At this point, it seems legitimate to compare symmetric powers and divided powers more closely.  Victory shall belong to the latter, and by far: they are much simpler, versatile and better behaved- for many reasons. We mention two of them.\\

 1)   Let $V$ and $W$ be two $A$-modules. On the one hand, $\Hom_A(\Sym_A^n(V),W)$ corresponds to symmetric $n$-multilinear forms $$F: V ^n \lra W.$$  Such a form is defined by an expression of the shape $$(v_1, \ldots, v_n) \mapsto F(v_1, \ldots, v_n).$$ It is a function of the $n$ variables $v_1, \ldots, v_n$. \\ On the other hand, $\Hom_A(\Gamma^n(V),W)$ corresponds to polynomial  laws from $V$ to $W$, which are homogeneous of degree $n$. Such a law is defined by an expression of the shape $$v \mapsto F(v).$$  Being a natural transformation between the functors $\underline V$ and  $\underline W$, the expression $F(v)$ has to functorially make sense for any commutative $A$-algebra $R$, and every $v \in V \otimes_A R$.  However,  it depends on a -single- variable $v$. In that sense, it is much easier to define than  a symmetric $n$-multilinear form.\\
 
 2) Let $m$ be positive integer. Symmetric powers (or tensor powers) of a module, which is of $m$-torsion, will remain of $m$-torsion, regardless to the ring of coefficients. This is far from being so for divided powers- a major fact which is at the heart of this paper. The underlying phenomena (in prime power torsion) will be studied extensively, starting from the next section. For instance, we shall see that, if $L$ is a free $(\Z/ p^n \Z)$-module of rank one, then $\Gamma^{p^s}_{\Z}(L)$ is  a free  $(\Z/ p^{n+s} \Z)$-module of rank one.\\
 

	\section{The Teichm\"uller representative, as a polynomial law.}
Among polynomial laws, there are fundamental ones, given by  Teichm\"uller representatives in truncated Witt vectors. Let us be more precise. We begin with a standard but extremely important Lemma, at the heart of $p$-adic theory.

\begin{lem} \label{Teichpoly}
Let $A$ be a commutative ring.\\
The map $$A \lra A/p^{n+m}A,$$ $$ x \mapsto x^{p^n},$$ factors through the quotient $A \lra A/p^m A$. Since this hold functorially in $A$, we get this way a polynomial law of $\Z$-modules
 $$\Z /p^m \Z \lra  \Z /p^{m+n} \Z ,$$ $$ x \mapsto x^{p^n}.$$

\end{lem}

\begin{dem}
By induction, it is enough to check the claim for $n=1$. In this case, for any $x,y \in A$, we  have the well-known congruence $$(x+p^m y)^p \equiv x^p $$ modulo $p^{m+1}A$, whence the claim.

\end{dem}

\begin{defi} \label{Teichpolydef}
The  polynomial law (of $\Z/p^{m+n}\Z$-modules)
 $$\Z /p^m \Z \lra  \Z /p^{m+n} \Z ,$$ $$ x \mapsto x^{p^n}$$ will be denoted by $\overline \tau_n$.

\end{defi}

\begin{rem}
Note that, for any $x \in k$,  we have $$ \overline \tau_n(x) =\tau(x^{p^n}) \in \W_{n+1}(k),$$ where $\tau$ is the usual  Teichm\"uller representative. Looking closely at this equality reveals the following.\\
In the theory of Witt vectors over perfect fields of characteristic $p$, the expression $\tau(x) \in \W_{n+1}(k)$ is well-defined because we can extract a (unique) $p^n$-th root of $x$ in our base ring $k$. The expression $\overline \tau_n(x) \in A$, in contrast, makes sense for any base ring $A$ of characteristic $p^{n+m}$, and any $x \in A/p^m$.
\end{rem}
 We  need a simple arithmetic Lemma.

\begin{lem}\label{coeffprimetop}
	
 Let $ a_1, \ldots, a_r$ be $r$ nonnegative integers. Put $n= a_1+ \ldots + a_r $. \\ The following assertions are true.\\ i) The number of carryovers in the base-$p$ addition $$n=( \ldots (a_1+a_2)+a_3)+ \ldots )+a_r$$ does not depend on the order of the $a_i$'s. It is equal to the $p$-adic valuation of the multinomial coefficient $ { n  \choose {a_1, \ldots, a_r}}$ .\\
 ii) We have $$ v_p( { pn \choose {pa_1, \ldots, pa_r}} )= v_p( { n \choose {a_1, \ldots, a_r}})$$ and  $$ v_p( { n \choose {a_1, \ldots, a_r}} ) \geq \mathrm{max}_i \{v_p(n)-v_p(a_i)\}.$$
iii)	 In the particular case where $r=2$ and $a_1+a_2=\iota p^m$, with $m \geq 0$ and $1 \leq \iota \leq p-1$,  the inequality in ii) is an equality. In other words, we have $$v_p( {{\iota p^m} \choose a_1, a_2 })=m-v_p(a_1)=m-v_p(a_2).$$ 
\end{lem}

\begin{dem}

For i), we can reduce to the case $r=2$ using the formula $$ { a_1+ \ldots + a_r \choose {a_1, \ldots, a_r}} = { a_1+ \ldots + a_r \choose {a_1+a_2, a_3, \ldots, a_r}}  {a_1 + a_2 \choose a_1, a_2 }.$$ The claim to prove is then a classical fact, which is also a nice elementary exercise left to the reader.\\ Assertion ii) is an easy consequence of i).\\
Let us now prove iii). The case $m=0$ is obvious; we thus assume that $m\geq 1$. By ii), we then note that  $$v_p( {{\iota p^{m+1}} \choose pa_1, pa_2 })=  v_p( {{\iota p^m} \choose a_1, a_2 }),$$ allowing us to reduce to the case  where $a_1$ and $a_2$ are prime-to-$p$. Write $$a_i=b_i+c_ip^m$$ for $i=1,2$,  with $b_i$   prime-to-$p$, and  $b_i \leq p^m-1$. The number $b_1+b_2$ is divisible by $p^m$, hence equals $p^m$. We thus have $$1+c_1+c_2=\iota \leq p-1,$$ from which we infer that the number of carryovers  in the base-$p$ addition of $a_1$ and $a_2$ equals that of $b_1$ and $b_2$, which is obviously $m$. The claim is proved.
	
\end{dem}

\begin{rem}
    The formula $$ v_p( { pa_1+ \ldots+ pa_r \choose {pa_1, \ldots, pa_r}} )= v_p( { a_1 +\ldots+a_r \choose {a_1, \ldots, a_r}})$$ suggests that the function $$ (a_1, \ldots, a_r) \mapsto v_p( { a_1 +\ldots+a_r \choose {a_1, \ldots, a_r}})$$ behaves like a height on the projective space $\P^{r-1}$, locally at $p$. This is perhaps worth investigating.
\end{rem}

\begin{lem}\label{torsionbornee}
Let $V$ be an $A$-module, such that $p^mV=0$, for some $m \geq 0.$ Let $n$ be a positive integer. Then $\Gamma_A^n(V)$ is of $p^{m+v_p(n)}$-torsion.
\end{lem}

\begin{dem}
	Note first the following two  obvious facts. \\
		i) Since $V$ is annihilated by $p^m$,   $\Gamma_A^n(V)$ is annihilated by $p^{mn}$. In particular, it is a $\Z_p$-module.\\
	ii) Let  $$n = a_1+ \ldots + a_r$$  be a decomposition of   $n$ into a sum of $r$ nonnegative integers. For $i=1 \ldots r$, let $v_i$ be an element of $V$. Then the (additive) order of $$[v_1]_{a_1} \ldots [v_r]_{a_r} \in \Gamma_A^n(V)$$ is at most the minimum of the orders of the elements $[v_i]_{a_i} \in  \Gamma_A^{a_i}(V)$.\\
	
 We now show that  $$p^{m+v_p(n)}[v]_n=0  \in  \Gamma_A^n(V),$$ for each $v \in V$.  Let $$n=a_0+a_1p+\ldots a_r p^r$$ be the base-$p$ expansion of $n$. In the equality $$[v]_{a_0} [v]_{a_1 p} \ldots  [v]_{a_r p^r }={n \choose {a_0, a_1p , \ldots, a_r p^r} }[v]_n \in  \Gamma_A^n(V),$$ the multinomial coefficient is prime-to-$p$, by Lemma \ref{coeffprimetop}.  We can thus assume that $n=\i p^r$, with $1\leq \i \leq p-1$.  In the formula $$[v]_{p^r}^\iota={\iota p^r \choose {p^r, \ldots, p^r} } [v]_n,$$  the multinomial coefficient is prime-to-$p$, by Lemma \ref{coeffprimetop} again. We thus reduce to the case $\i=1$, i.e.   $n=p^r$.
 
We then see that $$[v]_{p^{r-1}}^p={p^{r} \choose {p^{r-1} , \ldots, p^{r-1} } }[v]_n ,$$ and that the multinomial coefficient occuring in this formula has $p$-valuation one.  Induction on $r$ yields the result, since $V=\Gamma_A^1(V)$  is of $p^m$-torsion by assumption.\\

\end{dem}\\

\subsection{Divided powers of torsion $\W(k)$-modules.}
From now on, $k$ is a perfect field of characteristic $p$. \\
We will explain the first steps towards a canonical process for lifting torsion $\W(k)$-modules to modules of higher torsion. To do so, we use divided powers, which are arguably the most efficient "elementary" algebraic tool at disposal. \\

Before beginning, let us mention here the work of Kaledin- notably in the recent papers \cite{K} and \cite{K2}. His constructions are very close to what we produce here, though expressed in a  different language. We did not have time to explore the connections so far and it would be interesting to do so in the future.\\

Recall that we denote by $$\tau: k \lra \W(k)$$ the  Teichm\"uller representative (we set $\tau(0)=0$).  The map  $$\tau_{\vert k^\times}: k^\times \lra \W(k)^\times$$ is the unique multiplicative section of the quotient map $$\W(k)^\times \lra k^\times.$$

\noindent Let $n,m$ be  positive integers. Take $A$ to be $\W_{n+m}(k)$, the truncated Witt vectors of size $n+m$. By Lemma \ref{Teichpoly}, the formula $$\overline \tau_n: A/p^mA \lra A/p^{n+m}A=\W_{n+m}(k),$$ $$ x \mapsto x^{p^n},$$  defines a (multiplicative) polynomial law over $\Z$, homogeneous of degree $p^{n}$. It thus induces a group homomorphism $$ \Gamma^{p^{n}}_\Z(\W_m(k)) \lra  \W_{n+m}(k),$$ which we  denote by $T'_n$. Since  the Teichm\"uller representatives generate  $\W_{n+m}(k)$ additively, the map $T'_n$ is surjective. But obviously, the  polynomial law $\overline \tau_n$ might as well be considered as a polynomial law over $\W(k)$, hence giving rise to a $\W(k)$-linear homomorphism  $$T_n:    \Gamma^{p^{n}}_{\W(k)}(\W_m(k))  \lra \W_{n+m}(k).$$ 

\begin{lem}
	The map $T_n$ is an isomorphism.
\end{lem}

\begin{dem}
It is clearly surjective. But $\Gamma^{p^{n}}_{\W(k)}(\W_m(k)) $ is generated by $[1]_{p^n}$, as a $\W(k)$-module, and is killed by $p^{n+m}$, by Lemma \ref{torsionbornee}. It is thus a $\W_{n+m}(k)$-module, generated by one element. The claim follows.\end{dem}

Thus, the map $T'_n$ factors as $$T'_n:  \Gamma^{p^{n}}_\Z(\W_m(k)) \lra   \Gamma^{p^{n}}_{\W(k)}(\W_m(k))  \stackrel {T_n} \lra \W_{n+m}(k).$$ One can then infer a description of the kernel of $T'_n$, which in turn provides a rather simple, seemingly new,  natural recursive  definition of $\W_m(k)$ by generators and relations. It does not involve any intricate computation- assuming, of course, some familiarity with the spirit of divided powers.

\begin{prop}\label{newwitt}
Assume that $n=1$ in what precedes.\\
 Then the kernel of the natural surjection $$T'_1:  \Gamma^{p}_\Z(\W_m(k)) \lra \W_{m+1}(k)$$  is generated, as an Abelian group, by elements of the form $$[x]_1 [y]_{p-1}-[xy^{p-1}]_1 [1]_{p-1}, $$ with $x,y \in \W_m(k) .$
\end{prop}

\begin{dem}
We try to offer an  algorithmic proof. The reader who is not familiar with divided powers is advised to assume  $p=2$, to begin with.\\
First of all, it is not hard to see that these elements are in the kernel, since, on 'impure' symbols, we have $$ T'_1([x]_1 [y]_{p-1})=  p! X^i Y^j  \in \W_{m+1}(k),$$ for all $x,y \in \W_m(k)$, where $X,Y\in \W_{m+1}(k)$ are arbitrary lifts of $x$ and $y$, respectively. Denote by $I \subset \Gamma^{p}_\Z(\W_m(k))$ the Abelian group spanned the  elements $$[x]_1 [y]_{p-1}-[xy^{p-1}]_1 [1]_{p-1} .$$ \\
Pick an element $X \in \Ker(T'_1).$ A little computation shows that, modulo $I$, every element of $\Gamma^{p}_\Z(\W_m(k))$ is congruent to an element of the shape $$[a]_p+ [b]_1 [1]_{p-1}, $$ with $a ,b \in \W_{m}(k)$.  We can thus assume that $$X=[a]_p+ [b]_1 [1]_{p-1}. $$ Denote by $A\in \W_{m+1}(k)$ and $B\in \W_{m+1}(k)$ arbitrary lifts of $A$ and $B$, respectively. The relation $T'_1(X)=0$ translates as $$ A^p +p! B=0 \in \W_{m+1}(k),$$ which implies that $a$ is divisible by $p$, say $a=pa'$, for  $a'\in \W_{m}(k)$. From the fact that, modulo $I$, we have $$[a]_p= p^p [a']_p= {\frac {p^{p-1}} {(p-1)!}} [a']_1[a']_{p-1} \equiv {\frac {p^{p-1}} {(p-1)!}} [a'^p]_1[1]_{p-1}  \in \Gamma^{p}_\Z(\W_m(k)),$$  we are thus reduced to case  $a=0$. But then $pB=0 \in \W_{m+1}(k), $ implying $b=0 \in \W_{m}(k) .$  Hence, $I$ indeed equals  $\Ker(T'_1).$


\end{dem}

From now on, if $M$ is a torsion $\W(k)$-module, we shall put $$\Gamma^n(M):=\Gamma_{\W(k)}^n(M)$$ and $$\Sym^n(M)=\Sym_{\W(k)}^n(M) .$$ 

Note that these are polynomial functors, in  the category of torsion  $\W(k)$-modules. 

\begin{rem}
The preceding discussion shows that $\Gamma^{p^n}(\W_m(k))$, as a $\W(k)$-module, is generated by $[1]_{p^n}$ and is canonically isomorphic to $\W_{m+n}(k)$. We are now going to make this statement (a bit) more precise.
\end{rem}

Let $n,m$ be  positive integers.\\
Let $L$ be  a $\W_m(k)$-module which is free of rank one. Seeing it as a $\W(k)$-module,  we put $$\W_{n+m}(L):=\Gamma^{p^n}(L);$$  it is is a free $\W_{n+m}(k)$-module of rank one, whose construction is functorial in $L$. It comes equipped with the Teichm\"uller-like map (which is in fact  a polynomial law)$$\overline \tau_n: L \lra \W_{n+m}(L),$$ $$ v \mapsto [v]_{p^n}.$$  Note that, if $m=1$ and $L=k$, then  $\W_{n+1}(L)=\W_{n+1}(k)$, and $\overline \tau_n(x)=\tau(x^{p^n})$, as noted before.

\begin{lem}\label{gammadroite}

Let $i=p^n j$ be a positive integer, with $j$ prime to $p$.  Then the  formula $$ L \lra \W_{n+m} (L^{\otimes j}),$$ $$ v \mapsto [v ^{\otimes j}]_{p^n},$$ defines a polynomial law, which is homogeneous, of degree $i$. The induced $\W(k)$-linear map $$\phi: \Gamma^i(L) \lra  \W_{n+m} (L^{\otimes j})$$ is an isomorphism.
\end{lem}

\begin{dem}
Only the fact that $\phi$ is an isomorphism has, perhaps, to be checked. We may assume that $L=\W_m(k)$.   By lemma \ref{torsionbornee}, the $\W(k)$-module $\Gamma^n(\W_m(k))$, which is obviously generated by $[1]_n$, is of $p^{m+n}$-torsion, hence a $\W_{m+n}(k)$-module generated by one element. The map $\phi$ is obviously surjective, with target a free $\W_{m+n}(k)$-module. It is thus an isomorphism.

 \end{dem}

\begin{lem}
Let $n,m$  be  positive integers. 
Let $M$ be a  $\W_m(k)$-module.  Pick an element $x$, of order  $p^m$. Then the symbol $$ [x]_n \in \Gamma^n (M)$$ has order $p^{v_p(n)+m}$.
\end{lem}

\begin{dem}
By Lemma \ref{torsionbornee}, the symbol in question has order $\leq p^{v_p(n)+m}$. Now, pick a $\W_m(k)$-linear map $$f: M \lra \W_m(k),$$ sending $v$ to $1$.  By functoriality, it induces a $\W_{n+m}(k)$-linear map $$F: \Gamma^n (M)\lra \Gamma^n( \W_m(k)),$$ mapping $[v]_{n}$ to $[1]_n$. By Lemma \ref{gammadroite}, and by the fact that $1 \in \W_{m+n}(k)$ has order $p^{m+n}$, we  know that  $[1]_n$ has order  $ p^{v_p(n)+m}$. The claim follows.
\end{dem}

	\begin{lem}\label{symbolsgen0}
		Let $V$ be a $k$-vector space. Let $n$ be a positive integer,  lesser or equal to the cardinality of $k$. Then the symbols $ [v]_n $, for $v \in V$, generate $\Gamma^n(V)/p$ (as a $k$-vector space).
	\end{lem}

	\begin{dem}
	We can assume that $V$ is finite-dimensional.
	By a straightforward induction on the dimension $d \geq 2$ of $V$, it is enough to show that the natural map $$\bigoplus_{H \in \P(V^*)} \Gamma^{n}(H) /p \lra \Gamma^{n} (V) /p,$$ given by the sum of the inclusions $\Gamma^{n}(H) /p \lra \Gamma^{n} (V)/p,$ for all hyperplanes $H \subset V$, is surjective. Dually, letting $W:=V^*,$ we have to show that the natural map 
	 $$\Sym^{n}(W) \lra \oplus_{L \in \P(W)} \Sym^{n}(W/L) ,$$ given as the sum of the quotient maps, is injective. But, choosing a $k$-basis of $W$, an element of $\Sym^{n}(W)$ is just a homogeneous polynomial of degree $n$ in $d$ variables. The fact that it dies in  $\Sym^{n}(W/L)$ is equivalent to asking that it is divisible by $v$, where $v \in L$ is a nonzero vector. The statement now follows, since  $\P(W)$ has cardinality at least $\vert k \vert +1 \geq n+1$, and since a homogeneous polynomial of degree $n$, which is divisible by $n+1$ two by two non proportional linear factors, has to be zero.

	\end{dem}

\begin{lem}\label{symbolsgen}
	Let $M$ be a torsion $\W(k)$-module. Let $n$ be a positive integer, lesser or equal to the cardinality of $k$. Then the symbols $ [x]_n $, for $x \in M$, generate $\Gamma^n(M)$ (as a $\W(k)$-module).
\end{lem}

\begin{dem}
We can assume that $M$ is of finite type.
Consider the filtration $$\Gamma^n(M) \supset p\Gamma^n(M) \supset p^2 \Gamma^n(M) \supset \ldots \supset \{0\},$$ and note that the successive quotients are all quotients of $\Gamma^n(M)/p \simeq \Gamma^n_k(M/p)$. Apply induction, using Lemma \ref{symbolsgen0}, to get the result.
 
\end{dem}

We conclude this section by a concrete description of  divided power modules, using a basis. We first fix some useful notation.
\begin{defi}(Weighted partitions.)
Let  $n$ be a positive integer.  A  partition of $n$ is a decomposition  $$ A: (a_1+ \ldots +a_d= n)$$  of $n$ into a sum of $d$ nonnegative integers. The partition $A$ is said to be weighted, if we are given the extra data of  a $d$-tuple $$w:=(w_1, \ldots, w_d)$$ of  nonnegative integers, such that, for all $i=1 \ldots, d$,  $w_i$ is zero if $a_i$ is zero. \\
The (weighted) partition $(A,w)$ is said to be proper if all $a_i$'s are $\leq n-1$, or  Dirac  otherwise.
\end{defi}
\begin{defi}\label{notapratique}
	Let $M$ be a finite $\W_m(k)$-module. Choose a decomposition $$M = \bigoplus_{i=1}^d \W_{m-w_i}(k) e_i,$$ where $e_i \in M$ has order $p^{m-w_i}$  (the $w_i$'s are $ \leq m$, and uniquely determined by $M$). Put $$w=(w_1 , \ldots, w_d).$$ Denote by  $$M^\vee = \bigoplus_{i=1}^d \W_{m-w_i}(k) e_i^\vee$$  the dual decomposition, with $$<e_i,e_j^\vee>=p^{w_i} \delta_{i,j}. $$  Let $n$ be a positive integer. For each   partition  $$ A: (a_1+ \ldots +a_d= n),$$  put  $$ \tilde W(A,w):= \max\{v_p(n)-v_p(a_i)+ w_i, i=1 \ldots d\}$$  and  $$ [e]_A:=[e_1]_{a_1} \ldots [e_d]_{a_d} \in \Gamma^{n}(M).$$ 

\end{defi}

\begin{prop}\label{isobase}
	Let $M$ be a finite $\W_m(k)$-module. Let $n \geq 0$ be an integer.\\ We use the notation of Definition \ref{notapratique}.\\  Then the order of $[e]_A$ in the $\W_{m+v_p(n)}(k)$-module $\Gamma^n(M)$  is $p^{m+v_p(n)-\tilde W(A,w)}$, and  there exists a natural isomorphism of $\W_{m+v_p(n))}(k)$-modules $$\bigoplus_{A}   \W_{m+v_p(n)-\tilde W(A,w)}(k) \stackrel \sim \lra \Gamma^n(M),$$ $$1_A \mapsto [e]_A,$$ 
	where the sum is taken over all partitions $A$ of $n$, of size $d$.

\end{prop}

\begin{dem}
	The first statement follows  from Proposition \ref{polyfunctor} and Lemma \ref{gammadroite}, which implies that  $$\bigotimes_{i=1}^d \Gamma^{a_i}(\W_{m-w_i}(k))$$ is canonically isomorphic to $ \W_{m+v_p(n)- \tilde v(A,w)}(k)$.

\end{dem}

 \subsection{An alternate description of $\Gamma^p$ for vector spaces.}\label{alternategamma}

Here $k= \F_p$. Assume that $V=M \otimes_\Z \F_p,$ for $M$ a free $\Z$-module of finite rank. One readily checks that the map $$C: M \times M \lra \Sym_\Z ^p(M),$$ $$ (x,y) \mapsto \frac {(x+y)^p-x^p-y^p} p,$$ is a symmetric $2$-cocycle, for the trivial action of $M$ on  $ \Sym_\Z ^p(M)$. Indeed, this can be checked after extending scalars to $\Q$, where it is obvious: $c$ is then a trivial cocycle by definition! Reducing mod $p$, we obtain a symmetric cocycle $$c: V \times V \lra  \Sym_{k} ^p(V),$$ in fact given by $$ c(x,y)= \sum_1^{p-1}  \frac {(-1)^{i-1}} i  x^i y^{p-i}.$$ This cocycle defines an Abelian extension of $V$ by $\Sym^p_k(V)$. We leave it to the reader, as an instructive exercise, to check that this extension is canonically isomorphic to $\Gamma^p(V)$.





\section{The Frobenius and the Verschiebung.}\label{frob}

Recall that $k$ is a perfect field of characteristic $p$.\\
Let $A$ be a commutative ring of characteristic $p$. Denote by $$\frob_A :A \lra A,$$ $$x \mapsto x^p,$$ the Frobenius endomorphism of $A$. For any $A$-module $M$, put $$M^{(1)}:=M \otimes_A A,$$ where the tensor product is taken with respect to $\frob_A$. This notation is obviously coherent with the one used before.\\
	Moreover, if $B/A$ is a commutative algebra, we have a canonical isomorphism $$M^{(1)} \otimes_A B \stackrel \sim \lra  (M \otimes_A B)^{(1)}.$$ In other words, forming the twist by Frobenius commutes with extensions of commutative rings of characteristic $p$.\\
	Now, let $V$ be a $k$-vector space. By what precedes, the formula $$V \lra  V^{( 1)},$$ $$ v \mapsto v^{( 1)} := v\otimes 1,$$ actually defines a polynomial law, homogeneous of degree $p$. We shall refer to this law as the Frobenius law $\Frob_V$. It can be viewed as a morphism of affine $k$-spaces $$\A_k(V) \lra \A_k(V^{(1)}).$$

	Note that the Frobenius law exists only for $k$-vector spaces. For $n \geq 2$, an arbitrary (say, finite free) $\W_n(k)$-module does not come naturally equipped  with such a law. 

The next proposition contains the definition of the Frobenius and of the Verschiebung, borrowed from the theory of commutative group schemes in characteristic $p$.

\begin{prop}
 Let $V$ be a finite-dimensional $k$-vector space and let $n \geq 1$ be an integer. \\
 
 Then the formula $$V \lra \Gamma^n( V^{( 1)}),$$  $$v \mapsto  [v ^{(1)}]_n,$$ is a polynomial law of degree $np$, thus defining a  $\W(k)$-linear map $$\Frob: \Gamma^{np} (V) \lra  \Gamma^n( V^{( 1)}), $$ the Frobenius homomorphism (for divided powers). \\
 The polynomial law  $$V \lra \Gamma^{np}( V),$$  $$v \mapsto  p [v]_{pn}$$ canonically factors through the Frobenius law $V \lra V^{( 1)}$. The resulting polynomial law $$ V^{( 1)} \lra \Gamma^{np}( V)$$ is homogeneous of degree $n$, yielding a $\W(k)$-linear map $$\Ver: \Gamma^n( V^{( 1)} ) \lra  \Gamma^{np}( V),$$ $$ [v ^{(1)}]_n \lra p[v]_{pn},$$ the Verschiebung homomorphism (for divided powers).
 
\end{prop}

\begin{dem}

The first statement (defining the Frobenius map for divided powers) follows from the definition of the Frobenius law. For the second one, pick a basis $e_1 , \ldots, e_d$ of $V$. On the one hand, the Frobenius $V \lra  V^{(1)}$ then becomes the law $$k^d \lra k^d,$$ $$ (X_1,\ldots, X_d) \mapsto  (X_1^p,\ldots, X_d^p).$$ On the other hand, The polynomial law  (of $\W(k)$-modules)  $$V \lra \Gamma^{np}( V),$$  $$v \mapsto  p [v]_{pn}$$ then becomes the law $$ k ^d \lra  \Gamma^{np}( k^d),$$ $$(X_1,\ldots,X_d) \mapsto  p[X_1 e_1+ \ldots + X_d e_d]_{pn} $$ $$= \sum_{a_1 +\ldots + a_d = pn}  X_1 ^{a_1}\ldots   X_d ^{a_d} p [e_1]_{a_1} \ldots [e_d]_{a_d}$$  $$= \sum_{a_1 +\ldots + a_d = n}  X_1 ^{pa_1}\ldots   X_d ^{pa_d} p [e_1]_{pa_1} \ldots [e_d]_{pa_d},$$ where the first (resp. second) sum is taken over all decompositions of $pn$ (resp. of $n$) into the sum of $d$ nonnegative integers. Indeed, the symbols $[e_i]_a$, for $a$ not divisible by $p$, are of additive order $p$, hence all terms $ p [e_1]_{a_1} \ldots [e_d]_{a_d}$ vanish, as soon as one of the $a_i$'s is not divisible by $p$. The second part of the lemma, yielding the definition of the Verschiebung morphism for divided powers, is now obvious.

\end{dem}

\begin{lem}\label{formulaFV}
Let $V$ be a $k$-vector space. Let $a_1, \ldots, a_d,n$ be nonnegative integers, satistying $a_1+ \ldots+ a_d=np$. For $v_1, \ldots, v_d \in V$, the Frobenius $$\Frob: \Gamma^{np} (V) \lra  \Gamma^n( V^{( 1)}) $$ satisfies $$\Frob([v_1]_{a_1} \ldots [v_d]_{a_d} ) = 0,$$ if one of the $a_i$'s is not divisible by $p$. If all $a_i$'s are divisible by $p$, says $a_i=pb_i$, then $$\Frob([v_1]_{a_1} \ldots [v_d]_{a_d} ) = [v_1 ^{(1)}]_{b_1} \ldots [v_d ^{(1)}]_{b_d}.$$
Dually, the Verschiebung  $$\Ver: \Gamma^n( V^{( 1)} ) \lra  \Gamma^{np}( V)$$ satisfies   $$\Ver([v_1 ^{(1)}]_{a_1} \ldots [v_d^{(1)}]_{a_d} ) =p [v_1]_{pa_1} \ldots [v_d]_{pa_d}.$$

\end{lem}

\begin{dem}
We work in the polynomial ring $\W(k)[X_1, \ldots X_d]$. The relation $$ \Frob([X_1 v_1 + \ldots + X_d v_d]_{np})=[X_1^p( v_1^{(1)}) + \ldots + X_d^p ( v_d ^{(1)})]_{n}$$ holds by definition. But $$ [X_1 v_1 + \ldots + X_d v_d]_{np}=\Sigma_{a_1+ \ldots+ a_d=np}( X_1^{a_1} \ldots  X_d^{a_d} [v_1]_{a_1} \ldots  [v_d]_{a_d})$$ and  $$ [X_1^p( v_1 ^{(1)}) + \ldots + X_d^p( v_d ^{(1)})]_{n}=\Sigma_{b_1+ \ldots+ b_d=n}( X_1^{pb_1} \ldots  X_d^{pb_d} [v_1 ^{(1)}]_{b_1} \ldots  [v_d ^{(1)}]_{b_d}),$$ so that the first assertion follows by identifying the coefficients of the monomials occuring in those expansions. The proof for the Verschiebung is similar.

\end{dem}
\begin{coro}\label{Kerfrob}
   For any $s \geq 1,$ the kernel of  $$\Frob^s:\Gamma^{np^s}(V) \lra \Gamma^{n}(V^{(s)})$$ coincides with the $p^s$-torsion of $\Gamma^{np^s}(V)$.
\end{coro}

\begin{dem} This  follows from  Proposition \ref{isobase} and  Lemma \ref{formulaFV}.

\end{dem}

\begin{prop}\label{Frobsurj}
	Let $V$ be a $k$-vector space. Let $n,s \geq 1$ be  integers. Then the Frobenius  $$\Frob^s: \Gamma^{np^s} (V) \lra  \Gamma^n( V^{( s)}), $$ $$[v]_{np^s} \lra [v ^{(s)}]_n$$ is surjective. We have an exact sequence $$0 \lra \Gamma^{np^s} (V)[p^s]  \lra \Gamma^{np^s} (V) \stackrel {\Frob^s} \lra  \Gamma^n( V^{( s)}) \lra 0. $$  The Verschiebung  $$\Ver^s: \Gamma^n( V^{( s)} ) \lra  \Gamma^{np^s}( V),$$  $$ [v ^{(s)}]_n \lra p^s[v]_{np^s},$$ is injective. We have an exact sequence $$0 \lra  \Gamma^n( V^{( s)} ) \stackrel {\Ver ^s}\lra  \Gamma^{np^s}( V) \lra \Gamma^{np^s}( V) / p^s\lra 0. $$
\end{prop}
\begin{dem}
 We may assume that $V$ is finite-dimensional, with basis $e_1, \ldots, e_d$. Let $k'/k$ be a extension of perfect fields, with $k'$ infinite. Since the formation of divided powers commutes to the (faithfully flat!) base-change $\W(k')/ \W(k)$, we can assume that $k$ is infinite. In this case, Lemma \ref{symbolsgen} ensures that pure symbols are additive generators of divided powers. The surjectivity of $\Frob^s$  then directly follows from the description given in Lemma \ref{formulaFV}. That its kernel is the $p^s$-torsion is the content of Corollary \ref{Kerfrob}. By Lemma \ref{formulaFV}, it is clear that the image of $\Ver^s$ is $p^s\Gamma^{np^s}(V)$. It follows from the same Lemma, combined with Proposition \ref{isobase}, that $\Ver^s$ is injective.

\end{dem}

\begin{coro}\label{Frobhom}
   The Frobenius $$\Frob: \Gamma(V) \lra \Gamma(V^{(1)})$$ is a surjective homomorphism of $\W(k)$-algebras, with kernel $ \Gamma^+(V)[p]$.
\end{coro}

\begin{dem}
    This is now obvious.
\end{dem}

\section{Divided powers and Pontryagin duality.}\label{omega}
\subsection{Duality.}
Recall that, for every $k$-vector space $V$, we have canonical isomorphisms $$\Gamma^p_k(V) \simeq (V^{\otimes p})^{\mathcal S_p} \simeq \Sym_k^p(V^*)^*.$$

When working over a field  of characteristic zero, it is common (though somewhat misleading) to identify $\Sym_k^p(V^*)^*$ and $\Sym_k^p(V)$,  using what is called the 'symmetrizing operator'. Equivalently, in characteristic zero, the map $$ \Gamma_k^p(V) \lra \Sym_k^p(V),$$ $$ [v]_p \mapsto v^p,$$ is an isomorphism. It is of course far to be so in our context, where the perfect field $k$ has characteristic $p$. In other terms, the functor $\Gamma_k^p$ does not commute with duality of vector spaces. However, we will see that, for any $m \geq 1$, the functor $\Gamma^p=\Gamma_{\W_{m+1}(k)}^p$ \textit{does} commute with duality for a \textit{free} $\W_m(k)$-module $M$, in the sense that $\Gamma^p(M^\vee)$ and $\Gamma^p(M)^\vee$ are canonically isomorphic, as $\W_{m+1}(k)$-modules. \\
This phenomenon does unfortunately not extend to higher divided powers: the functors $\Gamma^{p^n}$, for $n \geq 2$, behave very badly with duality- except in dimension two.  This is a rather subtle fact, linked to intricate computations of  $p$-adic valuations	of scary  multinomial coefficients. To bypass this difficulty, we can choose to apply  $\Gamma^p$ $n$ times in a row, instead of applying the functor $\Gamma^{p^n}$ just once. In doing so, we lift the ($p$-adic valuation of the) torsion of the modules by one at each step. In some sense, this choice is justified by the basic computational fact that, for an ideal in a $\Z_p$-algebra, a divided power structure is uniquely determined by the operation $\gamma_p$. Note that, for non free $\W_m(k)$-modules, we  will have to take a \textit{quotient} of the functor $\Gamma^p$ anyway, in order to respect duality. \\
In what follows, we explore these two points of view: applying $\Gamma^{p^n}$, or applying $\Gamma^p$ $n$ times in a row. Forming their 'correct' quotients   (which commutes to Pontryagin duality) will give rise to medium and big Omega powers, respectively. Medium Omega powers will turn out to be a \textit{direct factor} of big Omega powers.\\

From now on, $m$  fixed positive integer.  The Pontryagin dual of a $\W_m(k)$-module $M$ shall be viewed as $$M^\vee=\Hom_{\W_m(k)}(M,\W_m(k)).$$

Let $M$ be a  $\W_m(k)$-module.  Then
the duality law $$\Delta : \overline M \times \overline {M^\vee} \lra \W_{m}(k)$$  of Definition \ref{duality} canonically factors through the projection $$\pi:  \W_{m+1}(k) \lra  \W_{m}(k).$$ Explicitly, the pairing $$\tilde \Delta_{ p} : \overline M \times \overline {M^\vee} \lra \W_{m+1}(k)$$ $$(m, \phi) \mapsto \overline \tau_1(\phi(m)) $$ is a polynomial law, bihomogeneous of bidegree $(p,p)$,  satisfying $$ \Delta_{p}=\pi \circ \tilde \Delta_{ p}. $$ It strongly depends on $m$. The integer $m$ being fixed in this section,  we shall abusively write  $\Delta$, or even $<.,.>$ for $\tilde \Delta_{ p}$. By the universal property of divided powers, it corresponds to a pairing of $\W_{m+1}(k)$-modules $$\Gamma^{p}_{\W_{m+1}(k)}(M) \times \Gamma^{p}_{\W_{m+1}(k)}(M^{\vee}) \lra \W_{m+1}(k),  $$  
given, on the level of pure symbols, by the formula $$<[\phi]_{ p},[v]_{ p}>=\overline \tau_1(\phi(v)) \in \W_{m+1}(k).$$
Similarly, for any nonnegative integer $n \geq 2,$ we get  a pairing of $\W_{m+n}(k)$-modules $$\Gamma^{p^n}_{\W_{m+n}(k)}(M) \times \Gamma^{p^n}_{\W_{m+n}(k)}(M^{\vee}) \lra \W_{m+n}(k),  $$  
$$<[\phi]_{ p^n},[v]_{ p^n}>=\overline \tau_n(\phi(v)) \in \W_{m+n}(k).$$ We shall denote this pairing by $\tilde \Delta_{ p^n}$, or again simply by $\Delta$ if the context is clear.

	\begin{lem}\label{deltabase}
	Let $n$ be a nonnegative integer.\\
	Let $M$ be a finite $\W_m(k)$-module.  Using the notation of Definition \ref{notapratique}, we  have a commutative diagram $$ \xymatrix {   \Gamma^{ p^n}(M) \times  \Gamma^{ p^n}(M^\vee)  \ar[d] ^\wr \ar[r]^{<.,.>}   &  \W_{m+n}(k) \ar@{=}[d] \\ ( \bigoplus_A  \W_{N(A)}(k)[e]_A) \times  (\bigoplus_B \W_{N(B)}(k) [e^\vee]_B\ar[r])  & \W_{m+n}(k) , }$$ where the vertical map on the left is the product of the isomorphisms given by Lemma \ref{isobase}, and the lower horizontal map is the pairing given by $$([e]_A ,[e^\vee]_B) \mapsto p^{(\sum_1^d w_i a_i)}\binom { p^n}{a_1,a_2, \ldots, a_d} \in \W_{m+n}(k),$$ if $A=B=(a_1, \ldots, a_d)$, or by $$([e]_A ,[e^\vee]_B) \mapsto  0$$  if $A \neq B$.
	\end{lem}
	
	\begin{dem}
		We work over the polynomial ring  $\W_{n+m}(k)[X_i, Y_i$, $i=1 \ldots d]$. By definition, $$<[X_1 e_1 +
		\ldots +X_d e_d]_{ p^n},[Y_1 e_1^* +
		\ldots +Y_d e_d^*]_{ p^n}>=(p^{w_1}X_1Y_1 + \ldots + p^{w_d}X_d Y_d)^{ p^n} .$$  Developping the lefthand side, we get that the coefficient of $X_1^{a_1} \ldots X_d^{a_d} Y_1^{b_1} \ldots Y_d^{b_d}$ is $<[e_1]_{a_1} \ldots [e_d]_{a_d},[e_1^*]_{b_1} \ldots [e_d^*]_{b_d}>$, whenever $a_i$ and $b_i$ are nonnegative integers such that  $a_1+ \ldots+ a_d= b_1+ \ldots+ b_d= p^n.$ Developping the righthand side, and identifying the coefficients, yields the result.
		
	\end{dem}

\begin{defi}
let $d$ be a positive integer. Pick a weighted  partition $(A,w)$ of $p^n$, of size $d$.  Put $$W( (A,w)):= \min \{m+n, v_p(\binom { p^n}{A}) +\sum_1^d w_i a_i\} .$$ 
\end{defi}

\begin{rem}\label{ptiterem}

Note that $$W((A,w)) \geq \tilde W((A,w)) . $$ If $n=1$ , then equality holds if, and only if, either $w=0$, or  $w_i$ is nonzero for a single index $i$, for which $a_i=1$.\\
If $d=2$ and $w=0$, equality holds for any $n$. This holds because, in this case, $v_p(\binom { p^n}{a_1,a_2})$ is precisely $n-v_p(a_1)=n-v_p(a_2)$.\\
In the context of Lemma \ref{isobase}, we know that the order of $[e]_A \in \Gamma^{p^n}(M)$ is  $p^{m+n- \tilde W((A,w))}$. We will now show that the order of (the class of) $[e]_A$ in $\Gamma^{p^n}(M)/ \Ker(\Delta)$ is  $p^{m+n-  W((A,w))}$.
\end{rem}

	\begin{lem}\label{kerdeltabase}
	Let $M$ be a finite $\W_m(k)$-module. Let $n$ be a nonnegative integer. \\ We use the notation of Definition \ref{notapratique}.\\
We  have a canonical isomorphism $$   \Gamma^{ p^n}(M)  / \Ker(\Delta) \stackrel \sim \lra \bigoplus_A \W_{m+n-W(A,w)}(k),$$ $$ [e]_A \mapsto 1_A,$$  where the direct sum is taken over all  partitions $A$ of $p^n$, of size $d$. \\It fits into a commutative diagram $$\xymatrix{ \Gamma^{ p^n}(M) \ar[r] \ar[d] &  \bigoplus_A \W_{m+n-\tilde W(A,w)}(k) \ar[d] \\  \Gamma^{ p^n}(M)  / \Ker(\Delta)  \ar[r] &  \bigoplus_A \W_{m+n-W(A,w)}(k)}  $$
	\end{lem}

	\begin{dem}
		Obvious from Lemma \ref{deltabase}.
	\end{dem}
	
	It is natural to ask whether the pairing $\Delta$ is non-degenerate; in other words, to ask whether $$ \Delta:  \Gamma^{ p^n}(M) \lra \Gamma^{ p^n} (M^\vee)^\vee$$ is an isomorphism of $\W_{m+n}(k)$-modules. We can now answer this question.
	
	\begin{prop}\label{gammaomega} 	Let $M$ be a finite $\W_m(k)$-module. 
	The  pairing $$ \Delta:  \Gamma^{ p}(M) \times   \Gamma^{ p}(M^\vee) \lra \W_{m+1}(k) $$ is  perfect if and only if $M$ is a free $\W_m(k)$-module.\\
	If $M$ is a free  $\W_m(k)$-module of rank at most two, then we have more: the pairing  $$ \Delta:  \Gamma^{ p^n}(M) \times   \Gamma^{ p^n}(M^\vee) \lra \W_{m+n}(k) $$ is  perfect, for any $n \geq 1$.
	\end{prop}
	
	\begin{dem}
We use Lemma \ref{kerdeltabase}. In the first case, we have to see that $\tilde W(A,w)=W(A,w)$ for all $A$, if and only if $w=0$. In the second case, we have to see that $\tilde W(A,0)=W(A,0)$ if $d=2$. This is clear, at the light of Remark \ref{ptiterem}.
		
	\end{dem}

\subsection{Medium and big Omega powers.}

We now present a construction of crucial importance in this paper: the so-called (big and  medium) Omega powers.  The choice of the name "Omega" is naturally inspired from topology: it is close to an algebraic loop space. 
\begin{defi}
	Let $M$ be a $\W_m(k)$-module. Let $n$ be a nonnegative integer. We define a $\W_{m+n}(k)$-module $\Gamma^{(n)}(M)$ by  $\Gamma^{(0)}(M)=M$, and by the recursive formula $$ \Gamma^{(n)}(M):=\Gamma^p(\Gamma^{(n-1)}(M)).$$
	\end{defi}
	
By Proposition \ref{gammaomega}, we now that $\Gamma^{(1)}$ commutes to Pontryagin duality, for \textit{free} $\W_m(k)$-modules only. Note that, if $M$ is a free   $\W_m(k)$-module, $\Gamma^{(1)}(M)$ is a free  $\W_{m+1}(k)$-module if and only if $M$ has rank one. The functor $\Gamma^{(2)}$ will thus never commute with duality, except for free modules of rank lesser than one. To define our  Omega functor (medium and big), we are thus naturally led to mod out the kernel of Pontryagin duality.

\begin{defi}(medium and big Omega functors.)\label{omega} Let $M$ be a (non necessarily finite) $\W_m(k)$-module. Let $n$ be a positive integer. We put
	 $$  \Omega^n_m(M):= \Gamma^{ p^n}_{\W_{m+n}(k)}(M) /\Ker(\Delta).$$ 
	
	It is the $n$-th medium Omega power of the  $\W_m(k)$-module $M$. It is a $\W_{m+n}(k)$-module.  
	We put $ \overline \Omega_m^{0}(M)=\Omega_m^0(M)=M$.\\ 
	We recursively define $$\overline \Omega_m^{n}(M):=  \Omega_{m+n-1}(\overline \Omega_m^{n-1}(M)) ;$$ it is a $\W_{m+n}(k)$-module as well.  	It is the $n$-th big Omega power of the  $\W_m(k)$-module $M$.\\
	We shall denote $ \Omega^n_m(M)$ (resp. $ \overline \Omega^n_m(M)$) simply by  $  \Omega^n(M)$ (resp. $ \overline \Omega^n(M)$), if the dependence in $m$ is clear. \\
	 For $x \in M$, if this creates no confusion,  we denote by $$(x)_n \in \Omega^n(M)$$ the class of the pure symbol $[x]_{p^n} \in  \Gamma^{ p^n}_{\W_{m+n}(k)}(M) .$
	
\end{defi}

\begin{rem}

We clearly have $\overline \Omega^1=\Omega^1$. For $n \geq 2,$ we will see in a moment that $\Omega^n$ appears canonically as a direct factor of $\overline \Omega^n$, in a duality-preserving way.
\end{rem}

\begin{rem}
    	Let $M$ be a $\W_m(k)$-module.  Let $n$ be a nonnegative integer. We have a canonical surjection $$  \Gamma^{(n)}(M) \lra \overline \Omega^n(M).$$
    
\end{rem}

\begin{rem}
If $M=\W_m(k)$ is a free $\W_m(k)$-module of rank one,   then $\overline \Omega^n (M) =  \Omega^n (M) = \Gamma^{ p^n}_{\W(k)}(M) $  is a free $\W_{m+n}(k)$-module of rank one.
\end{rem}
	
	\begin{rem}
	    The associations $$M \mapsto \overline \Omega^n(M)$$ and  $$M \mapsto  \Omega^n(M)$$  are  functors, from the category of  $\W_m(k)$-modules to that of  $\W_{m+n}(k)$-modules.  For $n=1$, they  behave slightly like the $p$-th symmetric power functor $\Sym^p$. However,  this analogy is quite bad, as we have  seen in Section \ref{PDVSSIM}.
	\end{rem}

The functor $\Omega^1$ is polynomial, with respect to the Tense Product. 

\begin{prop}\label{polytense}
	Let $M,N$ be $\W_m(k)$-modules. Then, we have a canonical isomorphism of $\W_{m+1}(k)$ -modules $$\overline \Omega^1(M \bigoplus N) \simeq \overline \Omega^1(M) \bigoplus \overline \Omega^1(N) \bigoplus \bigoplus_{i,j} \left(\overline \Sym_{m}^i(M) \overline \bigotimes_{m} \overline \Sym_{m}^j(M)  \right),$$ where the direct sum is taken over all proper partitions $i+j=p$.
	\end{prop}
	
	\begin{dem}
 By Proposition \ref{isobase}, we have a natural isomorphism $$\Gamma^p(M \bigoplus N) \simeq \Gamma^p(M) \bigoplus\Gamma^p(N) \bigoplus \bigoplus_{i,j} \left( \Gamma^i(M) \bigotimes \Gamma^j(N)  \right),$$ where the direct sum is taken over all proper partitions $i+j=p$. Let us write the dual decomposition $$\Gamma^p(M^\vee \bigoplus N^\vee) \simeq \Gamma^p(M^\vee) \bigoplus\Gamma^p(N^\vee) \bigoplus \bigoplus_{i,j} \left( \Gamma^i(M^\vee) \bigotimes \Gamma^j(N^\vee)  \right).$$ These are compatible with the duality pairing (with values in $\W_{m+1}(k)$). Notably, for $x \in M$, $y \in N$, $\phi \in M^\vee$ and $ \psi \in N^\vee$, we have $$<[x]_i \otimes [y]_j,[\phi]_i \otimes [\psi]_j> = {p \choose {i,j}} \phi(x)^i \psi(y)^j \in \W_{m+1}(k).$$ Since $i$ and $j$ are $\leq p-1$, the binomial coefficient in this formula has $p$-adic valuation one,  and the claim follows.
	\end{dem}

	\subsection{Medium Omega powers as a direct factor of big Omega powers.}
	We now investigate the previously evoked link between medium and big Omega powers.\\
Let $M$ be a $\W_m(k)$-module. Recall that, for each nonnegative integer $n$, we have, at our disposal, the $p$-th divided power operation $$\gamma_p : \Gamma^{p^n}(M) \lra   \Gamma^{p^{n+1}}(M).$$ It is a polynomial law, homogeneous of degree $p$. It can thus be viewed as a $\W(k)$-linear map  $$\Gamma^p(\Gamma^{p^n}(M)) \lra   \Gamma^{p^{n+1}}(M),$$  $$[X]_p \mapsto \gamma_p(X), $$ which we   denote by $\tilde \gamma_p$. In the reverse direction, the association  $$M \lra \Gamma^p(\Gamma^{p^n}(M)) $$ $$ x \mapsto [[x]_{p^n}]_p$$ defines a polynomial law, homogeneous of degree $p^{n+1}$. It thus yields a natural  $\W(k)$-linear map  $$\Gamma^{p^{n+1}}(M) \lra \Gamma^p(\Gamma^{p^n}(M)) ,$$ $$[x]_{p^{n+1}} \lra [[x]_{p^n}]_p. $$ We denote it by $\alpha_p$.\\

\begin{defi}
Let $M$ be a $\W_m(k)$-module. Let $n$ be a positive integer. We recursively define $\W(k)$-linear maps $$F_n: \Gamma^{(n)}(M) \lra \Gamma^{p^n}(M) $$ and $$G_n:  \Gamma^{p^n}(M) \lra \Gamma^{(n)}(M)  $$ by setting $$F_1=G_1=\Id,$$ $$ F_{n+1}=\tilde \gamma_p \circ \Gamma^p(F_n)$$  and $$G_{n+1}= \Gamma^p(G_n)\circ \alpha_p. $$ 
\end{defi}

\begin{defi}\label{strange1}
We put $c_1=1$ and, for each  integer $i\geq 2$, we put
$$c_{i}=\frac 1 {p!} {{p^{i}} \choose {p^{i-1},p^{i-1}, \ldots , p^{i-1}}} \in \N.$$ It is an integer, which is  a $p$-adic unit. For $n \geq 1$, we put $$ \mathcal C_n:= c_n c_{n-1}^p \ldots c_2^{p^{n-2}}.$$ It is an integer, which is a $p$-adic unit.
\end{defi}

\begin{lem}
    We have $$F_n \circ G_n=  \mathcal C_n \Id.$$  
\end{lem}

\begin{dem}
The case $n=1$ is obvious.  The general case is by induction on $n$, using the relation $$\gamma_p \circ \gamma_{p^n}=c_{n+1} \gamma_{p^{n+1}}$$ (cf. Proposition \ref{defigammai}).

\end{dem}

Hence, $F_n$ and $G_n$ present $\Gamma^{p^n}(M)$ as a direct factor of $ \Gamma^{(n)}(M)$, which is probably well-known. What is perhaps less standard, is that $F_n$ and $G_n$ are adjoint, for Pontryagin duality.

\begin{lem}\label{mediumbig}
  Let $M$ be a $\W_m(k)$-module. Let $n$ be a positive integer.  For every $X \in \Gamma^{p^n}(M)$ and every $\Phi \in \Gamma^{(n)}(M^\vee),$ we have the formula $$ <X, F_n(\Phi)>=\mathcal C_n <G_n(X), \Phi> \in \W_{m+n}(k).$$

\end{lem}
\begin{dem}
Induction on $n$. The case $n=1$ is obvious. For the induction step, pick $ x \in M$ and $\Phi \in \Gamma^{(n)}(M^\vee).$ We compute: $$ <[x]_{p^{n+1}}, F_{n+1}([\Phi]_p)>= <[x]_{p^{n+1}}, \tilde \gamma_p ( \Gamma^p(F_n)([\Phi]_p))>$$ $$=<[x]_{p^{n+1}},  \gamma_p (F_n(\Phi))> =c_{n+1}\overline \tau_1(<[x]_{p^n},F_n(\Phi)>),$$  where the last equality follows from point 2) of Lemma \ref{formuladual}. On the other hand, we have $$<G_{n+1}([x]_{p^{n+1}}),[ \Phi]_p>  =<\Gamma^p(G_{n})([[x]_{p^{n}}]_p), [\Phi]_p>$$ $$=<[G_n([x]_{p^{n}})]_p,[ \Phi]_p>= \overline \tau_1(<G_n([x]_{p^{n}}), \Phi> ).$$ Comparing the two expressions yields the result (plugging in the formula at the previous step).

\end{dem}





\begin{prop}\label{mediumbig}
 Let $M$ be a $\W_m(k)$-module. The linear map $G_n$ induces, by passing to the quotient, a canonical linear map $$\Psi^n_M: \Omega^n(M) \lra \overline \Omega^n(M) ,$$   compatible with the dualities on both sides. More precisely, we have $$ <\Psi^n_M(X),\Psi^n_{M^\vee}(\Phi)>= \mathcal<X,\Phi>,$$ for all $X \in \Omega^n(M) $ and all $\Phi \in \Omega^n(M^\vee)$. \\In particular, we have a canonical decomposition $$\overline \Omega^n(M) = \Omega^n(M) \bigoplus  \Omega^n(M^\vee)^\perp.$$ 
\end{prop}

\begin{dem}
The existence of $\Psi^n_M$ is a straightforward consequence of the adjunction formula of Lemma \ref{mediumbig}. The second formula is easily checked on pure symbols. The last assertion is a general fact.
\end{dem}

	\section{Functorial properties of Omega powers.}

	
	Let $m$ be a fixed positive integer.
	
	\subsection{Multilinearity.}\label{multilin}
	 Let $M$ and $N$ and $L$ be three   $\W_m(k)$-modules. 
	Let $$ B(.,.): M \times N \lra L$$ be a $\W_m(k)$-bilinear pairing. Let $n$ be a positive integer.
	 
The pairing $B$ induces a pairing $$ B_1:  L^\vee \times M \lra N^\vee,$$ $$ (\phi,x) \mapsto \phi(B(x,.)),$$  which gives rise to 
	 a  $\W_{m+n}(k)$-bilinear pairing $$\Gamma^{p^n}( B_1): \Gamma^ {p^n}(L^\vee)  \times \Gamma^{p^n}(M) \lra \Gamma^{p^n}(N^\vee),$$  $$( [\phi]_{p^n},[x]_{p^n} )\mapsto  [B_1(\phi,x)]_{p^n}.$$ In a similar way, the pairing $$ B_2:  L^\vee \times N \lra M^\vee,$$ $$ (\phi,y) \mapsto \phi(B(.,y)),$$   produces
	 a  $\W_{m+n}(k)$-bilinear pairing $$\Gamma^{p^n}( B_2): \Gamma^{p^n}(L^\vee)  \times \Gamma^{p^n}(N) \lra \Gamma^{p^n}(M^\vee),$$  $$ ([\phi]_{p^n},[y]_{p^n}) \mapsto  [B_1(\phi,y)]_{p^n}.$$ 
	 
	For $x \in M$, $y \in N$ and $\phi \in L^\vee$, it is straighforward to check that $$ <  \Gamma^{{p^n}}(B_1)([\phi]_{p^n},[x]_{p^n}),[y]_{p^n}> = <[B_1(\phi,x)]_{p^n}, [y]_{p^n}> = \overline \tau_n(\phi(B(x,y)))$$$$= <[B(x,y)]_{p^n},[\phi]_{p^n}> =<  \Gamma^{p^n}(B)([x]_{p^n},[y]_{p^n}),[\phi]_{p^n}> $$ $$ =<  \Gamma^{p^n}(B_2)([\phi]_{p^n},[y]_{p^n}),[x]_{p^n}>.$$
	 This shows that we actually have $$ <  \Gamma^{{p^n}}(B_1)(\Phi,X),Y> =< \Gamma^{{p^n}}(B)(X,Y),\Phi> $$ $$ =<  \Gamma^{{p^n}}(B_2)(\Phi,Y),X> \in \W_{m+n}(k),$$ for all $X \in \Gamma^{{p^n}}(M),$ $Y \in \Gamma^{{p^n}}(N)$ and $\Phi \in \Gamma^{{p^n}}(L^\vee).$ Indeed, this can be checked after extending scalars from $k$ to any perfect field extension of $k$ (by Proposition \ref{divbase}, for the base change $\W(k')/\W(k)$). We can hence assume that $k$ is infinite, in which case pure symbols generate divided power modules by Proposition \ref{symbolsgen}. This adjunction formula show that the pairing $\Gamma^{{p^n}}(B)$ is compatible with the duality. Consequently, it  passes to the quotient by its kernel, yielding a pairing of $\W_{ m+n}(k)$-modules $$ \Omega^n(B):  \Omega^n(M)  \times \Omega^n(N)  \lra   \Omega^n(L).$$ 
	 Note that the association 
	 $$B \lra \Omega^n(B) $$ does unfortunately not send "tense" pairings to tense pairings (in the sense of Section \ref{tenseproduct}).

	 \subsection{Omega powers of $\W_m(k)$-algebras}\label{Omegaalg}
	 Let $m,n$ be  positive integers.\\
	 Let $A$ be a (not necessarily finite-dimensional) $\W_m(k)$-algebra. We would like to canonically turn $\Omega^n(A)$ into a  $\W_{m+n}(k)$-algebra, with unit $(1)_n$, and multiplication given by $$(x)_n (y)_n = (xy)_n $$ on pure symbols. This is indeed possible: denoting by $\mu:A \times A \lra A$ the multiplication of $A$ (viewed as a $\W_m(k)$-bilinear pairing), the bilinear map $\Omega^n(\mu)$ of the preceding paragraph does the job.

	 \begin{prop}
	  Let $A$ be a $\W_m(k)$-algebra (in the usual sense), with multiplication $\mu:A \times A \lra A$. Then the $\W_{m+n}(k)$-module  $\Omega^n(A)$ can be canonically turned, via $\Omega^n(\mu)$, into a  $\W_{m+n}(k)$-algebra, with unit $(1)_n$, and  multiplication given by $$(x)_n (y)_n = (xy)_n $$ on pure symbols.
	 \end{prop}
	 
	 \begin{dem}
	 This is clear.
	 \end{dem}
	 
	 \begin{rem}
	 If $A$ is a Hopf algebra over $\W_m(k)$, we can wonder whether $\Omega^n(A)$ is naturally a Hopf algebra over $\W_{m+n}(k)$. With our current definition of $\Omega^n$, this is not the case.
	 \end{rem}
	 To begin with, we treat the instructive case of an \'etale algebra. Remember that the category of \'etale $k$-algebras is equivalent to that of  \'etale $\W_m(k)$-algebras. 
	 
	 \begin{lem}\label{decalg}
	Let $E$ be an \'etale $k$-algebra, of degree $d$. Denote by $l/k$ "the" Galois splitting field of $E$. Then the finite $\W_{n+m}(k)$-algebra $ \Omega^n(\W_m(E)) $ is isomorphic to  a finite  product of local $\W_{n+m}(k)$-algebras of the form $\W_{n_i}(k_i)$, where $1 \leq n_i \leq n+m$ is an integer, and $l/k_i/k$ is an intermediate field extension.
\end{lem}
	
	\begin{dem}	 We first make the following elementary observation. Let $R$ be a finite local  $\W_{i}(k)$-algebra $R$, such that its maximal ideal is $pR$, and  $i$ minimal (i.e. $p^{i-1} \neq 0$ in $ R$).  Then $R$ is canonically isomorphic to $\W_i(l)$, where $l=R/p$ is its residue field. Using Lemma \ref{Omegacommu}, we then see that the statement of the Lemma  is invariant under separable field extensions: we may thus assume that $k= 
	\overline k$ is algebraically closed. But then $E$ is isomorphic to the trivial \'etale algebra $k^d$, and the statement is straightforward (choose the basis of primitive idempotents).
	\end{dem}


	Now, the formula $$ \W_m(E) \lra \W_{n+m}(E),$$ $$ x \mapsto \overline \tau_{n}(x)$$ clearly defines a polynomial law of $\W_{n+m}(k)$-modules, which is homogeneous, of degree $ p^{n}$. Since it is multiplicative, the resulting $\W_{n+m}(k)$-linear map $$\rho=\rho_{E,n}: \Gamma_{\W_{n+m}(k)}^{ p^{n}}(\W_m(E)) \lra \W_{n+m}(E)$$ is actually a ring homomorphism.

\begin{lem}
    
    The homomorphism $\rho$ vanishes on $\Ker(\Delta)$.
\end{lem}

\begin{dem}
Base-changing to an algebraic closure of $k$, we can assume that $k$ itself is algebraically closed. Then, $E \simeq k^s$, and $\W_{n+m}(E) \simeq \W_{n+m}(k)^s$, as $\W(k)$-algebras. The map $\rho$ is then just given by  functoriality from the $s$ canonical projections $\pi_i: E \lra k$, and the claim becomes obvious.
\end{dem}

\begin{defi}\label{rho}
The homomorphism $\rho_{E,n}$ above induces, by passing to the quotient, a homomorphism (of $\W_{n+m}(k)$-algebras) $$ \Omega^n(\W_m(E)) \lra \W_{n+m}(E),$$  $$ (x)_n \mapsto \overline \tau_n(x)$$ which we still denote by $\rho_{E,n}$, or simply by $\rho$. 
\end{defi}
\begin{rem}
In the previous definition, $\overline \tau_n(\tau(x))$ is nothing but $\tau(x^{p^n}),$ where  $\tau: E \lra \W(E)$ is the usual Teichm\"uller representative.
\end{rem}

\begin{rem}
    Note that $\W_{n+m}(E)$ is a free $\W_{n+m}(k)$-module. Hence, by Lemma \ref{decalg}, the homomorphism $\rho$ can be identified with the projection onto a direct factor of the $\W_{n+m}(k)$-algebra  $\Omega^n(\W_m(E))$. In other words, $\Spec(\rho)$ is an open-closed immersion.
\end{rem}

	\subsection{Behaviour of Omega powers under field extensions.}

Let $k'/k$ be an extension of perfect fields of characteristic $p$.\\
Denote by $\tau$ (resp. $\tau'$) the Teichm\"uller representative for $k$ (resp. $k'$) and by $K'$ the field of fractions of $\W(k')$. Pontryagin duality $\Hom_{\W(k)}(.,K/\W(k))$  (resp. $\Hom_{\W(k')}(.,K'/\W(k'))$) will be denoted by $(.)^\vee$ (resp. $(.)^{\vee'}$) . \\ 
	 Omega powers of $\W_m(k)$-modules  (resp of $\W_m(k')$-modules) will be denoted by  $\Omega^n$  (resp.   ${\Omega'}^n$). \\
	 
	 \subsubsection{Extension of scalars.}
Let us first recall two properties of scalars extension, on the level of Witt vectors.

\begin{lem}\label{dualcomm}
	Pontryagin duality commutes with scalars extension, from $\W(k)$ to  $\W(k')$.\\
	More precisely, let $M$ be a $W(k)$-module. Put $M':= M \otimes_{\W(k)} \W(k').$\\ Then the canonical map $$\Hom_{\W(k)} (M,K / \W(k)) \otimes_{\W(k)} \W(k') \lra \Hom_{\W(k')} (M',K' / \W(k')),$$ $$ f \otimes x \mapsto ( m\otimes y \mapsto xyf(m)),$$  is an isomorphism.

\end{lem}
\begin{dem}
	This is clear.
\end{dem}

 Omega powers behave very nicely with respect to the extension $k'/k$.
\begin{lem}\label{Omegacommu}
	The formation of (medium  and big) Omega powers commutes to extending scalars from $k$ to $k'$. In other words, let $M$ be a (finite) $\W_m(k)$-module.  Put  $$M':= M \otimes_{\W_m(k)} \W_m(k').$$  We then have  canonical isomorphisms of $\W_{m+n}(k)$-modules $$ \Omega^n(M) \otimes_{\W_{m+n}(k)} \W_{m+n}(k') \simeq  {\Omega'}^n(M')$$ and $$ \overline \Omega^n(M) \otimes_{\W_{m+n}(k)} \W_{m+n}(k') \simeq  {\overline \Omega'}^n(M').$$ 
\end{lem}	 
\begin{dem}
	We know that the formation of divided powers commutes to extension of the base ring. Using Lemma  \ref{dualcomm}, we see that the duality arrow $$\Delta: \Gamma^{ p^n}(M) \lra  \Gamma^{p^n}(M^\vee) ^\vee$$ thus also commutes to extending scalars from $k$ to $k'$, in the sense that $\Delta \otimes_{\W(k)} \W(k')$ is canonically isomorphic to $\Delta'$. The claim follows.
\end{dem}	 
\subsubsection{Restriction of scalars.}
Assume now that $k'/k$ is finite, of degree $s$.\\
Let  $M$ be a $\W_m(k')$-module. We can also view it as a $\W_m(k)$-module. Applying the process of Section \ref{multilin} to the $\W_m(k)$-bilinear pairing $$ \W_m(k') \times  M'  \lra M'$$ $$ (\lambda , v') \mapsto \lambda v',$$ we can endow $\Omega^n(M')$  with a canonical structure of a $\Omega^n(\W_m(k'))$-module. On pure symbols, we have the formula $$(\lambda)_n(v')_n=(\lambda v')'_n .$$ \\Now, recall the homomorphism $$ \rho: \Omega^n(\W_m(k')) \lra \W_{n+m}(k')$$ of Definition \ref{rho}.

The natural quotient map $$\pi: \Omega^n(M') \lra \Omega'^n(M') ,$$ $$ (v')_n \mapsto (v')'_n$$ is compatible with $\rho$: we have $$ \pi( a.x)=\rho(a) \pi(x),$$ for all $a \in \Omega^n(\W_m(k'))$ and all $x \in \Omega^n(M')$. \\We thus have a canonical $\W_{m+n}(k')$-linear map $$ \Psi: \Omega^n(M') \otimes_\rho \W_{n+m}(k')  \lra \Omega'^n(M').$$

\begin{prop}
 
 The map $\Psi$ above is an isomorphism.
\end{prop}
\begin{dem}

Extending scalars to an algebraic closure of $k$, we can replace $k'/k$ by the trivial \'etale algebra $k^s/k$. The data of $M$ is now the data of $s$ $\W_m(k)$-modules $M_1,\ldots, M_s$, and both sides equal the direct sum of the $s$ modules $\Omega^n(M_i)$ .
\end{dem}

\section{Frobenius and Verschiebung, for Omega powers.}\label{FrobVerOmega}

	In the case $m=1$, Omega powers naturally inherit Frobenius and Verschiebung maps, from divided powers. We  now  explain how.\\ 
	Let $V$ be a $k$-vector space.  Let $n \geq 1$ be an integer.
	Note that the (surjective) linear maps   $$\Frob_V: \Gamma^{p^{n+1}}(V) \lra\Gamma^{p^{n}}(V^{(1)})$$ and  $$\Frob_{V^*}: \Gamma^{p^{n+1}}(V^*) \lra\Gamma^{p^{n}}(V^{*(1)})$$ satisfy the formula $$<\Frob_{V^*}([\phi]_{p^{n+1}}),\Frob_{V}([v]_{p^{n+1}}) > = <[\phi^{(1)}]_{p^{n}},[v^{(1)}]_{p^{n}}) > $$ $$= \overline \tau_n (\phi^{(1)} (v^{(1)}) )= \overline \tau_n(\phi(v))^p=<[\phi]_{p^{n+1}},[v]_{p^{n+1}} >,$$ modulo $p^{n+m}$. Hence, they respect the duality $\Delta$, and the following definition makes sense.
	\begin{defi}
	The Frobenius map $$\Frob: \Gamma^{p^{n+1}}(V) \lra\Gamma^{p^{n}}(V^{(1)})$$ yields, by passing to the quotient, a $\W(k)$-linear map $$\Omega^{n+1}(V) \lra \Omega^n(V^{(1)}).$$  It is the Frobenius, for (medium) Omega powers. By perfect duality, the dual of the Frobenius for $V^*$ is a  $\W(k)$-linear map $$ \Omega^n(V^{(1)}) \lra \Omega^{n+1}(V).$$ It is the Verschiebung, for (medium) Omega powers.
	\end{defi}

	As one can expect, we can define Frobenius and Verschiebung for big Omega powers, in a way compatible with the natural embedding. Here is how.\\

	Applying $\Gamma^{(n-1)}$ to the Frobenius map $$\Gamma^{(1)}(V)  \lra V^{(1)}, $$ $$[v]_p \mapsto v^{(1)}, $$ yields a surjective $\W(k)$-linear map $$ \Frob_V^{(n)}: \Gamma^{(n)}(V) \lra \Gamma^{(n-1)}(V^{(1)})$$ as $\Gamma^p$ commutes to Frobenius twist.\\
	Dually, we get a  $\W(k)$-linear map $$ \Frob_{V^*}^{(n)}: \Gamma^{(n)}(V^*) \lra \Gamma^{(n-1)}(V^{*(1)}).$$ For $X \in  \Gamma^{(n)}(V)$ and $\Phi \in  \Gamma^{(n)}(V^*)$, we check by induction on $n$ that $$p< X,\Phi> = < \Frob_V^{(n)}(X), \Frob_{V^*}^{(n)}(\Phi)> \in \W_{n+1}(k).$$ Here the righthand side, a priori belonging to $\W_{n}(k)$, is viewed as an element of $\W_{n+1}(k)$ via the inclusion (Verschiebung) 
 $$\W_{n}(k) \stackrel{1 \mapsto p}\lra  \W_{n+1}(k).$$
 We can also choose to write this equality as $$< X,\Phi> = < \Frob_V^{(n)}(X), \Frob_{V^*}^{(n)}(\Phi)> \in \W_{n}(k),$$ 
 where the lefthand side is taken modulo $p^{n}$... \\

Hence, $ \Frob^{(n)}$  yields by passing to the quotient a $\W(k)$-linear map $$  \overline \Omega^{n}(V) \lra  \overline \Omega^{n-1}(V^{(1)}).$$
	
	\begin{defi}
	
	The $\W_{n+1}(k)$-linear map $$  \overline \Omega^{n}(V) \lra \overline \Omega^{n-1}(V^{(1)})$$ that we have just defined is the Frobenius homomorphism,  for big Omega powers. It will simply be denoted by $\Frob_V$, or even by $\Frob$, if the context is clear.  \\
	Using the  duality between $\overline \Omega^{n}(V) $ and $\overline \Omega^{n}(V^*) ,$ we define the Verschiebung homomorphism  $$\Ver_V : \overline \Omega^{n-1}(V^{(1)}) \lra \overline \Omega^{n}(V)$$ to be dual to $\Frob_{V^*}$.
	\end{defi}
	\begin{rem}
	To be more precise, we can define the Verschiebung for big Omega powers for a finite-dimensional $V$ first, and then define it for $V$ arbitrary using a direct limit argument.
	\end{rem}
\begin{lem}
  The Frobenius and the Verschiebung  for (medium or big) Omega powers are adjoint operators, satisfying  
  $$\Ver \circ \Frob=p $$    and   $$\Frob \circ \Ver=p. $$ They are compatible with the natural inclusion $\Omega ^n \subset \overline \Omega^n.$
\end{lem}
\begin{dem}
We check the first part, for big Omega powers.
    That these operators are adjoint is clear from the definition of the Verschiebung. For $X \in  \Gamma^{(n)}(V)$ and $\Phi \in  \Gamma^{(n)}(V^*)$, the computation  $$<\Ver(\Frob(X)),\Phi>=<\Frob(X), \Frob(\Phi)>=p <X, \Phi> \in \W_{n+1}(k)$$ ensures that $\Ver \circ \Frob=p $. Since $\Frob$ is surjective, the other equality follows.
    \end{dem}

	\begin{prop}
		Let $V$ be a  $k$-vector space, and let $n$ be a positive integer. We have exact sequences $$\overline \KW_1(V)=\overline \KW_1(V,n): 0 \lra  \overline \Omega^{n-1}(V^{(1)}) \stackrel \Ver \lra  \overline \Omega^n(V) \lra \overline  \Omega^n(V)/p \lra 0$$ and $$\overline \KW_2(V)=\overline \KW_2(V,n): 0 \lra \overline  \Omega^{n-1}(V)[p]  \lra \overline \Omega^{n}(V)\stackrel \Frob \lra   \overline \Omega^{n-1}(V^{(1)}) \lra 0.$$ We shall refer to them as the first and second Kummer-Witt exact sequences for (big) Omega powers, respectively.  They are dual constructions, in the sense that, if $V$ is finite-dimensional, $\overline \KW_1(V^*)$ is canonically isomorphic to  $\overline \KW_2(V)^\vee$.\\
		
		We define $ \KW_2(V,n)$ and  $\KW_2(V,n)$, for medium Omega powers, in the same way.
	
	\end{prop}
	
	\begin{dem}
	Only the exactness of the sequences in question has perhaps to be checked. Note that $\overline \KW_1(V)$ is clearly exact on the right and in the middle. The injectivity of $\Ver$ follows by duality from the surjectivity of $\Frob$.
	\end{dem}

\section{The Transfer.}\label{sectrans}
In this section, $k$ is a finite field, of cardinality $q=p^r$. The goal of this section is to define the Transfer, which is a polynomial law "in the wrong direction". More precisely, if $W \lra V$ is an inclusion of $k$-vector spaces, of finite codimension $c$, we shall build a canonical polynomial law $$T_{W,V}:\A_k(V) \lra \A_k(W), $$ the Transfer,  enjoying nice properties.\\
As always, the letter $n$ denotes any positive integer. 
For every $k$-vector space $V$, we  have a canonical $k$-linear isomorphism $$V \simeq V^{(r)},$$ which we shall tacitly use to identify these two $k$-vector spaces.\\ The $r$-th Frobenius \textit{polynomial law}  $$V \lra V^{(r)}=V,$$ $$ v \mapsto v^{(r)}$$ will be denoted by $F_V$, or simply by $F$ if no confusion arises. It is homogeneous, of degree $p^r$.

\begin{defi} \label{det}
	Let $ V$ be a finite-dimensional $k$-vector space, of dimension $d$.  We put $$ \Det^n(V):=\Omega^{n}(\Det(V));$$ it is a free $\W_{n+1}(k)$-module of rank one.
 	\end{defi}
\subsection{Laws in one variable.}
 
  \begin{defi}
   	Let $ V$ be a finite-dimensional $k$-vector space, of dimension $d$. Let $$S=\{s_0 < s_1 < \ldots <  s_{m-1} \} \subset \{0, \ldots, d\}$$ be any subset, of cardinality $m$. The formula $$ \A_{k}(V) \lra \A_{k}( \Lambda^m(V))$$ $$v \mapsto F^{s_0}(v) \wedge F^{s_1}(v) \wedge \ldots   \wedge  F^{s_{m-1}}(v)$$  defines a  polynomial law, which is homogenous,  of degree $$q^{s_0}+q^{s_1}+ \ldots+ q^{s_{m-1}}.$$
   It is \textit{the exterior power in one variable}, with respect to $V$ and $S$. We denote it by $ \overline \lambda^S_V$, or simply by $\overline \lambda^S$, if the dependence in $V$ is clear. \\If $S=\{0, \ldots, m-1\}$, we denote  $ \overline \lambda^S$ by $ \overline \lambda^m$.\\
   The law $\overline \lambda_V^d$ will be denoted by $\det_V^1$. It is the determinant in one variable.
\end{defi}
 
 \begin{rem}\label{lambdavanish}
 The locus where the law $ \overline \lambda^m$ vanishes is exactly the (finite) union of all linear subvarieties of $\A_k(V)$, which are of dimension stricly less than $m$.
 \end{rem}
 
 The determinant in one variable $\det_V^1$ is, in fact, the product of all nonzero $k$-linear forms  on $V$ (up to scalar multiplication). Let us make this statement more precise. We thank Ofer Gabber for  an interesting discussion, which helped us clarify the exposition.

\begin{defi}
	Let $ V$ be a finite-dimensional $k$-vector space, of dimension $d \geq 2$. For each $k$-rational hyperplane $H \subset V$, denote by $$\pi_H: V \lra V/H$$ the $k$-linear projection. Put $$ \overline \Det (V):= \bigotimes_{H \subset V} (V/H),$$ where the tensor product is taken over all hyperplanes $H \subset V$. It is a one-dimensional $k$-vector space.\\ Denote by $\mathrm{\overline {\det}}^1$ the polynomial law $$ \A_k(V) \lra \A_k (\overline  \Det V),$$  $$ v \mapsto    \otimes_{H \subset V} (\pi_H(v)).$$  It is homogeneous, of degree  $1+q+ \ldots + q^{d-1}= \vert  \P_k(V) \vert$.
	
	\end{defi}
	
\begin{prop}\label{productlin}
	Let $ V$ be a finite-dimensional $k$-vector space, of dimension $d \geq 2$. There exists a canonical  isomorphism $$\theta:\overline \Det (V) \lra \Det (V)$$ of one-dimensional $k$-vector spaces,  such that $$ \theta \circ \overline \det{}^1= \det{}^1.$$ 
    
\end{prop}

\begin{dem}
Choose coordinates $V \simeq k^d$. Then $ \det^1$ is given by a polynomial $$P \in k[X_1, \ldots, X_d],$$ which is homogeneous, of degree $1+q+ \ldots + q^{d-1}= \vert  \P_k(V) \vert$. For each hyperplane $H \subset V$, let $$L_H \in k[X_1, \ldots, X_d]$$ be a linear polynomial, with kernel $H$.\\
    Let $H \subset V$ be a $k$-rational hyperplane. It is clear that the composite $$\A_k(H) \stackrel {can} \lra \A_k(V) \stackrel {\det^1} \lra \A_k(\Det(V))$$ identically vanishes (indeed, $\Lambda_k^d(H)=0$). Thus, the polynomial $P$ has to be divisible by $L_H$. Since these linear polynomials, for various $H$, are two by two coprime, $P$ has to be divisible by their product $$Q:= \Pi_{H \subset V} L_H,$$ which is a homogeneous polynomial of the same degree as $P$. Hence, $P=Q$ up to a nonzero scalar. The statement of the proposition is, obviously, the canonical translation of this fact.
\end{dem}
\begin{exo}
    Let $k'/k$ be a finite field extension, of degree $n$. Let $V$ be a $d$-dimensional $k$-vector space. Show that $\det^1_V$ identically vanishes on $k'$-rational points if and only if $n <d$.
\end{exo}
 \begin{exo}
   	Let $ V$ be a finite-dimensional $k$-vector space, of dimension $d$.  For $i=1 \ldots d$, denote by  $$F_i: \A_{k}(V) \lra \A_{k}( \Det(V))$$  the polynomial law $\lambda^d_{\{0,1, \ldots, \hat i, \ldots d\}}$, where the symbol $\hat i$ means that $i$ is omitted.	Then  the morphism of affine $k$-varieties $$F: \A_{k}(V) \lra \A_{k}( \Det(V))^d $$ $$ v \mapsto (F_1(v), \ldots, F_d(v))$$ is \'etale exactly outside the (finite) union of all $k$-rational hyperplanes of $V$.
\end{exo}
\subsection{The Transfer, as a polynomial law.}
We now proceed one step further towards the definition of the Transfer. In the exterior algebra, it is easy to define such an operation, in a $k$-linear way. We apologize for the choice of the terminology "exterior transfer" in what follows- it is a bit pompous... 
 \begin{defi}
 
 	Let $ V$ be a finite-dimensional $k$-vector space, of dimension $d$. Let $W \subset V$ be a $k$-linear subspace, of codimension $c$.  Consider the exact sequence $$ 0 \lra W \lra V \lra V/W \lra 0$$ and its $k$-dual sequence $$ 0 \lra (V/W)^*=W^\perp \lra V^* \lra W^* \lra 0.$$
 	Pick an integer $m\geq c$.
 Then the wedge product $$ \Lambda^c(W^\perp)  \otimes_k  \Lambda^{m-c}(V^* )\lra \Lambda^{m}(V^*),$$ $$(x,y) \mapsto x \wedge y $$passes to the quotient by the arrow $V^* \lra W^*$, yielding an injective $k$-linear map $$ \Det(W^\perp) \otimes_k \Lambda^{m-c}(W^*) \lra \Lambda^{m}(V^*).$$ Its $k$-dual is a surjective map $$\Lambda^{m}(V) \lra \Det(V/W) \otimes \Lambda^{m-c}(W) ,$$ which we denote by $\lambda T^m_{W,V}$. \\
 It is the exterior transfer, from $V$ to $W$.
 \end{defi}
 
 \begin{rem}
The linear map  $\lambda T^m_{W,V}$ is explicitly given by the formula $$v_1 \wedge \ldots \wedge v_{m} \mapsto \sum_{I=\{i_1 < \ldots < i_c\}} \epsilon(I) (\pi (v_{i_1}) \wedge  \ldots  \wedge \pi (v_{i_c}) )\otimes (\rho(v_{j_1}) \wedge \ldots \wedge \rho(v_{j_{m-c}})),$$ where $\pi : V \lra V/W$ is the quotient map,  $\rho: V \lra W$ is -any- $k$-linear retraction of the inclusion $W \lra V$, and the sum ranges over all subsets $$I=\{ i_1 < \ldots < i_c\} \subset \{ 1, \ldots, m\}, $$ with complement $ I^c=\{ j_1 < \ldots < j_{m-c}\}$. \\The number $\epsilon(I) \in \{1,-1\}$ is a sign, which is not hard to compute.
 \end{rem}

\begin{lem}\label{lambdacomp}
    Let $$Z \subset W \subset V$$ be three finite-dimensional $k$-vector spaces. Denote by $c$ (resp $c'$) the codimension of $W$ in $V$ (resp. of $Z$ in $W$). Let $m \geq c+c'$ be an  integer. Through the canonical isomorphism $$ \Det(V/Z) \simeq  \Det(V/W) \otimes_k  \Det(W/Z), $$ the exterior transfers satisfy the 'cocycle' condition $$ \lambda T^{m-c}_{Z,W} \circ \lambda T^m_{W,V} = \lambda T^{m}_{Z,V},$$  as linear maps $$ \Lambda ^{m}(V) \lra  \Det(V/Z) \otimes \Lambda^{m-c-c'}(Z).$$
\end{lem}

\begin{dem}
Looking at the definition of the exterior transfer, the $k$-dual statement of this Lemma boils down to the associativity of the wedge product.

\end{dem}

We can now define the Transfer, as a polynomial law.

\begin{prop} 
Let $V$ be a (finite-dimensional) $k$-vector space. Let $W \subset V$ be a $k$-linear subspace, of codimension $c$. Denote by $$ \pi: V \lra V/W$$ the projection .\\ Then, there exists a unique polynomial law $$ T_{W,V}: \A_k(V) \lra \A_k(W),$$ such that $$ \lambda T_{W,V}^{c+1} \circ \overline \lambda_V^{c+1} =  (\det{}_{V/W}^1 \circ \pi)\otimes T_{W,V},$$  as polynomial laws $$ \A_k(V) \lra \A_k(\Lambda^{c+1}_k(V)) \lra \A_k(\Det(V/W)  \otimes_k W).$$ It is homogeneous, of degree $q^c$.

\end{prop}

\begin{dem}
Uniqueness is clear: we simply have to see that $(\det{}_{V/W}^1 \circ \pi)$ divides $ \lambda T_{W,V}^{c+1} \circ \overline \lambda_V^{c+1}$ (as polynomial laws).
At the light of Proposition \ref{productlin}, it suffices to show that the law $ \lambda T_{W,V}^{c+1} \circ \overline \lambda_V^{c+1}$ identically vanishes on all $k$-rational hyperplanes of $V$, which contain $W$.  Let $H \subset V$ be such a hyperplane. Then the composite $$\A_k(H) \lra \A_k(V) \stackrel {\overline \lambda_V^{c+1}} \lra \A_k(\Lambda^{c+1}(V))$$ takes values in $\A_k(\Lambda^{c+1}(H))$. But  $ \lambda_{W,V}^{c+1}$ vanishes on $\A_k(\Lambda^{c+1}(H))$. To see this, first use Lemma \ref{lambdacomp} to reduce to the case $W=H$. It then becomes obvious, by definition of  $\lambda T_{W,V}^{c+1}$.
\end{dem}

\begin{defi} \label{defitrans}
Let $V$ be a (finite-dimensional) $k$-vector space. Let $W \subset V$ be a $k$-linear subspace, of codimension $c$. The polynomial law $$ T_{W,V}: \A_k(V) \lra \A_k(W)$$ constructed in the previous Proposition  is the Transfer, from $V$ to $W$. It is homogeneous, of degree $q^c$.

\end{defi}

\begin{rem}
    Naively speaking, the $T_{W,V}$ can be interpreted as 'the extension by zero' of the inclusion $W \lra V$, to the whole $V$. Proposition \ref{Frobtrans} makes this statement precise.
\end{rem}
\begin{lem}\label{transcomp}
    Let $$Z \subset W \subset V$$ be three finite-dimensional $k$-vector spaces. Then the Transfers satisfy the 'cocycle' condition $$ T_{Z,W} \circ T_{W,V} = T_{Z,V},$$  as polynomial laws $$ \A_k(V) \lra  \A_k(Z).$$
\end{lem}

\begin{dem}
This is clear, from the definition of the Transfer, together with Lemma \ref{lambdacomp} and Proposition \ref{productlin}.
\end{dem}

In codimension one, the Transfer is given by a very simple formula.
\begin{lem}\label{transcod1}
Let $H \subset V$ be a hyperplane inclusion, with $V$ a finite-dimensional $k$-vector space. Let $\pi: V \lra k$ be a nonzero linear form with kernel $H$. Then the Transfer $$ T_{H,V}: \A_k(V) \lra \A_k(H)$$ is given by the formula $$ v\mapsto F(v)-\pi(v)^{q-1}v. $$
\end{lem}

\begin{dem}
We identify $V/H$ and $k$, through $\pi$.\\
The composite of the exterior power in one variable $$ \lambda^2: \A_k(V) \lra \A_k(\Lambda^2_k(V)) $$ $$v \mapsto v \wedge F(v) $$ with the  exterior transfer $$  \A_k(\Lambda^2_k(V))  \lra \A_k( H),$$ $$ v \wedge w \mapsto \pi(v) w - \pi(w)v$$ is easily computed to be $$ v \mapsto \pi(v) F(v)-\pi(v)^q v.$$ Indeed, $\pi$ is defined over $k$, hence commutes with $F=\Frob^r$. By  definition of the Transfer, dividing by $\pi(v)$ yields the result.


\end{dem}

\begin{defi}

	    	Let $V$ be a (finite-dimensional) $k$-vector space, of dimension $d$. Let $c$ be a positive integer, with $ c \leq d-1$. The law $$ T^c=T^c_V: \A_k(V) \lra \A_k( \bigoplus_{W \subset V} W)$$ $$ v \mapsto (T_{W,V}(v))_{W\subset V},$$ where the direct sum is taken over all $c$-codimensional $k$-linear subspaces $W \subset V$, will be call the Total Transfer for $V$, in codimension $c$.
	\end{defi}
	
		\begin{prop}\label{Frobtrans}
	
		Let $V$ be a (finite-dimensional) $k$-vector space. Let $W \subset V$  be a $k$-linear subspace, of codimension $c \geq 1$. \\
 Then the composite $$\A_k(W)  \stackrel {can} \lra \A_k(V) \stackrel {T_{W,V}} \lra \A_k( W) $$ equals $F^c=\Frob^{rc}$. 

\end{prop}

\begin{dem}
From Lemma \ref{transcomp}, we can assume $c=1$ by induction. The formula then clearly follows from the formula given in Lemma \ref{transcod1}.
\end{dem}

	The next Proposition is much more remarkable. 
	\begin{prop}\label{frobintegral}(Frobenius Integral Formula.)
	
		Let $V$ be a finite-dimensional $k$-vector space, of dimension $d \geq 2$. Let $c$ be an integer, with $1 \leq c \leq d-1.$  \\
 Then the composite $$\Phi:= \A_k(V) \stackrel {T^c_V} \lra \A_k(\bigoplus_W W) \lra \A_k(V)$$ equals $F^c=\Frob^{rc}$, where the second map is given by the (finite!) sum of the inclusions $W \lra V$.\\

\end{prop}
\begin{dem}
Induction on $c$. Let us first deal with the case $c=1$. For a $k$-rational hyperplane $H \subset V,$ the Transfer $T_{H,V}$ is a polynomial law of degree $q$. By the universal property of divided powers, together with Lemma \ref{transcod1}, it is given by the $k$-linear map $$ \Gamma_k^q(V) \lra H,$$ $$[v]_q \mapsto   F(v)-\pi_H(v)^{q-1}v.$$ Pick a  ($k$-rational) $v \in V$. Then $F(v)=v$, and $T_{H,V}(v)$ equals $0$ if $v \notin H$ (in which case  $\pi_H(v)^{q-1}=1$), or equals $v$  if $v \in H$ (in which case  $\pi_H(v)^{q-1}=0$). Since  the number of hyperplanes $H$ containing $v$ is congruent to $1$ modulo $p$, summing over all $H$ shows that $\Phi(v)=v=F(v)$.  Now, $F=\Frob^r$ and $\Phi$ are both polynomial laws of the same degree $q$. Since we know, by Lemma \ref{symbolsgen}, that $k$-rational symbols $[v]_q$ generate $ \Gamma_k^q(V)$, we can indeed conclude that $\Phi=F.$\\

 For the induction step, look at the composite $$ \A_k(V) \stackrel  {T^1_V} \lra \A_k(\bigoplus_{H \subset V} H)  \stackrel  {\sum T_{W,H}} \lra \A_k(\bigoplus_{W \subset H \subset V} W) \lra  \A_k(\bigoplus_{H \subset V} H) \lra \A_k(V), $$ where the direct sums are taken over all hyperplanes $H \subset V$, and all inclusions $W \subset H \subset V$ of a $c$-codimensional $W$  into an hyperplane $H$, respectively (and where the last two arrows on the right are the canonical linear surjections). \\ On the one hand, using Lemma \ref{transcomp}, together with the fact that the cardinality of a projective space over a finite field is congruent to $1$ modulo $p$, we see that this composite equals $\Phi$. On the other hand, the composite of the two middle arrows equal $F^{c-1}$ by induction, so that, using the case $c=1$, the composite of all four arrows equals $F \circ F^{c-1}=F^c$ (note that $T_V^1$ obviously commutes with $F$). The proof is complete.

\end{dem}

\begin{defi}
Let $V$ be a (finite-dimensional) $k$-vector space. Let $W \subset V$  be a $k$-linear subspace, of codimension $c \geq 1$. By the universal property of divided powers, there exists a unique $\W(k)$-linear map $$ \Gamma^{p^{n+rc}}(V) \lra \Gamma^{p^n}(W), $$ $$[v]_{p^{n+rc}} \mapsto [T_{W,V}(v)]_{p^n}. $$
We shall denote it by $\Gamma T^n_{W,V}$. It is the Transfer, for divided powers.

\end{defi}

Let $H \subset V$ be a $k$-rational hyperplane.\\
In view of the preceding definition, it is natural to wonder whether we can  define, by passing to the quotient,  a descending transfer $$ \Omega^{i+r}(V) \lra \Omega^{i}(H),$$ for $i \geq 0$. It is doable for $i \leq 1$. After a few unsuccessful attempts to do so for $i \geq 2$, we noticed that the difficulty can be bypassed, by considering only the submodule of $\Omega^n(V)$ generated by pure symbols. Since $k$ is finite, this submodule, for $n$ large, is much smaller that $\Omega^n(V)$. In particular,  its rank (as a $\W(k)$-module) is bounded by the cardinality of the finite projective space $\P_k(V)$, whereas that of the whole space  $\Omega^n(V)$ grows (a priori doubly exponentially!) to infinity with $n$.\\
At the present moment, we do believe that the submodule of  $\Omega^n(V)$ generated by pure $k$-rational symbols is the right object, for applications to Galois cohomology- and perhaps to other areas. It is the small Omega power functor. We  elaborate on this new object in the next section.

\section{Small Omega powers.}\label{small}
In this section, $k$ is a finite field, of cardinality $q=p^r$. We denote by $m$ a positive integer.\\
\begin{defi}(Small Omega powers.)
Let $M$ be a $\W_m(k)$-module. Let $n$ be a positive integer.\\
We define $$\underline \Omega_m^n(V) \subset \Omega_m^n(V)$$ to be the $\W_{n+m}(k)$-submodule spanned by all pure symbols $(x)_n$, with $x \in M$.\\
It is the $n$-th Small Omega power of $M$. We will simply denote it by $\underline \Omega^n(M)$, if the dependence in $m$ is understood.

\end{defi}

\begin{rem}

The small Omega power $\underline \Omega^n$ is a functor, from the category of $\W_m(k)$-modules to that of $\W_{n+m}(k)$-modules. Contrary to medium and big Omega powers, it is clear from the definition that small Omega powers do not commute to extending scalars to a larger finite field.

\end{rem}

We will now show that small Omega powers naturally occur as a direct summand of medium Omega powers. This can be compared to the occurence of medium Omega powers as a direct summand of big Omega powers. We will need  a computation in finite fields, which presents similarities with Gauss sums.
\subsection{Small Omega powers as a direct factor of medium Omega powers.} 	Let $k'/k$ be 'the' finite field extension, of degree $s$. The field $k'$ has $q^s$ elements.\\ 	\begin{defi}
		We put $$\kappa':=\Hom_{k'}(k',k).$$ It is the $k$-linear dual of $k'$, viewed as a $k$-vector space.
	\end{defi}
	Denote by $\tau: k' \lra \W(k')$ the Teichm\"uller representative of $k'$. Its restriction to $k$ is the Teichm\"uller representative $\tau$ of $k$.
	
	\subsubsection{A funny  computation in finite fields.}
	Denote by $\tau: k' \lra \W(k')$ the Teichm\"uller representative of $k'$. Its restriction to $k$ is the Teichm\"uller representative  of $k$.
	 Denote by $$\tr: k' \lra k$$ the trace map.  We know that the Galois group of the extension $k'/k$ is cyclic of order $s$, generated by the Frobenius $x \mapsto x^q$. Hence, for $z \in k'^*,$ we have $$\tr(z)= \sum_{i=0} ^{s-1} z^{q^i}.$$ 

	\begin{lem}\label{strange2}
	Put $$C=C(k,k'):=  \sum_{z \in k'^*} \tau(z)^{-1} \tau(\tr(z))   \in \W(k').$$ Then $C$ belongs to $\Z_p$, and $C$ is congruent to $-1$ modulo $p$.

	\end{lem}
	
	\begin{dem}
The fact that $C$ belongs to $\Z_p$ is clear: $C \in \W(k')$ is invariant by the Frobenius of $\W(k')/ \Z_p$ .\\ Modulo $p$, we have  $$ C \equiv  \sum_{z \in k'^*} z^{-1} \left( \sum_{i=0} ^{s-1} z^{q^i} \right)=-1 \in \W(k')/p=k'.$$ Indeed, for any integer $N$, the quantity $$ \sum_{y \in  k'^*} y^N$$ vanishes, except when $(q^s-1)$ divides $N$, in which case its value is $-1$.

	\end{dem}

	\begin{defi}
	The number $C=C(k,k') \in \Z_p^\times$  of the previous Lemma will be called the conductor of the extension $k'/k$. \\
	For each linear form $f \in \kappa'$,  denote by $y \in k'$ the unique element such that $$f(.)=\tr(y.).$$ We put $$C(f):= \frac 1 C \tau(y).$$
	\end{defi}
	\begin{prop}\label{fun1}
Let $C$ be the  conductor of $k'/k$.\\

	For every $x \in k'^*,$ we have the formula $$ \tau(x)= \frac 1 C \sum_{y \in k'^*} \tau(y)^{-1} \tau(\tr(xy))   \in \W(k').$$ 
	\end{prop}	
\begin{dem}
    This is clear from the previous Lemma, setting $z=xy$.
\end{dem}

	\begin{coro}\label{fun2}

	For every $x \in k'^*,$ we have the formula $$\tau(x)=  \sum_{ f \in \kappa'} C(f) \tau(f(x))  \in \W(k').$$
	\end{coro}	

\begin{rem}
    The authors are grateful to Pierre Colmez for helping us to clarify the exposition of the previous formula.
\end{rem}
	
\subsubsection{ Perfect duality for small Omega powers.}

We first need a reformulation of Lemma \ref{fun2}. We use freely the notation of the preceding subsection. \\

\begin{lem}
    Let $V$ be a (finite-dimensional) $k$-vector space. Let $n$ be a positive integer.  Let $k'/k$ be a finite field extension, of degree $s$.  Put $V':=V \otimes_k k'$.  Pick an element $\phi' \in V'^{* '} (= V^* \otimes_{k} k').$ For each $k$-linear form $f \in \kappa'$, denote by $f(\phi') \in V^*$ the  composite $$ V \stackrel {x \mapsto x \otimes 1} \lra V' \stackrel {\phi'} \lra k' \stackrel {f} \lra k.$$ We then have the relation $$<X,(\phi')_n>=  \sum_{ f \in \kappa'} C(f) < X, (f \circ \phi')_n>  \in \W_{n+1}(k'),$$ for all $X \in \underline \Omega^n(V)$.

\end{lem}
\begin{dem}

It is enough to check this when $X$ is a pure symbol $(v)_n$, for $v \in V$. The formula follows from Lemma \ref{fun2}, applied (modulo $p^{n+1}$) to $x:=\phi'(v) ^{p^n}\in k'$.
\end{dem}

 The next Proposition  is a key.
\begin{prop}\label{smallmedium}
Let $M$ be a (finite) $\W_m(k)$-module. Let $n$ be a positive integer.

Consider the natural embeddings $$ \underline \Omega^n(M) \lra \Omega^n(M) $$ and $$ \underline \Omega^n(M^\vee) \lra \Omega^n(M^\vee) .$$  The  perfect duality  $$ \Omega^n(M) \times \Omega^n(M^\vee) \lra \W_{m+n}(k)$$ yields  by restriction a duality $$ \underline \Omega^n(M) \times \underline \Omega^n(M^\vee) \lra \W_{m+n}(k).$$
This duality is perfect.

\end{prop}

\begin{dem}
We reduce to the case $m=1$,  by induction.
We  have to show the following. Let $X \in  \underline \Omega^n(M)$  be orthogonal to $\underline \Omega^n(M^\vee)$. Then $X$ is orthogonal to the whole of $ \Omega^n(M^\vee)$ (and hence vanishes). To do so, let $k'/k$ be a finite field extension, such that $k'$ has cardinality greater than $p^n$. Put $M':=M \otimes_{\W(k)} \W(k')$. Denote by $\Omega'$ the (small or medium) Omega powers of $\W(k')$-modules. By Lemma \ref{symbolsgen0}, the inclusion $\underline \Omega'^n(M') \subset \Omega'^n(M')$ is an equality. Since we now that the formation of medium Omega powers commutes to base change, it is enough to show that, for every $\phi' \in (M')^\vee$, we have $$ <X, (\phi')_n>=0 \in \W_{n+1}(k').$$  This follows from the preceding Lemma, since $X$ is orthogonal to $\underline \Omega^n(M^\vee)$.

\end{dem}

\begin{rem}
Note that the perfect duality $$ \underline \Omega^n(M) \times \underline \Omega^n(M^\vee) \lra \W_{m+n}(k)$$ of the previous Proposition is given on pure symbols by $$< (x)_n,(\phi)_n>=\overline \tau_n(\phi(x)). $$ 

\end{rem}

\begin{coro}
   
  With the notation of the preceding Proposition, we have a natural direct sum decomposition $$ \underline \Omega^n(M) = \underline \Omega^n(M) \bigoplus \underline \Omega^n(M^\vee)^\perp .$$
\end{coro}

\begin{dem}
Clear.
\end{dem}

The next Lemma is helpful for defining the Transfer for small Omega powers, in the next section.

	\begin{lem}\label{omegainj} Let   $M$ be  a  $\W_m(k)$-module. Let $n$ be a positive integer. Then the canonical map $$\underline \Omega^n(M) \lra \bigoplus_{L  }\underline \Omega^n(L)$$ is injective, where the direct sum is taken over all split surjective linear maps $$M \lra L $$  of $M$ onto a (not necessarily free) $\W_m(k)$-module of rank one.

	\end{lem}
	
	\begin{dem}
	   We can assume that $M$ is a $\W_m(k)$-module of finite-type. By perfect duality (Proposition \ref{smallmedium}), it is then equivalent to show that the canonical map $$\bigoplus_{L^\vee } \underline \Omega^n(L^\vee) \lra \underline \Omega^n(M^\vee)$$  is surjective, where the direct sum is taken over  all split injections $L^\vee \subset M^\vee$.  This holds (almost) by definition: small Omega powers are generated by pure symbols.
	\end{dem}

\section{The Transfer, for small Omega Powers.}\label{sectranssmall}

In this section, we study the -examplary- compatibility of the Transfer with small Omega powers.\\
Here $k$ is a finite field, of cardinality $q=p^r$.

\begin{prop}

		Let $V$ be a (finite-dimensional) $k$-vector space. Let $W \subset V$  be a $k$-linear subspace, of codimension $c \geq 1$. Let $n$ be a positive integer. Then the Transfer $$\Gamma T^n_{W,V}: \Gamma^{p^{rc+n}}(V) \lra \Gamma^{p^n}(W) $$ is compatible with the formation of small Omega powers. It thus induces a $\W(k)$-linear map  $$ \underline\Omega^{rc+n}(V) \lra\underline\Omega^{n}(W). $$ 

\end{prop}

\begin{dem}
By induction (see Lemma \ref{transcomp}), we can assume that $c=1$.
The statement is clear if $V$ is two-dimensional: we know that the duality $\Delta$  on $\Gamma^{p^{r+n}}(V)$ is perfect in this case. In general, let $H \subset W$ be a $k$-rational hyperplane. Then $H$ is of codimension two in $V$. We have a commutative diagram $$ \xymatrix{ \Gamma^{p^{r+i}}(V) \ar[r]^{\Gamma T_{W,V}} \ar[d] &  W \ar[d] \\ \Gamma^{p^{r+i}}(V/H) \ar[r]^{\Gamma T^i_{W/H,V/H}} & \Gamma^{p^i}(W/H),}$$ where the vertical maps are induced by the canonical surjections. Forming the direct sum $$ \xymatrix{ \Gamma^{p^{r+i}}(V) \ar[r]^{\Gamma T^i_{W,V}} \ar[d] & \Gamma^{p^i}( W) \ar[d] \\ \bigoplus_H \Gamma^{p^{r+i}}(V/H) \ar[r]^{\oplus \Gamma T^i_{W/H,V/H}} &  \bigoplus_H \Gamma^{p^i}(W/H),}$$ over all hyperplanes $H$, yields the result. Indeed, the composite arrows vanish on $\Ker(\Delta)$ by the two-dimensional case, and we can apply  Lemma \ref{omegainj} (the vertical arrow on the right induces an injection on small Omega powers) .

\end{dem}

The following Definition thus makes sense.
\begin{defi}(Ascending and Descending Transfer for small Omega powers.)\label{transmedium}
	Let $V$ be a (finite-dimensional) $k$-vector space. Let $n$ be a positive  integer. \\Let $W \subset V$  be a $k$-linear subspace, of codimension $c \geq 1$.  The  $\W(k)$-linear map $$\underline \Omega^{rc+n}(V) \lra \underline \Omega^n(W) $$ of the preceding Proposition  is the Descending Transfer, from $V$ to 
	$W$, for small Omega powers. We denote it by $\underline {DT}^n_{W,V}$.\\
	
	Dually, let $$\pi : V \lra W $$ be a surjection between (finite-dimensional) $k$-vector spaces, such that $\Ker(\pi)$ has dimension $c$. Using Pontryagin duality, the dual of $\underline {DT}^n_{W^*,V^*}$ is a $\W(k)$-linear map $$ \underline \Omega^n(W) \lra \underline\Omega^{n+rc}(V) .$$  It is the Ascending Transfer, from $W$ to 
	$V$, for small Omega powers. We denote it by $\underline {AT}^n_{V,W}$.\\

\end{defi}

\begin{rem}\label{formulaDT}
    For $v \in V$, we have $$\underline{DT}_{W,V}((v)_{n+rc})= (v)_n$$ if $v \in W$, and  $$\underline{DT}_{W,V}((v)_{n+rc})= 0 $$ if $v \notin W$. The surprising fact is that this simple formula on symbols indeed defines a $\W(k)$-linear map!
\end{rem}

The next Proposition gives a simple formula for the ascending transfer.

\begin{prop}\label{formulaAT}
	 Let $$\pi : V \lra W $$ be a surjection between  $k$-vector spaces, such that $\Ker(\pi)$ has dimension $c$.\\
 Let $n$ be a positive integer.  Put $$X_\pi := \sum_{w \in \Ker(\pi)} (w)_{n+cr}\in \Omega^{n+cr}(V).$$  (Note that $X_\pi=0$ if $k$ has at least $3$ elements.) \\Pick an element $w \in W$. Then we have the formula $$\underline{AT}_{V,W}((w)_n)= -X_W+ \sum_{v \in \pi^{-1}(w)} (v)_{n+cr}.$$ 
\end{prop}

\begin{dem}

Denote by $$\tau: k \lra \W(k)$$ the Teichm\"uller representative, and by $$\pi^* : W^* \subset V^*$$ the inclusion of the $c$-codimensional subspace which is dual to $W$. Pick an arbitrary linear form $\phi \in V^*$. By definition of the ascending transfer, we have $$ <\underline{AT}_{V,W}((x)_n),(\phi)_{n+cr}>= p^{rc}<(x)_n, \underline{DT}_{W^*, V^* }( (\phi)_{n+cr})> \in  \W_{n+rc+1}(k).$$  Denote this quantity by $a$.\\ Note that  we used here the natural injection (Verschiebung) $$\W_{n+1}(k) \stackrel {1 \mapsto p^{rc} }\lra \W_{n+rc+1}(k).$$ Put  $$a':= <-X_\pi+\sum_{v \in \pi^{-1}(w)} (v)_{n+cr}, (\phi)_{n+cr}> \in \W_{n+rc+1}(k).$$ By perfect duality, it suffices to show that $a=a'$: indeed, the symbols $(\phi)_{n+cr}$ span the $\W(k)$-module $\underline \Omega^{n+cr}(V^*)$.

We distinguish two cases.\\
Case i): The linear form does not belong to $W^*$.\\
Case ii) The linear form $\phi$  belongs to $W^*$.\\
By Remark \ref{formulaDT}, $\underline{DT}_{W^*, V^*}( (\phi)_{n+cr})$ is equal to zero  in Case i),   and to $(\phi)_n$ in Case ii). Hence, $a=0$ in Case i), and  $$a=  p^{rc}\tau(\phi(x))^ {p^{n}}\in  \W_{n+rc+1}(k)$$  in Case ii). Clearly, we have $$a'= ( -\sum_{v \in \Ker(\pi)} \tau(\phi(v))^{ p^{n+rc}}+ \sum_{v \in \pi^{-1}(w)} \tau(\phi(v))^{ p^{n+rc}})$$ $$ =  ( -\sum_{v \in \Ker (\pi)} \tau(\phi(v))^{ p^{n}}+ \sum_{v \in \pi^{-1}(w)} \tau(\phi(v))^{ p^{n}})$$  (remember that $x^q=x$ for every $x \in k$).  Assume that we are in Case i), i.e. that $\phi$ does not vanish on $\Ker(\pi)$. Then $$ \phi_{\vert \pi^{-1}(w)}: \pi^{-1}(w) \lra k$$ and $$ \phi_{\vert \Ker(\pi)}:\Ker(\pi) \lra k$$
are both surjective maps, between $k$-affine spaces. Hence, the cardinality of the fiber of any element of $k$ by these two maps is the same, from which we get $a'=0$.
Assume now that we are in case ii). Then $\phi$ vanishes on $\Ker(\pi)$, and the $p^{rc}$ other terms occuring in the sum defining $a'$ are all equal to $\tau(\phi(w))^{ p^{n}}$. Again, we conclude that $a=a'$ and the Proposition is proved.
   
\end{dem}

\subsection{The Transfer, as a contravariant functor.}

\begin{defi}
		Let $n$ be a positive integer. 
	Let $$f :V \lra W $$ be a linear map between finite-dimensional $k$-vector spaces. Then $f$ factors canonically as the composite $$V \twoheadrightarrow V / \Ker(f)= \Im(f) \hookrightarrow W.$$  Denote by $\rho:=\delta(\Im(f))$ the rank of $f$.\\
	
	We denote by $\underline T^n(f)$ (or simply by $\underline T(f)$) the composite $$\underline \Omega^{n+r\delta(W)}(W) \stackrel{\underline{DT}_{ \Im(f),W}} \lra \underline \Omega^{n+r\rho}(\Im(f)) = \underline \Omega^{n+r\rho}(V/ \Ker(f))   \stackrel{\underline{AT}_{V,V/ \Ker(f)}} \lra \underline \Omega^{n+r\delta(V)}(V).$$
	It is the Transfer, for small Omega powers.
	
\end{defi}

\begin{prop} \label{formulatransfer}
	Let $n$ be a positive integer. Let $$f :V \lra W $$ be a linear map between finite-dimensional $k$-vector spaces.
	Put $$X_f:=  \sum_{v \in \Ker(f)} (v)_{n+r\delta(V)} \in \underline \Omega^{n+r\delta(V)}(V).$$ Note that we always have $X_f=0$ if $k$ has at least $3$ elements.\\
	Then $$\underline T(f):\underline \Omega^{n+r\delta(W)}(W) \lra \underline\Omega^{n+r\delta(V)}(V)$$ is given, on pure symbols, by the formula $$(w)_{n+r\delta(W)} \mapsto -  X_f+\sum_{v \in f^{-1}(\{ w \})} (v)_{n+r\delta(V)}.$$
	
\end{prop}

\begin{dem}
	It is easy to check that, if $g: W \lra Z$ is another linear map between finite-dimensional $k$-vector spaces, then the formula of the Proposition is true for $g \circ f: V \lra Z$, if it is true for both $f$ and  $g$. Hence, it suffices to check the formula if $f$ is injective or surjective.\\
	The formula is obviously true if $f$ is injective,  by the very definition of the descending transfer (Remark \ref{formulaDT}). That the formula is true if $f$ is surjective is the content of Proposition \ref{formulaAT}.
\end{dem}

\begin{prop}\label{smallomegafun}
	Let $n$ be a positive integer. The association $$V \mapsto \underline \Omega^{n+r\delta(V)}(V),$$ $$ f \mapsto \underline T(f)$$
	is a contravariant functor, from the category of finite-dimensional $k$-vector spaces to that of $\W(k)$-modules.
	
\end{prop}

\begin{dem}
	This follows  from the expression of $T(f)$ given in Proposition \ref{formulatransfer}.
\end{dem}

\begin{prop}
	Let $$f: V \lra W$$ be a linear map between finite-dimensional $k$-vector spaces. Then the following is true.
	
	i) If $f$ is injective, then the composite $$\underline\Omega^{n+r \delta(W)} (V) \stackrel {\underline \Omega(f)}\lra\underline \Omega^{n+r \delta(W)}(W)  \stackrel {\underline T(f)}\lra \underline \Omega^{n+r \delta(V)}(V) $$ equals $\Frob^{r(\delta(W)-\delta(V))}.$

	ii) If $f$ is surjective, then the composite $$\underline \Omega^{n+r \delta(W)} (W) \stackrel {\underline T(f)}\lra \underline \Omega^{n+r \delta(V)}(V)  \stackrel {\underline \Omega(f)}\lra \underline \Omega^{n+r \delta(V)}(W) $$ equals $\Ver^{r(\delta(V)-\delta(W))}.$
\end{prop}

\begin{dem}
	Computation, using Proposition \ref{formulatransfer}.
\end{dem}
\begin{prop} \label{formulatransfer2}
	Let $n$ be a nonnegative integer. Let $$\xymatrix{ V \ar[r]^{f_2} \ar[d]^{f_1}& W_1 \ar[d]^{g_2} \\ W_2 \ar[r]^{g_1} & Z } $$ be a cartesian diagram in the category of finite-dimensional $k$-vector spaces. Then we have $$ \underline \Omega(f_2) \circ \underline T(f_1) = \underline T( g_2)\circ \underline \Omega (g_1),$$
	
	 as $\W(k)$-linear maps $\underline \Omega^{n+r\delta(W_2)}(W_2) \lra \underline \Omega^{n+r\delta(V)}(W_1)$
	
\end{prop}

\begin{dem}
	Using Proposition \ref{formulatransfer}, we  compute that both composites are given, on symbols, by the formula $$(x)_. \mapsto \sum_{y \in W_1, g_2(y)=g_1(x)} (y)_.,$$ and the claim is proved.
\end{dem}

We conclude this section with an important Proposition. Is content is that the modules $\underline \Omega^n(V)/p^m$ are \textit{induced from dimension one}, if $n-m$ is large enough.

\begin{prop}\label{induced}
	Let $V$ be a finite-dimensional $k$-vector space, of dimension $\geq 2$. Let $n,m$ be positive integers, satisfying $$n-m \geq r(\delta(V)-1)-1.$$ For each line $L \subset V$, functoriality of small Omega powers yields a canonical map $$\underline \Omega^n(L)/p^m \stackrel {f_L} \lra    \underline \Omega^n(V)/p^m.$$ Then the  map $$f: \bigoplus_{L \in \P_k(V)} \underline \Omega^n(L)/p^m \stackrel {\sum f_L} \lra   \underline \Omega^n(V)/p^m $$ is an isomorphism of free $\W_m(k)$-modules.
\end{prop}

\begin{dem}
    The map $f$ is obviously surjective: $\underline \Omega^n(V)$ is generated by symbols $(v)_n$, for $v$ ranging through the nonzero vectors of $V$, and such a $v$ belongs to a unique $L$! For each line $L\hookrightarrow V$, we have the transfer $$\underline{DT}_{L,V}: \underline \Omega^{n}(V) \lra  \underline\Omega^{n-r(\delta(V)-1)}(L).$$ But $\underline \Omega^{n-r(\delta(V)-1)}(L)$ is a free $\W_{p^{n-r(\delta(V)-1)+1}}(k)$-module of rank one. Under our assumption on $m$, the quotient $\underline \Omega^{n-r(\delta(V)-1)}(L)/p^m$ is then a free $\W_m(k)$-module of rank one. Using the canonical isomorphism $L \simeq L^{(r)},$ we see (via the Frobenius) that there is a canonical isomorphism $$\underline \Omega^{n-r(\delta(V)-1)}(L)/p^m \simeq \underline \Omega^{n}(L)/p^m,$$ through which the transfer can be seen, modulo $p^m$, as a  $\W(k)/p^m$-linear map $$g_L: \underline \Omega^{n}(V)/p^m \lra \underline \Omega^{n}(L)/p^m,$$ sending a symbol $(v)_n$ to zero if $v \notin L$, or to $(v)_n$ if $v \in L$.\\ The sum of the $g_L$'s is a  $\W(k)/p^m$-linear map $$g:\underline  \Omega^{n}(V)/p^m \lra \bigoplus_{L \in \P(V)} \underline \Omega^{n}(L)/p^m,$$ which is seen to the inverse of $f$ (check this on pure symbols).
    
\end{dem}

\subsection{The  Integral Formulas for the Frobenius and the Verschiebung.}\label{integral}
The purpose of this section is to state and prove  Integral Formulas for the Frobenius and the Verschiebung, for small Omega powers. These formulas say  that averaging all descending (resp. ascending) transfers over all linear subspaces of a given dimension yields the Frobenius (resp. the Verschiebung). These formulas are reminiscent of motivic integration.

Recall that $\Gr_k(m,n)$, the Grassmannian of $m$ dimensional subspaces of $k^n$,  has cardinality $$ \vert \Gr_k(m,n) \vert =\frac{(q^n-1) (q^{n-1}-1) \ldots (q^{n-m+1}-1)} {(q^m-1) (q^{m-1}-1) \ldots (q-1)}.$$ In particular, all these numbers are congruent to $1$ modulo $p$.

\begin{prop}\label{frobint}(Frobenius Integral Formula.)
  Let $n$ be a positive integer. Let $V$ be a nonzero finite-dimensional $k$-vector space. Let $1 \leq m \leq \delta(V)-1$ be an integer. For each linear subspace $W \in \Gr(\delta(V)-m,V),$ denote by $$i_W: W \hookrightarrow V$$ the canonical inclusion. We have the formula $$\frac 1 {\vert  \Gr_k(m,\delta(V)-1) \vert} \sum_{W \in \Gr(m,V)}\underline \Omega(i_W) \circ \underline T(i_W) = \Frob^{mr}, $$ as $\W(k)$-linear maps $$\underline \Omega^{n+mr}(V) \lra \underline \Omega^{n}(V).$$
\end{prop}

\begin{dem}
    This is an easy computation, using the formula for  the transfer given by Proposition \ref{formulatransfer}. Pick a nonzero vector $v \in V$. For $W \in \Gr(\delta(V)-m,V)$,  the quantity $\underline \Omega(i_W) \circ  \underline T(i_W)((v)_{n+rm})$ equals zero if $v \notin W$, or $(v)_n$ otherwise. It is clear that the set of subspaces $W \in \Gr(\delta(V)-m,V)$ containing $v$ is in bijection with $\Gr(\delta(V)-m-1,V/<v>)$, hence has cardinality $\vert \Gr_k(m,\delta(V)-1) \vert$. The  formula follows.

\end{dem}

\begin{prop}\label{verint}(Verschiebung Integral Formula.)
   Let $n$ be a positive integer. Let $V$ be a nonzero finite-dimensional $k$-vector space. Let $1 \leq m \leq \delta(V)-1$ be an integer. For each linear subspace $W \in \Gr(m,V),$ denote by $$\pi_W: V \twoheadrightarrow  V/W$$ the quotient map. We have the formula $$\frac 1 {\vert  \Gr_k(m-1,\delta(V)-1) \vert} \sum_{W \in \Gr(m,V)} \underline T(\pi_W)\circ \underline \Omega(\pi_W)  = \Ver^{mr}, $$ as $\W(k)$-linear maps $$\underline \Omega^{n}(V) \lra \underline\Omega^{n+mr}(V).$$
\end{prop}

\begin{dem}
    This follows, by Pontryagin duality, from the Frobenius Integral Formula for $V^*$.

\end{dem}

	\section{Axiomatizing Hilbert's Theorem 90.}\label{axioH90}
	In this section, we explain a possible way to axiomatize the consequences of Hilbert's Theorem $90$ (Kummer theory) for the cohomology of profinite groups. The key notion here is that of cyclotomic modules, and smooth profinite groups.
	
\subsection{The notion of $n$-surjectivity.}\label{equiv}
In this section, $k$ is a perfect field of characteristic $p$, and $G$ is a profinite group.\\

\begin{defi}\label{definsurj}
Let $n \geq 1$ be an integer. Let $$f: M \lra N$$ be a  morphism of $(\W(k),G)$-modules. We say that $f$ is $n$-surjective (resp. $n$-injective) if the following holds. For every open subgroup $G' \subset G$,  the map $$f_*:H^n(G',M) \lra H^n(G',N)$$ is surjective (resp. injective).
\end{defi}

\begin{rem}
	Let $n \geq 0$ be an integer. Let $$ \E: 0 \lra A \stackrel i \lra B \stackrel \pi \lra C \lra 0$$ be an exact sequence of $(\W(k),G)$-modules. \\Then $\pi$ is $n$-surjective if and only if $i$ is $(n+1)$-injective. Indeed, using the associated long exact sequences in cohomology,  both conditions are equivalent to the vanishing of the connecting homomorphism (Bockstein) $$H^n(G',C) \lra H^{n+1}(G',A),$$ for every  open subgroup $G' \subset G$.\\
	
	\end{rem}

The next Lemma states that $n$-surjectivity is preserved by pullback and pushforward of exact sequences.

\begin{lem}\label{pushpull}
Let $n \geq 0$ be an integer. Let $$ \E: 0 \lra A \stackrel i \lra B \stackrel \pi \lra C \lra 0$$ be an exact sequence of $(\W(k),G)$-modules. Let $$f:A \lra A'$$ and $$g:C' \lra C$$ be morphisms of $(\W(k),G)$-modules. Denote by  $$ \E': 0 \lra A' \stackrel {i'} \lra B' \stackrel {\pi'} \lra C' \lra 0$$ the exact sequence $f_*(g^*(\E))$. If $\pi$ is $n$-surjective, then so is $\pi'$.
\end{lem}

\begin{dem}
	Easy diagram chase. 
\end{dem}
\begin{exo}
	Prove the preceding Lemma without using the connecting map, i.e. without invoking cohomology groups in degree $n+1$.
\end{exo}
\begin{rem}($0$-surjectivity).
 Let $$ \E: 0 \lra A \stackrel i \lra B \stackrel \pi \lra C \lra 0$$ be an exact sequence of $(\W(k),G)$-modules. Then $\pi$ is $0$-surjective if and only if it possesses a $G$-equivariant -\textit{set-theoretic}- section. This situation appears in the study of rationality questions for algebraic tori, through the use of (co)flasque resolutions of their character lattice. We refer here to the work of Endo and Miyata \cite{EM}, of Colliot-Th\'el\`ene and Sansuc \cite{CTS}, and of Voskresenskii \cite{Vos}.
\end{rem}

\begin{exo}
	It is clear that a split surjection is $n$-surjective for every $n$. In this exercise, we show that the converse implication is false in general. This exercise, though by no means easy, provides good  pratice for understanding the ideas developped in this paper.\\
	 For simplicity, we assume here that $k=\F_p$.\\
	Let $X$ be a finite $G$-set.  Let $$V \subset \F_p ^X$$ be a sub-$(\F_p,G)$-module. Put $$W:=  \F_p ^X /V.$$
	Let $M$ be a $(\Z/p^2 \Z,G)$-module, which is free of rank one as a $\Z /p^2 \Z$-module, together with an isomorphism  of $(\F_p,G)$-modules $$M/p \simeq \F_p.$$\\
	The exact sequence $$0 \lra  (pM) ^X \lra M ^X  \lra (M/pM) ^X \lra 0$$ can thus be viewed as an exact sequence  $$0 \lra  \F_p ^X \lra  M ^X  \lra \F_p ^X \lra 0.$$  Pulling it back by the inclusion $V \lra \F_p ^X$ and pushing it forward by the surjection $\F_p ^X \lra W$ yields an extension $$\E: 0 \lra W \lra E \stackrel \pi \lra V \lra 0.$$

i) Show that $\E$ is an exact sequence of $(\F_p,G)$-modules.\\
ii) Show that $\E$  depends neither on $M$ nor on the choice of the isomorphism  $M/p \simeq \F_p$ (up to isomorphism of short exact sequences of  $(\F_p,G)$-modules).\\
iii) Show that $\pi$ is $0$-surjective.\\

From now on, we assume that $G$ is "the" absolute Galois group of a field $F$ of characteristic not $p$, containing the $p$-th roots of unity for simplicity. We make no extra assumption on $X$.\\

iv) Using Kummer theory, show that $\pi$ is $1$-surjective. \textit{Hint: choose $M=\mu_{p^2}$.}\\
v) Using the Bloch-Kato conjecture, show that  $\pi$ is $n$-surjective, for every $n \geq 1$.\\
vi) Give an example (of $F$ and $X$) where $\E$ is not split.\\

\end{exo}

\subsection{Cyclotomic modules and smoothness.}
\begin{defi}
We put $$ \overline \N:=\N_{\geq 1} \cup \{ \infty\},  $$ and $\W_\infty (k)=\W(k).$
\end{defi}

\begin{defi}
Pick an element $d \in  \overline \N$.
	Let $\T$ be a  free $\W_{d+1}(k)$-module  of rank one.\\
	For $i$ a non negative integer,  we put $$\T(i)= \T^{\otimes^i_{\W_{d+1}(k) }}.  $$   For negative $i$, we  put $$\T(i)= \Hom_{\W_{d+1}(k)}(\T(-i) , \W_{d+1}(k)).$$ 
	For any   $(\W_{d+1}(k),G)$-module $M$, we put $$M(i)=\T(i) \otimes_{\W_{d+1}(k)} M, $$ the dependence in $\T$ being implicit. \\

\end{defi}

\begin{defi}(Cyclotomic module.) \label{defismoothmodule}
Let $n \geq 0$ and $d \in \overline N$ be integers. Let $\mathcal T$ be  a free $\W_{d+1}(k)$-module of rank one, endowed with a continuous $\W_{d+1}(k)$-linear action of $G$.\\ The module $\T$ is said to be $n$-cyclotomic (relative to $k$ and $G$)  if the following condition holds.\\

For every integer $s \geq 1,$ the quotient map $$\mathcal T/p^{s+1} \lra \mathcal T/p$$  is $n$-surjective. \\
If $\T(n)$ is $n$-cyclotomic for every $n \geq 1$, we shall say that $\T$ is cyclotomic (relative to $k$ and $G$). \\
The integer $d$ is the depth of $\T$. It will be denoted by $\delta(\T).$

\end{defi}

\begin{rem}
 If $\T$ has finite depth, it suffices of course to require $n$-surjectivity for $s=\delta(\T)$ in the previous Definition.
\end{rem}

\begin{rem}
The preceding Definition has an interest only if $n \geq 1$. 
\end{rem}

	\begin{rem}
	A  cyclotomic $G$-module is given by a continuous character $$\chi: G \lra \W_{\delta(T)+1}(k)^\times,$$ which shall, in our theory, play the r\^ole of the cyclotomic character in Galois theory. Indeed, we will see in a moment that Kummer theory (a consequence of Hilbert's Theorem 90) implies that the cyclotomic character at $p$ of a field of characteristic not $p$ is $1$-cyclotomic, in our sense. That it is in fact cyclotomic is the main content of the norm-residue isomorphism theorem (the Bloch-Kato conjecture).
\end{rem}

\begin{exo}
Let $\T$ be an $n$-cyclotomic $G$-module. Show that the quotient map $$\mathcal T/p^{s+1} \lra \mathcal T/p^s$$  is $n$-surjective, for every $s \geq 1$. \\ 
\end{exo}

\begin{lem}\label{extsmoothmodule}
Let $n \geq 0$ be an integer. Let $k'/k$ be a field extension.\\ Let $\mathcal T$ be a $n$-cyclotomic $G$-module over $k$. Put $$\mathcal T ':= \mathcal T \otimes_{\W(k)} \W(k').$$  Then $\mathcal T'$ is $n$-cyclotomic over $k$'.

\end{lem}

\begin{dem}
   This is clear, since $\W(k')$ is a free $\W(k)$-module.
\end{dem}

\begin{rem}
In the Definition of a cyclotomic module, it might be worth allowing the case where $\T$ is free of rank $\geq 2$: this would enable restriction of scalars for finite field extensions $k'/k$. We shall not consider this possibility here. 
\end{rem}

\begin{lem}\label{cyclo2res}
Let  $k'/k$ be an (arbitrary) field extension. Let $n$ be a positive integer. \\ Assume that $G$ is a pro-$p$-group, and that there exists an  $n$-cyclotomic $G$-module over $k'$, of depth $1$. Then, there exists  an  $n$-cyclotomic $G$-module over $k$, of depth $1$.

\end{lem}

\begin{dem}

Let $\T'$ be an  $n$-cyclotomic $G$-module over $k'$, of depth $1$; in particular, $\T'$ is a free $\W_2(k')$-module of rank one. Consider the exact sequence $$ 0 \lra \T'/p \simeq p \T' \lra \T' \lra \T'/p \lra 0.$$ Since $G$ is a pro-$p$-group,  it acts trivially on the one-dimensional $k'$-vector space $\T'/p$.  As an exact sequence of $(\W(k'),G)$ modules, the preceding sequence can thus be rewritten as  $$ \E': 0 \lra k' \lra \T' \lra k'\lra 0.$$ Choose a linear form $\phi \in \Hom_k(k',k),$ such that $\phi(1)=1$.  Denote by $i: k \lra k'$ the canonical inclusion. Then $\phi_*(i^*(\E))$ is an  sequence of $(\W(k),G)$ modules of the shape $$ \E: 0 \lra k \lra \T \lra k \lra 0,$$ where $\T$ is a free $\W_2(k')$-module of rank one, equipped with an action of $G$. From Lemma \ref{pushpull}, it follows that $\T$ is  $n$-cyclotomic, qed.
\end{dem}

\begin{rem}
The preceding Lemma can probably be generalized to cyclotomic modules of  arbitrary depth. 
\end{rem}

\begin{defi}(Smooth profinite group.) \label{defismooth}
Let $n$ and $d$ be   positive integers. The group $G$ is said to be $(d,n)$-smooth (resp. $d$-smooth)  relative to $k$,  if there exists an $n$-cyclotomic (resp. cyclotomic)  $G$-module over $k$, of depth $d$.\\
The group $G$ is said to be $n$-smooth (resp. smooth)  relative to $k$,  if there exists an $n$-cyclotomic (resp. cyclotomic)  $G$-module over $k$, of infinite depth.
\end{defi}


The fundamental example of $1$-smoothness is that of absolute Galois groups. It can be extended to a broader class of Galois groups, as follows.

\begin{prop}\label{H90}
 Let $E/F$ be an extension of fields of characteristic not $p$.  Assume that the multiplicative group $E^\times$ is $p$-divisible (i.e. the map $E^\times \stackrel {x \mapsto x^p} \lra E^\times$ is onto), and contains all $p$-th roots of unity (hence also all roots of unity of order a power of $p$). Put $$G:= \Gal(E /F)$$ and $$\mu= \varprojlim_n \mu_{p^n}(E).$$ Then $\mu$  is a  cyclotomic $G$-module (relative to $k= \F_p$).
\end{prop}

\begin{dem}
 
	By our assumptions on $E$, we have a diagram $$\xymatrix{
		1 \ar[r]   &  \mu_p(E) \ar@{=}[d]  \ar[r] & \mu_{p^{s+1}}(E) \ar [d] \ar[r]^{x \mapsto x^p} &  \mu_{p^{s}} (E) \ar[r]  \ar[d] & 1 \\ 1  \ar[r] &   \mu_{p} (E) \ar[r]  & E^\times \ar[r]^{x \mapsto x^{p}}& E^\times \ar[r] & 1 ,}$$ given by classical Kummer theory. The surjection in the lower line is obviously $1$-surjective by Hilbert's Theorem 90 for $\G_m$. By Lemma \ref{pushpull}, the surjection in the upper line is $1$-surjective as well, yielding the result.
	\end{dem}

Is is natural ask whether all $1$-cyclotomic modules occur this way. We did not investigate this question, but we expect a positive answer- possibly given by a simple construction. We now formulate it precisely.

\begin{pb}
Let $G$ be a profinite group. Let $M$ be a $1$-cyclotomic $G$-module of infinite depth, for $k=\F_p$. Let $l$ be zero or  a prime number distinct from $p$. Find an extension $E/F$ of fields of characteristic $l$, such that  the multiplicative group $E^\times$ is $p$-divisible,  contains all $p$-th roots of unity, and such that the following holds.\\

There is an isomorphism $$\phi: G \lra \Gal(E/F)$$ of profinite groups, and an isomorphism $$\psi :M \lra  \varprojlim_n \mu_{p^n}(E)$$ of $\Z_p$-modules, such that $$\psi(g.m)=\phi(g).\psi(m),$$
for all $g \in G$ and $m \in M$.

\end{pb}

\subsection{The Smoothness Conjecture.}

Recall that $G$ is a profinite group.

  \begin{defi}\label{defisymbol}
   Let $n\geq 1$ be an integer. Let  $L$ be a one-dimensional  $(k,G)$-module. \\ Cohomology classes in the image of the natural cup-product map $$ H^1(G, L)^ n \lra H^n(G,L^{\otimes n})$$ are called symbols (relative to $L$). \\If $H \subset G$ is an open subgroup , the image of a symbol in $H^n(H,L^{\otimes n})$  by the corestriction (norm) $$ \mathrm{Cor}:  H^n(H,L^{\otimes n}) \lra H^n(G,L^{\otimes n}) $$ is called an $H$-quasi-symbol (relative to $L$). \\A class which can be written as a sum $a_1 + \ldots +a_N$, where the $H_i$'s are open subgroups of $G$, and $a_i$ is an $H_i$-quasi-symbol, will be called a quasi-symbol (relative to $L$).
  \end{defi} 
  
  
  \begin{rem}\label{remsmooth}
  	Assume that  $\T$ is a $1$-cyclotomic $G$-module (of any depth). Let $n$ and $s$ be positive integers. It is  straightforward that any symbol (and hence any quasi-symbol) in  $H^n(G,\T(n)/p) $ can be lifted to a class in  $H^n(G,\T(n)/p^s) $.
  \end{rem}

  \begin{defi}\label{defiweakBK}
  	  We say that $G$ has the weak Bloch-Kato property (at $p$) if the following holds. For every integer $n\geq 1$ and for every open subgroup $H \subset G$, every   class in $H^n(H,k)$ is a  quasi-symbol (relative to $L=k$).
  \end{defi} 
  
  The following Remark is elementary, but important.
  
   \begin{rem}
   	Let $G_p$ be a pro-$p$-Sylow of $G$.
   	By the standard restriction-corestriction argument, it is straightforward to prove the following two assertions.\\
   i) The group $G$ has the weak Bloch-Kato property at $p$ if and only if $G_p$ has it.\\
   	ii)  In the preceding definition, we may replace the trivial $(k,G)$-module $L=k$ by -any- one-dimensional $(k,G)$-module $L$.
   \end{rem}

We now state the Smoothness Conjecture, which we plan to prove in a future work.

\begin{conj}\label{smoothconj}
Let $G$ be a profinite group. If $G$ is $1$-smooth, then it  has the weak Bloch-Kato property. In particular, it is smooth.

\end{conj}

\begin{rem}
The Smoothness Conjecture implies the (surjectivity part of the) Bloch-Kato conjecture, using a classical input from Milnor $K$-theory: the Lemma of Rosset and Tate. It implies that a quasi-symbol in $H^n(\Gal(F_{sep}/F), \mu_p^{\otimes n})$ is in fact a sum of symbols. Note that this Lemma, whose proof uses Euclidean division for polynomials, is of highly effective nature.
\end{rem}

We conclude this section with an instructive exercise.

\begin{exo}
	Let $k$ be an arbitrary perfect field of characteristic $p$.\\
Let $G$ be a finite $p$-group.\\
i) Assume that $\vert G \vert \geq 3.$ Show that $G$ is not smooth.\\
ii) Assume that $p=2$ and $G= \Z /2 \Z$. Show that $G$ is $1$-smooth. What are the possible cyclotomic modules for $G$?\\

\end{exo}

\subsection{Exact sequences of Kummer type.}

\begin{defi}

Let $a$ and $b$ be  positive integers.\\ The extension (of $(\W(k),G)$-modules with trivial $G$-action)   $$  0 \lra \W_{b}(k) \stackrel {1 \mapsto p^a} \lra \W_{a+b}(k) \lra  \W_{a}(k)\lra 0$$  will be called the elementary Kummer extension, of type $(a,b)$. We shall denote it by $\K_{a,b}$. The integer $a+b-1$ is called the depth of  $\K_{a,b}$. \\ We denote $\K_{d,1}$ simply by $\K_d$. It is the elementary Kummer extension, of depth $d$.
\end{defi}
\begin{rem}\label{KDpull}

Let $a$ and $b$ be  positive integers.\\ Then Pontryagin duality exchanges $\K_{a,b}$ and $\K_{b,a}$.\\
The diagram  $$\xymatrix{0  \ar[r] & \W_{b}(k) \ar[r] \ar@{=}[d] &  \W_{a+b}(k) \ar[r] \ar[d]^{1 \mapsto p} & \W_{a}(k)  \ar[r] \ar[d]^{1 \mapsto p} & 0 \\ 0 \ar[r] & \W_{b}(k) \ar[r]  &  \W_{a+b+1}(k) \ar[r]  &  \W_{a+1}(k)   \ar[r] & 0 }$$  is clearly a pullback diagram. This shows that $\K_{a,b}$ is a pullback of $\K_{a+1,b}$. Dually, $\K_{a,b}$ is a pushforward of $\K_{a,b+1}$.
\end{rem}
\begin{defi}\label{Kummertype}
Let $d$  be a positive integer.\\ 
We denote by $ \K_d(G)$ the smallest class of   extensions   $$0 \lra A \lra B \lra C \lra 0 $$ of $(\W(k),G)$-modules,  containing the  extension $$ 0 \lra k  \lra \W_{d+1}(k) \lra  \W_{d}(k)\lra 0,$$  and stable by the following operations.\\

i) Arbitrary finite direct sums, pullbacks and pushforwards (of extensions of $(\W(k),G)$-modules).\\
ii) Induction from open subgroups: if $H \subset G$ is an open subgroup, and if $$0 \lra A \lra B \lra C \lra 0 $$ belongs to  $ \K_d(H)$, then $$0 \lra \Ind_H^G(A) \lra \Ind_H^G(B) \lra \Ind_H^G(C) \lra 0 $$ belongs to  $ \K_d(G)$.\\
iii) Composition (on the right): if  $f: A \lra B \lra 0$ and $g: B \lra C \lra 0$ are the epimorphisms of  extensions belonging to $\K_d(G)$, then $$ 0 \lra \Ker(f \circ g) \lra A \stackrel {f \circ g} \lra C \lra 0$$ belongs to $\K_d(G)$.\\

Extensions belonging to $\K_d(G)$ are said to be of Kummer type, of depth $\leq d$. \\

We put $$\K(G) = \bigcup_{d \geq 1} \K_d(G). $$ Extensions belonging to $\K(G)$ are said to be of Kummer type.\\
Epimorphisms or monomorphisms fitting into an exact sequence of Kummer type, will also be called of Kummer type.
\end{defi}

\begin{rem}\label{Kabin}
Using Remark \ref{KDpull} and property iii) of the previous Definition, we see that the Kummer extension of type $(a,b)$ belongs to  $\K_{a+b-1}(G)$, for every positive integers $a$ and $b$.
\end{rem}

\begin{lem}
Let $d$ be a positive integer. The following assertions are true.\\

a) The class $\K_d(G)$ is  stable by composition on the left. In other words, if  $f: 0 \lra A \lra B $ and $g: 0 \lra B \lra C $ are the monomorphisms of  extensions belonging to $\K_d(G)$, then $$ 0 \lra  A \stackrel {f \circ g} \lra C  \stackrel h \lra \Coker {f \circ g} \lra 0$$ belongs to $\K_d(G)$ as well.\\
b) The class $\K_d(G)$ is stable under Pontryagin duality.

\end{lem}

\begin{dem}
Point a) is obviously Pontryagin dual to point iii) of the definition of $\K_d(G)$, whereas point i) and ii) are self dual (Pontryagin duality exchanges pullbacks and pushforwards, and Pontryagin duality commutes to induction from open subgroups). By Remarks \ref{KDpull}   and \ref{Kabin},  we thus  see that b) follows from a).  We now prove a).\\
Forming the pullback of $$0 \lra B \stackrel g \lra C \lra C/B \lra 0 $$ by the natural quotient map $$ C/A  \lra C/B$$ yields the diagram $$\xymatrix{0  \ar[r] & B  \ar[r] \ar@{=}[d] & C \bigoplus B/A \ar[r]^{can} \ar[d] & C /A   \ar[r] \ar[d] & 0 \\ 0 \ar[r] & B  \ar[r]  & C  \ar[r]  & C /B    \ar[r] & 0, }$$ where $can$ is the sum of the canonical inclusion and of the canonical surjection.  By point i) of the definition of $\K_d(G)$ (for pullbacks), the upper row is of Kummer type. By assumption on $f$, and by point i) again (but for direct sums), the natural surjection $$  C \bigoplus B \lra C  \bigoplus B/A$$  belongs to $\K_d(G)$ as well. By point iii), we  see that the composite surjection $$ C \bigoplus B \lra C  \bigoplus B/A \lra C/A  $$  belongs to $\K_d(G)$. Noting that it equals the composite $$C \bigoplus B \stackrel \Sigma \lra  C \stackrel h \lra C/A,$$ we finally conclude (using point i), for pushforwards this time) that $h$ belongs to $\K_d(G)$.
\end{dem}

\begin{rem}
A surjection of Kummer type  can be intuitively thought of as 'a surjection through which cohomology classes can be lifted', in the spirit of model categories. This is made precise in the Proposition below.
\end{rem}

\begin{prop}\label{Kummersurj}
Let $n$ and 
$d$ be  positive integers. Assume that $G$ is $(d,n)$-smooth. Let $\T$ be an $n$-cyclotomic $G$-module, of depth $d$. Let $$ 0 \lra A \lra B \stackrel f \lra C \lra 0$$ be a exact sequence of Kummer type, of depth $\leq d$.  Then the sequence $$ 0 \lra A(n) \lra B(n) \stackrel {f(n)} \lra C(n) \lra 0$$ is $n$-surjective.

\end{prop}
\begin{dem}
This property holds, by the definition of an $n$-cyclotomic $G$-module of depth $d$, for the exact sequence  $\K_d$. It remains to be checked that it is stable under the operations i), ii), iii) of the definition of an exact sequence of Kummer type. For point i), use  Lemma \ref{pushpull}. For point ii), use the definition of $n$-surjectivity. For point iii), there is almost nothing to do.
\end{dem}

\begin{lem}\label{crosslemma} (The Cross Lemma.)\\
Let $$\xymatrix{& 0 \ar[d] &   0 \ar[d] &  0 \ar[d] \\ 0 \ar[r]  & A_1 \ar[r] \ar[d]& B_1 \ar[r]^h \ar[d] &  C_1 \ar[r] \ar[d] & 0 \\0 \ar[r]  & A_2 \ar[r] \ar[d]^f & B_2 \ar[r] \ar[d]^g &  C_2 \ar[r] \ar[d] & 0 \\ 0 \ar[r]  & A_3 \ar[r] \ar[d] & B_3 \ar[r] \ar[d] & C_3 \ar[r]  \ar[d]& 0 \\ & 0 &   0  &  0}$$  be a commutative diagram of $(\W(k),G)$-modules, with exact rows and columns. Assume that $g$ and $h$ are of Kummer type, of depth $\leq d$. Then so is $f$.

\end{lem}

\begin{dem}
Exercise in homological algebra for the reader.
\end{dem}
\section{About Hilbert's Theorem 90.}\label{Hilbert}
 	The authors now want make a brief digression, to stress  the importance of  Hilbert's Theorem $90$. It is, by the way, the favorite Theorem of the second author of this paper, who is a big fan of descent statements. The theory developped  here shows that this Theorem,  contrary to what one could expect, is perhaps -the- key ingredient to a 'short' proof of the Bloch-Kato conjecture, over a field $F$ of characteristic not $p$. Indeed, the Stable Lifting Theorem in the next section,  is the starting point of a machinery that applies Hilbert's Theorem 90 for $\G_m$ ceaselessly, not only to the base field $F$ itself, but also to a vast amount of finite extensions of $F$.  Furthermore, we are tempted to make the following analogy.  Adopting the point of view of Grothendieck's descent theory, the main content of (the classical version of) Hilbert's Theorem $90$ for $\GL_n$ is to convert into  cohomological information ($H^1(F,\GL_n)=1$) the highly non canonical fact that, over a field, every vector space possesses a basis. This is perfectly in the spirit of this paper: studying intrinsic  properties of divided powers for modules over Witt vectors. Choosing a basis  for these is often misleading- except, perhaps, in some proofs.\\
 	
 	Since Hilbert's Theorem 90 is central in this paper, we decided to discuss, in this section, some of its most significant algebraic incarnations. They are probably folklore for some mathematicians. They make precise the following philosophical statement: two finite linear data over a local ring $A$, which become isomorphic after a faithfully flat extension of $A$, are already isomorphic over $A$.  Before proceeding any further, we wish to remind the reader that Hilbert's Theorem 90 (for $\G_m$) is actually due to Kummer for cyclic field extensions, and that its generalization to arbitrary Galois extensions is due to Noether.\\
 
 We begin by an elementary correspondence, which is the set-theoretic version of the equivalence between line bundles and $\G_m$-torsors.
 	
\begin{lem}\label{gmtors}
	Let $S$ be a (not necessarily commutative, unital) ring. Then there is an equivalence between (left) $S$-modules $L$ which are free of rank one, and sets $X$ equipped with a (left) simply transitive action of the multiplicative group $S^\times$. In one direction, it is given by associating to $L$ its set of generators: $$L \mapsto X:=\{x \in L, L=Sx \}.$$ In the other direction, it is given by  $$ X \mapsto (S \times X) /S^\times,$$ where we mod out the free action of $S^\times$ given by $$\lambda.(s,x)= (s\lambda^{-1} , \lambda.x).$$
\end{lem}

 	\begin{dem}
 		This is clear.
 	\end{dem}
 	\begin{lem}\label{H90lem}
 		Let $A$ be a Noetherian local ring. Let $A'/A$ be a  faitfully flat extension of commutative rings, and let $S$ be a  $A$-algebra, which is finite as an $A$-module.   Let $M$ be an $S$-module.  Put $S':=S \otimes_A A'$ and $M':=M \otimes_A A'$.  If $M'$ is a free $S'$-module of rank one, then $M$ is a free $S$-module of rank one.
 	\end{lem}

 	\begin{dem}
 		Let $\kappa$ be the residue field of $A$. Put $\overline M := M \otimes_A \kappa$, $\overline S:= S \otimes_A \kappa$. Assume that $\overline M$ is a free $\overline S$-module of rank one. Then, by Nakayama's Lemma, the lift of a generator of $\overline M$ (as an $\overline S$-module) to $M$ will be a generator of $M$ (as an $S$-module). Hence, we are reduced to the case where $A$ is a field. Another similar application of Nakayama's Lemma shows that we may mod out the Jacobson radical of $S$, and assume that $S$ is a semi-simple algebra. Hence, $S$ is isomorphic direct product of matrix rings of the form $M_{n_i}(D_i)$, where $D_i$ are division $A$-algebras. We may thus assume that $S=M_n(D)$  for $D$ a division $A$-algebra.  But then, by Morita equivalence, $M$ is isomorphic to a sum of $r$ copies of the simple module $D^n$. Since $M \otimes_A A'$ is free of rank one as an $S'$-module, we must have $r=n$ by dimension count, and $M$ is free of rank one. 
 		
 	\end{dem}
 	
 	\begin{rem}
 	Assume, in what precedes, that $S$ is finite and locally free as an $A$-module.	 Then, the  group of invertible elements in $S$ is representable by the affine $A$-group scheme $GL_1(S)$ (which  is an open subscheme of $\A_A(S)$), and Grothendieck's descent theory asserts that $GL_1(S)$-torsors over $\Spec(A)$, for the fppf topology, correspond to $S$-modules $M$ as  in the previous Lemma. We thus get $$H^1(\Spec( A), GL_1(S)) = \{*\},$$ where  cohomology is taken with respect to the fppf topology. This statement is known as Grothendieck-Hilbert's Theorem 90.\\
 	\end{rem}
 	
 	\begin{prop}\label{H90gen}
 		Let $A$ be a Noetherian local ring. Let $A'/A$ be a  faitfully flat extension of commutative rings, and let $R$ be an $A$-algebra. Let $N$ be an $R$-module, which is finite  as an $A$-module. Let $M_1$, $M_2$ be two $R$-submodules of $N$.  Put $R'=R \otimes_A A'$, $N'=N \otimes_A A'$, $M_1'=M_1 \otimes_A A'$ and $M_2'=M_2 \otimes_A A'$. Assume there exists $f' \in GL_{R'}(N')$ such that $f'(M_1')=M_2'$. Then there exists $f \in GL_{R}(N)$ such that $f(M_1)=M_2$.
 	\end{prop}
 	
 	\begin{dem}
 		
 		Put $$S:=\{ f \in \End_R(N), f(M_1) \subset M_1\};$$ it is 	an $A$-algebra. It is a subalgebra of $\End_A(N)$. Writing $N$ as a quotient of a free module $A^n$, $\End_A(N)$ then occurs as a sub-$A$-module of $N^n$, which is finite by assumption. Hence, $S$ itself is a finite $A$-module. Put $S':=S \otimes_A A'$. By faithful flatness, we get that the canonical morphism $$S' \lra \{ f' \in \End_{R'}(N'), f'(M'_1) \subset M'_1\}$$ is an isomorphism. The set $$X:=\{f' \in  \GL_{R'}(N'), f'(M'_1)=M'_2\}$$
 		is endowed with a simply transitive action of the multiplicative group $S'^\times$. As such (see Lemma \ref{gmtors}), it canonically corresponds to a free $S'$-module  of rank one $M'$, given by the set-theoretical formula $$M'=(X \times S') /{S' }^\times. $$  But the $S'$-module $M'$, viewed as an $A'$-module, is endowed with a canonical descent data for the faithfully flat morphism $A'/A$. By descent, we get an $A$-module $M$, which is in fact a locally free $S$-module of rank one. To prove the Proposition is equivalent to proving that $M$ is actually a free $S$-module of rank one (to give a generator of the $S$-module $M$ is equivalent to giving  $f \in GL_{R}(N)$ such that $f(M_1)=M_2$). We conclude by applying Lemma \ref{H90lem}.
 		
 	\end{dem}
 	
 	\begin{coro}
 		Let $A$ be a Noetherian local ring. Let $A'/A$ be a  faitfully flat extension of commutative rings, and let $R$ be an $A$-algebra. Put $R':=R \otimes_A A'$. Let $N,M$ be two $R$-modules, one of which is finite as an $A$-module. Assume that  $M \otimes_A A'$ and $N \otimes_A A'$ are isomorphic as $R'$-modules. Then $M$ and $N$ are isomorphic as $R$-modules.
 			\end{coro}
 			
 			\begin{dem}
 		 To see this, just apply the Proposition to $M$ and $N$, viewed as $R$-submodules of $M \bigoplus N$.  
 		 \end{dem}
 		 
 		 \begin{rem}
 	Specializing to linear representations, we get the following statement. Two finite-dimensional linear representations of an abstract group $G$ over a field $F$, which become isomorphic over an extension $E/F$, are already isomorphic over $F$. Note that this holds, in particular, in the modular case (i.e.  where $F$ has characteristic $p$ and $G$ is a finite $p$-group).
 	\end{rem}

 	\begin{rem}
 		In all what precedes, the Noetherian assumptions may probably be dropped. They are  here to simplify the proofs.
 	\end{rem}

\section{The Stable Lifting Theorems.}\label{stablesec}

In this section, $k$ is a finite field of cardinality $q=p^r$ and $G$ is a profinite group. 
 We use small Omega powers here: we believe they are better behaved than medium (or big) Omega powers, for applications in Galois theory. \\

Let us explain how we intend to apply small Omega powers to prove Lifting Theorems in Galois cohomology- with an explicit proof of the Bloch-Kato conjecture as a main motivation. \\
Assume that $G$ is $s$-smooth, and denote by $\T$  a fixed $s$-cyclotomic $G$-module, of infinite depth. \\

Let $V$ be a $(k,G)$-module, and let $n$ be a nonnegative integer.
Almost by definition of smoothness, we have that the (twist of) Frobenius $$\Frob(s): \underline \Omega^{n+1}(V)(s) \lra \underline \Omega^{n}(V^{(1)})(s)$$ is $s$-surjective, \textit{if $V$ is one-dimensional}.  Indeed, we can reduce to the case where $G$ is a pro-$p$-group. For such a $G$, any one-dimensional $V$ is isomorphic to the trivial $(k,G)$-module $k$, and the statement becomes nothing but the definition of $s$-surjectivity. It is then legitimate to wonder whether the same holds for an arbitrary $(k,G)$-module $V$. The next exercise shows that it is not the case in general. This is deeply related to the notion of $R$-equivalence, due to Manin.

\begin{exo}\label{Requi}
Let $F$ be an infinite field of characteristic not $p$, with separable closure $F_{sep}/F$. Assume for simplicity that $F$ contains the $p^2$-th roots of unity. We denote by $H^1_{et}(.,.)$ the first \'etale cohomology groups.\\
In this exercise and in this exercise only,  $s=1$, $G:=\Gal(F_{sep}/F)$ and $$\T:=\Z /p^2 \Z.$$ By Proposition \ref{H90}, we know that $\T$ is a $1$-cyclotomic $G$-module. \\In view of the assumptions, we can remove all twists (by Frobenius, and by roots of unity).\\
 Let $V$ be a finite commutative algebraic $F$-group of multiplicative type, killed by $p$. We shall identify $V$ with the $(\F_p,G)$-module $V(F_{sep})$. \\

A cohomology class $c \in H_{et}^1(\Spec(F),V) =H^1(G,V)$ is said to be (elementarily) $R$-trivial if the following holds. There exists an open subvariety $U \subset \A^1_F$, containing $0$ and $1$, and a class $$C \in  H_{et}^1(U,V),$$ whose specialization at $0$ (resp at $1$) is trivial (resp. equals $c$).\\

1) Assume that $V$ has dimension one. Show that every class in $c \in H^1(G,V)$ is $R$-trivial.\\
\textit{Hint: reduce to the case $V=\mu_p$, and use Kummer theory.}\\
2) Assume that the map  $$\Frob:\underline  \Omega^1(V) \lra V$$ is $1$-surjective, for every $V$ as above. Show that every element in  $H^1(G,V)$ would then be $R$-trivial.\\
\textit{Hint: induction on the dimension of $V$, using the Frobenius Integral formula \ref{frobintegral}, and Shapiro's Lemma.}\\
3) Using the work of Colliot-Th\'el\`ene and Sansuc (\cite{CTS}),  give an example of  a field $F$ and of a $V$ as above, such that not every class in $H_{et}^1(\Spec(F), V)$ is $R$-trivial. Conclude.\\
4) Let $c\in H^1(G,V)$ be a Galois cohomology class, which is $R$-trivial. It is true that $c$ is in the image of $$ \Frob_*:H^1(G, \underline \Omega^1(V) )\lra H^1(G,V) ~~? $$
\textit{Hint: we don't know the answer...}
\end{exo}

The following  Lifting Theorem  brings hope that our approach will shortly lead to an 'elementary' proof of the Bloch-Kato conjecture. Note that, if this Theorem was true for $n=0$ and $V$ arbitrary,  Galois cohomology would be much simpler, and the Bloch-Kato conjecture would follow  quite easily.

\begin{thm}\label{StableLift}(First Stable Lifting Theorem.)

	Let $V$ be a $d$-dimensional $(k,G)$-module. Let $n$ be a positive integer, with $$n \geq r(d-1)-1.$$ Let $t$ be an arbitrary positive integer.\\  
	The following assertions are true.\\
	
	1) The Frobenius homomorphism  $$ \Frob^t: \underline \Omega^{n+t}(V) \lra  \underline \Omega^{n}(V^{(t)})$$  is of Kummer type, of depth $\leq n+t$. \\
2) If 	$\T$  is an $s$-cyclotomic $G$-module of  depth $n+t$, then  the twist  $$\Frob(s): \underline \Omega^{n+t}(V)(s) \lra  \underline \Omega^{n}(V^{(t)})(s)$$ is $s$-surjective. 
	
\end{thm}

\begin{dem}
We prove part 1); part 2) will follow by Proposition \ref{Kummersurj}.\\

	The assertions are clear if $V$ is one-dimensional, by the very definition of surjections of Kummer type. We can thus assume that $d  \geq 2$.\\
	 Now, consider the commutative diagram $$ \xymatrix{0 \ar[r] & \underline \Omega^{n+t} (V)[p^t] \ar[d]^g \ar[r]  & \underline \Omega^{n+t} (V) \ar[r]^{\Frob^t} \ar[d]&  \underline \Omega^{n} (V^{(t)})  \ar[d]  \ar[r] & 0    \\ 0 \ar[r] &  \bigoplus_H \underline \Omega^{n+t} (V/H)[p^t] \ar[r] &  \bigoplus_H \underline  \Omega^{n+t} (V/H) \ar[r] ^{\oplus \Frob_{V/H}} &  \bigoplus_H \underline \Omega^{n} ((V/H)^{(t)})  \ar[r] & 0  ,} $$   where the direct sums are taken over all $k$-rational hyperplanes $H \subset V$, and the vertical arrows are obtained by functoriality from the quotient maps $V \lra V/H$. The Pontryagin dual of  the map $g$ is the map $f$ of Proposition \ref{induced} (applied to $V^*$, $n+t$ and $m=t$). Since $n \geq r(d-1)-1$, the same Proposition asserts that $f$, hence $g$, is  an isomorphism. \\
Point 1) now follows, from the definition of a surjection of Kummer type: the lower row is induced from dimension one (the open subgroups involved are the stabilizers of hyperplanes of $V$).
\end{dem}

The next Theorem is more precise: it asserts that the cohomology of a smooth profinite group, with values in a (twist of a) $G$-module of the type  $\underline \Omega^n(V)$, is \textit{induced from dimension one}, if $n$ is large enough.
\begin{thm}(Second Stable Lifting Theorem.) \label{Lifting2}
	Let $V$ be a $d$-dimensional $(k,G)$-module, with $d \geq 2$.\\  Let $n$ be a positive integer, with $$n \geq r(d-1)-1.$$ 
 For each line $L \in \P(V),$ functoriality of Omega powers yields a canonical injection $$ \underline \Omega^n(L) \lra  \underline \Omega^n(V).$$ Then the following assertions are true.\\
 
 1) The (surjective) $G$-equivariant map $$\bigoplus_{L \in \P(V)} \underline \Omega^n(L)  \lra  \underline\Omega^n(V)$$ is  
 of Kummer type (of depth $\leq n+r(d-1)$).\\
 2) If 	$\T$  is an $s$-cyclotomic $G$-module of  depth $n+r(d-1)$, then the twisted $G$-equivariant map
 $$\bigoplus_{L \in \P(V)}  \underline\Omega^n(L)(s)  \lra  \underline \Omega^n(V)(s)$$ is $s$-surjective.
	
\end{thm}

\begin{dem}
We prove part 1); part 2) will follow by Proposition \ref{Kummersurj}.\\
By the previous Theorem, for $t=r(d-1)$, we see that the map $$\Frob^{r(d-1)}:  \underline\Omega^{n+r(d-1)}(V) \lra  \underline\Omega^{n}(V)$$ is of Kummer type, of depth $\leq n+r(d-1)$. By the Frobenius Integral formula \ref{frobint}, this map factors as $$   \underline\Omega^{n+r(d-1)}(V) \lra \bigoplus_{L \subset V}   \underline\Omega^{n}(L) \lra \underline \Omega^{n}(V),$$ where the direct sum is taken over all lines $L \subset V$, and the first map is the sum of the Transfers, for all inclusions $L \subset V$.  The claim follows, by the definition of a surjection of Kummer type.
\end{dem}

The preceding Theorem has a very concrete consequence, for cohomology classes with values in two-dimensional Galois representations, over $\F_p$.

\begin{coro}\label{concrete}
   Let $F$ be a field of characteristic not $p$. Denote by $F_{sep}/F$ a separable closure of $F$.  Let $V$ be a two-dimensional Galois representation of $G:=\Gal(F_{sep}/F)$ over $\F_p$. For each $v \in V$, denote by $G_v \subset G$ the stabilizer of $v$.  Denote by $\phi_v$ the composite  $$ H^1(G_s, \F_p) \lra   H^1(G_s, V)\stackrel {\Cor_{G_s}^G} \lra H^1(G, V),$$ where the first map is induced by functoriality from the $G_s$-equivariant map $$ \F_p \stackrel {1 \mapsto v} \lra V.$$
   Then the map $$\bigoplus_{v \in V} H^1(G_s, \F_p) \stackrel {\oplus \phi_v} \lra   H^1(G, V)$$ is surjective.\\
   In other words, classes in $H^1(G,V)$ are 'induced from dimension one'.

\end{coro}
\begin{dem}

The statement  is  the concrete form of Theorem \ref{Lifting2}, for $\mathcal T= \mu_{p^2}$ (which is $1$-smooth by Kummer theory), $s=1$, $r=1$, $d=2$ and $n=0$.
\end{dem}

\begin{exo}

Show that, for $p=2$, the preceding Corollary is true for \textit{any} profinite group $G$, and any two-dimensional $(\F_2,G)$-module.
\end{exo}

\section{An application to $p$-adic deformation theory.}

We finish this paper with an application of our point of view to $p$-adic deformation theory. It  uses very few of the theory of Omega powers that we have explained before. In truth, it uses only the divided power functor $\Gamma_\Z ^p$, for $\Z$-modules of $p$-primary torsion.\\
Recall that we have seen that $p$-typical Witt vectors of level $n$ may be defined, in a pretty elementary way, as a quotient of a divided power module over $\Z$ (cf. Proposition \ref{newwitt}).\\
In this section, we give another elementary application of  this point of view, to the problem of lifting $\F_p$-algebras to flat $(\Z/p^2 \Z)$-algebras. 

\subsection{The case of perfect $\F_p$-algebras.}
\begin{defi}\label{pnice}
A $p$-nice ring is a commutative ring $R$,  satisfying the following three conditions.\\

 \noindent i) The ring $R$ is $p$-adically complete (i.e. complete with respect to the ideal $pR \subset R$).\\
\noindent ii) $p$ is not a zero-divisor in $R$, i.e. $R$ is torsion-free.\\
\noindent iii) The $\F_p$-algebra $R/p$ is perfect, i.e.  its Frobenius  $$\Frob: R/p \lra R/p$$ $$x \mapsto x^p$$ is an isomorphism.

\end{defi}

Recall the following result in $p$-adic deformation theory, well-known to experts.

\begin{prop}\label{deformation}

The reduction functor $$ \Phi: \{p-\mathrm{nice \> rings} \} \lra \{\mathrm{perfect \>} \F_p-\mathrm{algebras} \}, $$ $$ \mathcal A \mapsto \mathcal A/p$$

is an equivalence of categories (morphisms being ring homomorphisms on both sides).

\end{prop}

\begin{dem}
    The usual proof of this Proposition is through Illusie's cotangent complex, as explained in the work of Scholze (\cite{S}, Theorem 5.11 and Theorem 5.12).
\end{dem}

It is possible to give a direct elementary proof of the preceding Proposition, using merely the functor  $\Gamma^p_\Z$. Rather than giving more details, we prefer to state and prove a statement  which is much more general- but, for the time being, for mod $p^2$ liftings only.

   \subsection{Descent for the arrow $\Z/ p^2 \Z \lra \Z/p \Z$.}
   
   In  Proposition \ref{deformation}, 
  the crucial assumption is that the Frobenius map of $R/p$ is surjective. The fact that it is injective (i.e. that $R/p$ is reduced) is secondary. The authors believe that it is important to adapt this Proposition to the case of (possibly non reduced) $\F_p$-algebras whose Frobenius map is surjective. In a naive sense, these algebras are 'Frobenius-smooth' objects (existence of lifts by Frobenius), whereas perfect $\F_p$-algebras are "Frobenius-\'etale" objects (existence and uniqueness of lifts by Frobenius). To tackle this question, today's trend is  to use Scholze's theory of perfectoid spaces, in which these  algebras typically occur. We here suggest a first step towards an alternate approach, in the spirit of our paper: Proposition \ref{wrinklelift}. Its content is that deforming an $\F_p$-algebra $A$, whose Frobenius is surjective, to a flat $(\Z/p^2 \Z)$ -algebra $\mathcal A$, is equivalent to endowing the kernel of the Frobenius of $A$ with a partial (level $p$) divided power operation $\gamma_p: I \lra A$.   Moreover, the Frobenius of $A$ lifts to $\mathcal A$ if, and only if, $\gamma_p$ takes values in $I$. 
  
  \begin{defi}
Let $R$ be a commutative ring in which $(p-1)!$ is invertible. For any integer $i$ with $0 \leq i \leq p-1$ and $x\in R$, we set 
$$\gamma_i(x):= \frac  1 {i !} x^i .$$
 
\end{defi}

  \begin{defi}
  A $2$-wrinkled ring (relative to $p$) is the data of a pair $(A, \gamma_p)$,  consisting of  an $\F_p$-algebra $A$, whose Frobenius is surjective, and a map $$\gamma_p:\Ker(\Frob_A) \lra A,$$  such that the relations $$\gamma_p (ax)=a^p \gamma(x)$$ and $$ \gamma_p(x+x')=\sum_{i+i'=p} \gamma_i(x)\gamma_{i'}(x'),$$  hold for all $a \in A$ and all $x,x'$ in $\Ker(\Frob)$. \\
  
  $2$-wrinkled rings obviously form a category: a morphism $(A,\gamma_p) \lra (A',\gamma'_p)$ is a ring homomorphism $\phi: A \lra A'$, such that $\gamma'_p \circ \phi=\phi  \circ \gamma_p.$
  \end{defi}
  \begin{rem}
  Let $(A, \gamma_p)$ be a $2$-wrinkled ring. Put $I:=\Ker(\Frob_A)$.\\ 
  It is not hard to see that $\gamma_p$ vanishes on $I^2$, and that it is in fact given by a unique \textit{polynomial law of $A$-modules}, which is homogeneous of degree $p$, from  $I/I^2$ to $A$.
 \end{rem}
  \begin{defi}
  Denote by $\mathcal F_2$ the forgetful functor, from the category of $2$-wrinkled rings to that of $\F_p$-algebras.\\
  A $2$-liftable ring (relative to $p$) is an $\F_p$-algebra $A$, which lies in the essential image of $\mathcal F_2$.
  \end{defi}
  
  \begin{rem}
    Note that a non-reduced (i.e. non-perfect) $2$-liftable ring is not Noetherian (a surjective endomorphism of a Noetherian ring is an isomorphism...).
\end{rem}

   \begin{defi}\label{pflat}
 A $2$-flat ring (relative to $p$) is a commutative $(\Z /p^2 \Z)$-algebra $R$,  satisfying the following conditions.\\

\noindent  i) $R$ is a flat (i.e. free) $\Z / p^2 \Z$-module. \\
\noindent ii) The  Frobenius map of $R/p$ is surjective.\\ 
The $2$-flat rings form a category, with morphisms been ring homomorphisms.\\

\end{defi}

  Let $\mathcal A$ be a $2$-flat ring. Put $$A:=\mathcal A/p . $$
  
  We shall now see that $A$ can be given the structure of a $2$-wrinkled ring in a canonical way. \\
  
  Put  $$ I:=\Ker(\Frob_A).$$ 

   Let $x$ be an element of the ideal $I$. Let $X \in \mathcal A$ be any lift of $x$. The quantity $X^p \in \mathcal A$ does not depend on the choice of $X$ (cf. Lemma \ref{Teichpoly}). By assumption, there exists $Y \in \mathcal A$ such that $$ X^p =pY.$$ Denote by $y \in A$ the reduction of $Y$. Since $\mathcal A$ is a free $(\Z /p^2 \Z)$-module, $y$ does not depend on the choice of $y$. We then put $$\gamma_p(x):= \frac 1 {(p-1) !} y= -y \in A.$$ It is not hard to see that $(A,\gamma_p)$ is a $2$-wrinkled ring. 

  We have in fact built a functor $$\Psi_2: \{ 2- \mathrm{flat \> rings}\} \lra  \{ 2- \mathrm{wrinkled  \> rings}\},$$ by the formula $$ \Psi_2( \mathcal A ):= (A, \gamma_p).$$ 
   
   We can then prove the following result, which is a descent statement for the quotient map $$\Z/p^2 \Z \lra \Z/ p \Z.$$ It  generalizes Proposition \ref{deformation} (for mod $p^2$ deformations). Note that this Proposition may seem strange at first glance: the quotient map $\Z/p^2 \Z \lra \Z/ p \Z$ is not quite flat, and descent statements in algebraic geometry are often the privilege of faithfully flat morphisms. However, the categorical data that allows descent here is not at all the usual one (it is non-linear).  
   
   \begin{prop}\label{wrinklelift}
   The functor $\Psi_2$ is an equivalence of categories. \\In particular, every $2$-wrinkled ring admits a unique lift to a $2$-flat ring.
   \end{prop}
   
   \begin{dem}
   
   Let $\mathcal A$ be a $2$-flat ring.  Put $A=\mathcal A/p$. Denote by $$\gamma_p : \Ker (\Frob_A) \lra A $$ the $p$-th divided power operation constructed above. \\ We know, by Lemma \ref{Teichpoly}, that the polynomial law of $\Z$-modules $$ \mathcal A \lra \mathcal A,$$ $$X \mapsto X^p $$ factors through the quotient map $\pi: \mathcal A \lra A,$ yielding by the universal property of divided powers a group homomorphism
  $$ F: \Gamma^p_{\Z}(A) \lra \mathcal A,$$ $$[\pi( X)]_p \mapsto X^p .$$  Note that pure symbols generate $ \Gamma^p_{\Z}(A) =  \Gamma^p_{\Z/p^2 \Z}(A)$, by Lemma \ref{symbolsgen}. The group $\Gamma^p_{\Z}(A)$   bears a natural  ring structure (formula on pure symbols: $[x]_p[y]_p=[xy]_p$), for which $F$ is a ring homomorphism. Since the Frobenius of $A$ is surjective, it is easily checked that $F$ is onto. Let $X_1, \ldots, X_m$ be elements of $\mathcal A$, with reductions $x_1, \ldots, x_m$ in $A$.  One has $$ [x_1]_p+ \ldots + [x_m]_p \in \Ker(F)$$ if and only if $$X_1^p+ \ldots +X_m^p=0 \in \mathcal A. $$ Rewrite this equality as $$(X_1 + \ldots+ X_m)^p = p \sum_{a_1, \ldots, a_m} C_{a_1, \ldots, a_m} X_1^{a_1} \ldots X_m^{a_m} \in \mathcal A, $$  where the sum ranges over all proper  partitions of $p$, i.e. partitions $$p=a_1+ \ldots+a_m,$$ with $0 \leq a_i \leq p-1$ for all $i$, and where $$C_{a_1, \ldots, a_m}:= \frac 1 p { p \choose {a_1, \ldots, a_m}}  \in \N. $$ 
  
  By the very definition of $\gamma_p$, this is equivalent to the combination of the two equalities $$ x_1+  \ldots +x_m \in \Ker(\Frob_A)$$ and $$-\gamma_p(x_1 + \ldots+ x_m) = \sum_{a_1, \ldots, a_m} C_{a_1, \ldots, a_m} x_1^{a_1} \ldots x_m^{a_m} \in A. $$   The isomorphism $ \Gamma^p_{\Z}(A) / \Ker(F) \simeq \mathcal A $ thus yields a canonical presentation of $\mathcal A$, depending only on $A$ and $\gamma_p$. We infer that $\Psi_2$ is fully faithful. \\
   It remains to be shown that it is essentially surjective. To prove this, we first reset notation. Let $(A,\gamma_p)$ be an arbitrary $2$-wrinkled ring. We denote by $$ \mathcal I \subset  \Gamma^p_{\Z}(A)$$ the subset consisting of elements $X$ which can be written as $$X=  [x_1]_p+ \ldots + [x_m]_p,$$ where $x_1, \ldots, x_m \in A$  are such that $$ x_1+  \ldots +x_m \in \Ker(\Frob_A)$$ and $$-\gamma_p(x_1 + \ldots+ x_m) = \sum_{a_1, \ldots, a_m} C_{a_1, \ldots, a_m} x_1^{a_1} \ldots x_m^{a_m} \in A, $$ where the sum ranges over all proper partitions of $p$. From the equality $$ \gamma_p(x+x')=\sum_{i+i'=p} \gamma_i(x)\gamma_{i'}(x'),$$ which holds for all $x, x' \in \Ker(\Frob_A)$, we see that $ \mathcal I$ is in fact a subgroup of $\Gamma^p_{\Z}(A)$. It is thus an ideal of $\mathcal I$. We put $$\mathcal A:= \Gamma^p_{\Z}(A) /\mathcal I.$$
   
   The surjection of rings $$f: \Gamma^p_{\Z}(A) \lra A, $$ $$ [x]_p \mapsto x^p,$$ clearly factors through $\mathcal I$, yielding a surjective ring homomorphism $$\pi: \mathcal A \lra A. $$ Pick an element $X \in \Ker (\pi).$ Write it as  $$X=  [x_1]_p+ \ldots + [x_m]_p,$$ where $x_1, \ldots, x_m \in A$. We have $$ x_1+  \ldots +x_m \in \Ker(\Frob_A).$$  Choose $x_{m+1} \in A$ such that $$x_{m+1}^p=\gamma_p(x_1+  \ldots +x_m) \in A.
  $$ Then one has $$  [x_1]_p+ \ldots + [x_m]_p +p  [x_{m+1}]_p \in \mathcal I$$ (verification left to the reader). We conclude that $X \in p \mathcal A $. Therefore, we have $$\Ker(\pi)=p\mathcal A.$$ Now, pick an element $Y \in \mathcal A[p].$ Write it as $$Y=  [y_1]_p+ \ldots + [y_m]_p,$$ where $y_1, \ldots, y_m \in A$. From the equality  $$p[y_1]_p+ \ldots + p[y_m]_p \in \mathcal I,$$ we see (by definition of $\mathcal I$) that  $$y_1^p+ \ldots+ y_m^p =0 \in A.$$ Hence $Y$ belongs to $\Ker(\pi)=p\mathcal A$. Altogether, we see that $$\mathcal A[p]=p\mathcal A,$$ i.e. that $\mathcal A$ is a free $(\Z /p^2 \Z)$-module. The ring $\mathcal A$ is thus a $2$-flat ring, with a canonical isomorphism $\mathcal A/p \simeq A.$ From the definition of $\mathcal I$, it is straightforward to check that $$\Psi(\mathcal A)=(A, \gamma_p), $$ which finishes the proof.

   \end{dem}
  
  \begin{rem}
  The descent statement of the previous Proposition has a clear  analogue in the classical context of  (quasiprojective) complex varieties, as follows.\\
  It is standard that a complex variety may be view as a real variety, of double dimension. This is a simple form of Weil's restriction of scalars. It is also standard that the data of an anti-involution on a  complex variety $Y$ is equivalent to giving  a   real variety $X$ and an isomorphism $Y \simeq X \times_\R \C$ of complex varieties. This is a simple (linear) descent statement.\\
      Let us now switch to the $p$-adic setting. Then Greenberg's functor allows  to consider a scheme over $\Z /p^2 \Z$ as a scheme over $\F_p$, of double dimension.  This is, in some sense, a non-linear analogue of Weil's scalar restriction, in which $\Z /p^2 \Z$ (resp. $\F_p$) plays the r\^ole of $\C$ (resp. $\R$).\\
      Proposition \ref{wrinklelift} is then an analogue of the simple descent statement above,  for the  morphism  $\Z /p^2 \Z \lra \Z /p \Z$. But roles are  exchanged in this second analogy:  $\Z /p^2 \Z$ now plays that  of $\R$, whereas $\F_p$ plays that of $\C$! In a daring poetic sense, Proposition \ref{wrinklelift} is both a descent result and a lifting result: it just depends in which direction you choose to look...
  \end{rem}

  We are grateful to Luc Illusie for his remarks, which led to the following improvement.
  
  \begin{prop}
Let $\mathcal A$ be a $2$-flat ring, corresponding to the $2$-wrinkled ring $(A,\gamma_p)$ (cf. Proposition \ref{wrinklelift}). \\Denote by $I \subset A$ the kernel of the Frobenius homomorphism.  Then the following conditions are equivalent.\\

i) The Frobenius of $A$ admits a (unique) lift to a (surjective) endomorphism of the ring $\mathcal A$.\\
ii) The  divided power operation $\gamma_p: I \lra A$ takes values in $I$.\\

If these conditions are fulfilled, then there exists a unique structure of $PD$-ideal on $I$, with $\gamma_p$ as $p$-th divided power operation.
\end{prop}

\begin{dem}
By the preceding Proposition, i)  holds if and only if the Frobenius of $A$ commutes with $\gamma_p$. This is obviously equivalent to $\gamma_p$ taking its values in $I$, qed.\\
The last assertion follows from the fact that a divided power structure on an ideal, in our context, is uniquely determined by the data of $\gamma_p$, see for instance Stacks Project, Tag 07H4, Lemma 23.5.3.
\end{dem}

  \begin{exo}
  Let  $A$ be a $2$-liftable ring. Put $$I:= \Ker(\Frob_A).$$ 
  
  1) Show that  the set of maps $$\gamma_p: I \lra A$$ such that $(A, \gamma_p)$ is a $2$-wrinkled ring is  a principal homogeneous space of $\Hom_{\Frob}(I/I^2, A)$. \\
 2) How is 1) connected  to Illusie's theory of the cotangent complex? \\
 
 Assume now that $I=I^2$,  and  that $(A,\gamma_p)$ is  a $2$-wrinkled ring.
 
 3) Show that  $\gamma_p=0.$ \\
 4) Deduce that $I^p=0$, hence that $A$ is a perfect  $\F_p$-algebra.
  \end{exo}

The generalisation of Proposition \ref{wrinklelift} to higher level descent statements is left to future considerations. 

\section*{ Acknowledgements}
We are  grateful to Patrick Brosnan for his support, and for spotting a mistake in  the last part of the  previous version of this paper, one year ago.  We owe thanks to Michel Brion for precious comments on the first (correct) part of the previous version. We thank J\'an Min\'a\v{c} for his kind and enthusiastic support. We thank Pierre Guillot and Fabien Morel for interesting discussions.


\begin{thebibliography}{2}
\bibitem[CTS]{CTS}\textsc{Colliot-Th\'el\`ene, J.-L.,  Sansuc, J.-J.}--- \textit{La $R$-\'equivalence sur les tores}, Ann. sci. \'ENS \textbf{10} (1977), no. 2,  175-229. \\

\bibitem[EM]{EM}\textsc{S. End\^o, T. Miyata}--- \textit{Integral representations with trivial first cohomology groups}, Nagoya Math. J. \textbf{85} (1982), 231–-240. \\

\bibitem[Fe]{Fe} \textsc{D. Ferrand}.--- \textit{Un foncteur norme}, Bull. Soc. Math. France   \textbf{126} (1998),  no. 1, 1-49. \\
\bibitem[FFSS]{FFSS} \textsc{V. Franjou, E. Friedlander, A. Scorichenko, A. Suslin}.--- \textit{General linear and functor cohomology over finite fields}, Ann. of Math. \textbf{150} (1999),  no. 2, 663-728. \\
\bibitem[GS]{GS} \textsc{P. Gille, T. Szamuely}.--- \textit{Central simple algebras and Galois cohomology},  Cambridge Studies in Advanced Mathematics  \textbf{101} (2006),  Cambridge University Press. \\
\bibitem[J]{J} \textsc{J. C. Jantzen}.--- \textit{Representations of algebraic groups. Second edition.},  Math. surveys and monographs  \textbf{107}, AMS (2003). \\
\bibitem[K1]{K} \textsc{D. Kaledin}.--- \textit{Witt vectors as a polynomial functor},   preprint, available on the arXiv server.\\
\bibitem[K2]{K2} \textsc{D. Kaledin}.--- \textit{Witt vectors, commutative and non-commutative},   preprint, available on the arXiv server.\\
\bibitem[Ro]{Ro} \textsc{N. Roby}.--- \textit{Lois polynomes et lois formelles en th\'eorie des modules}, Ann. Sci. \'ENS (3)  \textbf{80} (1963), 213-348. \\

\bibitem[St]{St} \textsc{G. S. St\r ahl}.--- \textit{An Intrinsic Definition of the Rees Algebra of a Module}, Proc. Edinburgh Math. Soc. (2), to appear. \\

\bibitem[S]{S} \textsc{P. Scholze}.--- \textit{Perfectoid spaces}, Pub. Math. IH\'ES  \textbf{116} (2012), 245-313.  \\
\bibitem[Se]{Se} \textsc{J.-P. Serre}.--- \textit{Galois cohomology}, Springer-Verlag, 2002. 
\bibitem[Sm]{Sm} \textsc{L. Smith}.--- \textit{An algebraic introduction to the Steenrod algebra}, Geo. \& Top. Monographs  \textbf{11} (2007), 327-348.  \\
\bibitem[Vo]{Vo} \textsc{V. Voevodsky}.--- \textit{On motivic cohomology with $\Z/l$-coefficients}, Ann. of Math.  \textbf{174} (2011), 401-438.\\
\bibitem[Vos]{Vos} \textsc{V.E. Voskresenskii}.--- \textit{Algebraic groups and their birational invariants}, Trans.  Math. Mono. AMS \textbf{179} (1998). 
\\



\end{thebibliography}
\end{document}